\documentclass[12pt]{article}

\pdfoutput=1

\usepackage[english]{babel}

\usepackage{amsfonts, amssymb, amsmath, amsthm, epsf,cite}
\usepackage{booktabs}

\usepackage[usenames,dvipsnames]{color}
\usepackage{color}
\usepackage{graphicx}
\usepackage{float}

\usepackage{fullpage}

\newtheorem{theorem}{Theorem}[section]
\newtheorem{observation}{Observation}[section]

\newtheorem{lemma}{Lemma}[section]

\newtheorem{corollary}{Corollary}[section]
\newtheorem{condition}{Condition}[section]

\newtheorem{approach}{Practical approach}[section]

\theoremstyle{definition}
\newtheorem{definition}{Definition}[section]

\newtheorem{remark}{Remark}[section]

\numberwithin{equation}{section} \numberwithin{figure}{section}
\numberwithin{figure}{section}
\numberwithin{table}{section}

\renewcommand{\phi}{{\varphi}}

\newcommand{\Vest}{V_{\rm est}}
\newcommand{\Var}{\rm Var}
\newcommand{\tVest}{\tilde{V}_{\rm est}}
\newcommand{\LerSp}{{\rm LerSp}}

\newcommand{\cN}{{\mathcal N}}

\newcommand{\cL}{{\mathcal L}}

\newcommand{\bC}{{\mathbf C}}

\newcommand{\bS}{{\mathbf S}}

\newcommand{\by}{{\mathbf y}}

\newcommand{\bz}{{\mathbf z}}

\newcommand{\btau}{{\boldsymbol{\tau}}}
\newcommand{\bmu}{{\boldsymbol{\mu}}}
\newcommand{\btheta}{{\boldsymbol{\theta}}}

\newcommand{\bbR}{{\mathbb R}}

\newcommand{\bbE}{{\mathbb E}}
\newcommand{\bbV}{{\mathbb V}}

\let\phi=\varphi

\providecommand{\keywords}[1]{\textbf{\textit{Keywords.}} #1}

\title{Gradient conjugate priors and multi-layer neural networks}

\author{Pavel  Gurevich\thanks{Free University of Berlin, Arnimallee 3, 14195 Berlin, Germany; RUDN University, Miklukho-Maklaya 6, 117198 Moscow, Russia; email: gurevich@math.fu-berlin.de},
Hannes Stuke\thanks{Free University of Berlin, Arnimallee 7, 14195 Berlin, Germany; email: h.stuke@fu-berlin.de}\thanks{Equal contribution.}}

\begin{document}

\maketitle

\begin{abstract}
The paper deals with learning probability distributions of observed data by artificial neural networks. We suggest a so-called gradient conjugate prior (GCP) update appropriate for neural networks, which is a modification of the classical Bayesian update for conjugate priors. We establish a connection between the gradient conjugate prior update and the maximization of the log-likelihood of the predictive distribution. Unlike for the Bayesian neural networks, we use deterministic weights of neural networks, but rather assume that the ground truth distribution is normal with unknown mean and variance and learn by the neural networks the parameters of a prior (normal-gamma distribution) for these unknown mean and variance. The update of the parameters is done, using the gradient that, at each step,  directs towards minimizing the Kullback--Leibler divergence from the prior to the posterior distribution (both being normal-gamma). We obtain a corresponding dynamical system for the prior's parameters and analyze its properties. In particular, we study the limiting behavior of all the prior's parameters and show how it differs from the case of the classical full Bayesian update. The results are validated on synthetic and real world data sets.
\end{abstract}


\tableofcontents

\section{Introduction}

Reconstructing  probability distributions of observed data by artificial neural networks is one of the most essential parts of machine learning and artificial intelligence~\cite{BishopBook,MurphyBook}. Learning  probability distributions not only allows {one} to predict the behavior of a system under consideration, but to also quantify the uncertainty with which the predictions are done. Under the assumption that the data are normally distributed, the most well studied way of reconstructing  probability distributions is the Bayesian learning of neural networks~\cite{MacKay1992}. One treats the weights of the network as normally distributed random variables, prescribes their prior distribution, and then finds the posterior distribution conditioned on the data. The main difficulty is that neither the posterior, nor the resulting predictive distributions are given in a closed form. As a result, different approximation methods have been developed~\cite{Neal95,Jylaenki14,Hinton93, Welling2011, Blundell2015, Kingma2015, GalGhahramani2015, HernandezLobato15, GalThesis16, LiTurner2016, HernandezLobatoLi2016, LiuWang2016, LiGal2017, Louizos2017}.
However, many of them have certain drawbacks related to the lack of scalability in data size or the neural network complexity, and are still a field of ongoing research. Furthermore, Bayesian neural networks often  assume homoscedastic variance in the likelihood (i.e., same for all samples) and rather learn uncertainty due to lack of data (epistemic uncertainty).
Among other methods for uncertainty quantification, there are the {delta method}~\cite{WildSeber,18Hwang, 19Veaux}, the {mean-variance estimate}~\cite{Nix}, and {deep ensemble methods}~\cite{Lakshminarayanan16,Lakshminarayanan17}.
A combination of the Bayesian approach (using the dropout variational inference) with the mean-variance estimate was used in~\cite{KendallGal17}, thus allowing for a simultaneous estimation of epistemic and aleatoric (due to noise in data) uncertainty. A new method based on minimizing a joint loss for a regression network and another network quantifying uncertainty was recently proposed in~\cite{GurHannesLU}. We refer to~\cite{KhossraveReview, Titterington04} and the recent works~\cite{Myshkov16, GalThesis16,Lakshminarayanan16,Lakshminarayanan17,GurHannesLU} for a comprehensive comparison of the above methods and further references to research on the Bayesian learning of neural networks.

We study an alternative approach to reconstructing the ground truth probability distribution based on what we call a {\em gradient conjugate prior (GCP) update}. We are interested in learning conditional probability distributions $q(y|x)$ of targets\footnote{Throughout this paper, we denote random variables by bold letters and the arguments of their probability distributions by the corresponding non-bold letters.} $\by\in\bbR$ corresponding to data samples $x\in\bbR^m$, using artificial neural networks ({\em supervised learning}). For brevity, we will often omit the dependence of distributions on $x$. Thus, assuming that the ground truth distribution $q(y)$ of a random variable $\by$ (corresponding to observed data) is Gaussian with unknown mean and precision, we let neural networks  learn the four parameters of the normal-gamma distribution that serves as a prior for the mean and variance of $\by$. We emphasize that, unlike for Bayesian neural networks, the {\em weights of the neural networks are deterministic} in our approach.  Given a parametrized prior,  one has the predictive distribution in the form of a (non-standardized) Student's t-distribution $p_{\rm pred}(y)$, whose parameters are explicitly determined by the outputs of the neural networks. For further details, we refer to Sec.~\ref{secMotivation}, which includes a graphical model visualization in Fig.~\ref{figGraphicalModel} and a comparison with Bayesian neural networks in Table~\ref{tableBayesVsGCP}.

Given an observation $(x_n,y_n)$ $(n=1,\dots,N)$, the classical Bayesian update yields the posterior distribution for the mean and variance of $\by_n$. This posterior appears normal-gamma as well~\cite{BishopBook}. However, {\em one cannot update its parameters directly} because they are represented by the outputs of the neural networks. Instead, one has to update the weights of the neural networks. We suggest to make a gradient descent step in the direction of minimization of the Kullback--Leibler (KL) divergence from the prior to the posterior (see the details in Sec.~\ref{secMotivation}). This is the step that we call the GCP update. After updating the weights, one takes the next observation $(x_{n+1},y_{n+1})$ and repeats the above update procedure. One cycles over the whole training data set until convergence of the log-likelihood of predictive distribution
\begin{equation}\label{eqLogLikelihoodPredictive}
\frac{1}{N}\sum\limits_{n=1}^N \ln p_{\rm pred}(y_n|x_n).
\end{equation}

In the paper, we provide a detailed analysis of the dynamics given by the GCP update. Intuitively, one might think that the GCP update, after convergence, yields the same result as the classical CP update. Surprisingly, this is not the case: the {\em parametrized normal-gamma distribution does not converge to the Bayesian posterior} (see Remark~\ref{remVestAlphaInfinity}). Nevertheless, the predictive distribution does converge to the ground truth Gaussian distribution $q(y)$. This is explained by the observation, which we prove in Sec.~\ref{subsecGSPandPredictiveDistribution}: the {\em GCP update is actually equivalent to maximizing by gradient ascent the log-likelihood~\eqref{eqLogLikelihoodPredictive} of the predictive distribution}. As the number of observations tends to infinity the {\em GCP update becomes also equivalent to minimizing by gradient descent the KL divergence from the predictive distribution $p_{\rm pred}(y)$ to the ground truth distribution $q(y)$}.  We show that these equivalences hold in general, even if the prior is not conjugate to the likelihood function. Thus, we see that the GCP method estimates aleatoric uncertainty.

We emphasize that, although in our approach the approximating distribution gets pa\-ra\-met\-ri\-zed (as the predictive distribution in the mean-variance approach~\cite{Nix} or the approximating latent variables distribution in variational autoencoders~\cite{KingmaWelling2014}), the way we parametrize and optimize and the way we interpret the result is different, as shown in Fig.~\ref{figGraphicalModel} and summarized in Table~\ref{tableBayesVsGCP}.

Now let us come back to our original assumption that $q(y)$ is a normal distribution and $p_{\rm pred}(y)$ is a Student's t-distribution. The latter appears to be overparametrized (by four parameters instead of three). We keep it overparametrized in order to compare the dynamics of the parameters under the classical CP update and under the GCP update. Reformulation of our results for Student's t-distribution parameterized in the standard way by three parameters will be  straightforward. There is a vast literature on the estimation of parameters of Student's t-distribution, see, e.g., the overview~\cite{Nadarajah08} and the references therein. Note that, in the context of neural networks, different samples correspond to different inputs of the network, and hence they belong to  different Student's t-distributions with {\it different unknown parameters}. Thus, the maximization of the likelihood of Student's t-distribution with respect to the weights of the networks is one of the most common methods. In~\cite{Uludag07}, the possibility of utilizing evolutionary algorithms for maximizing the likelihood was explored experimentally. Another natural way is to use the gradient ascent with respect to the weights of the network. As we said, the latter is equivalent to the usage of the GCP update. In the paper, we obtain a dynamical system for the prior's parameters that approximates the GCP update (as well as the gradient ascent for maximization of Student's t-distribution). We study the dynamics of the prior's
parameters in detail, in particular analyzing their convergence properties. Our approach is illustrated with synthetic data and validated on various real-world data sets in comparison with other methods for learning probability distributions based on neural networks. To our best knowledge, {\em neither the dynamical systems analysis of the GCP (or gradient ascent for maximizing the likelihood of Student's t-distribution), nor a thorough comparison of the GCP with other methods has been carried out before}.

As an interesting and useful consequence of our analysis, we will see how the GCP interacts with the output outliers in the training set (a small percentage of observations that do not come from the assumed normal distribution $q(y)$). The outliers prevent one of the prior's parameters ($\alpha$, which is related to the number of degrees of freedom of $p_{\rm pred}(y)$) from going to infinity. On one hand, this is known~\cite{Lucas96,Scheffler08} to allow for a  better estimate of the mean and variance of $q(y)$, compared with directly using the maximization of the likelihood of a normal distribution. However, on the other hand, this still leads to overestimation of the variance of $q(y)$. To deal with this issue, we obtain an explicit formula (see~\eqref{eqEstMeanVarNetwork}) that allows one to correct the estimate of the variance and recover the ground truth variance of $q(y)$. To our knowledge, {\em such a correction formula was not derived in the literature before.}

%
%

The paper is organized as follows. In Sec.~\ref{secMotivation}, we provide a detailed motivation for the GCP update, explain how we approximate the parameters of the prior distribution by neural networks, establish the relation between the GCP update and the predictive distribution, and formulate the method of learning the ground truth distribution from the practical point of view. Section~\ref{secDynamics4} is the mathematical core of this paper. We derive a dynamical system for the prior's parameters, induced by the GCP update, and analyze it in detail. In particular, we obtain an asymptotics for the growth rate of $\alpha$ and find the limits of the other parameters of the prior. In Sec.~\ref{secDynamicsFixedAlpha}, we study the dynamics for a fixed $\alpha$. We find the limiting values for the rest of the parameters and show how one can recover the variance of the ground truth normal distribution $q(y)$. In Sec.~\ref{secRoleFixedAlpha}, we clarify the role of a fixed $\alpha$. Namely, we compare the sensitivity to output outliers of the GCP update with that in minimizing the standard squared error loss or maximizing the log-likelihood of a normal distribution. Furthermore, we show how $\alpha$ controls the learning speed in clean and noisy regions. In Sec.~\ref{secNN}, we illustrate the fit of neural networks for synthetic and various real-world data sets. Section~\ref{secConclusion} contains a conclusion and an outline of possible directions of further research. Appendices~\ref{appendix}--\ref{appendixLemmaDKalpha} contain the proofs of auxiliary lemmas from Sec.~\ref{secDynamics4}. In Appendix~\ref{appendixHyperparameters}, we present the values of hyperparameters of different methods that are compared in Sec.~\ref{secNN}.

\section{Motivation}\label{secMotivation}

\subsection{Estimating normal distributions with unknown mean and precision}

Assume one wants to estimate unknown mean and precision (the inverse of the variance) of normally distributed data $\by|x$. We remind that $\by$ is conditioned on $x\in \bbR^m$, but we often omit this dependence in our notation. We will analyze scalar $\by$ and refer to Sec.~\ref{secConclusion} for a discussion of multivariate data.  One standard approach for estimating the mean and precision is based on conjugate priors. One assumes that the mean and precision are random variables, $\bmu$ and $\btau$ respectively, with a joint prior given by the normal-gamma distribution
\begin{equation}\label{eqNormalGamma}
p(\mu,\tau|m,\nu,\alpha,\beta) = \dfrac{\beta^\alpha \nu^{1/2}}{\Gamma(\alpha)(2\pi)^{1/2}} \tau^{\alpha-1/2} e^{-\beta\tau} e^{-\frac{\nu\tau(\mu-m)^2}{2}},
\end{equation}
where
$
m\in \bbR,\ \nu>0,\ \alpha>1,\ \beta>0.
$

The marginal distribution for $\bmu$  is a non-standardized Student's t-distribution with
\begin{equation}\label{eqEVmu}
\bbE[\bmu] = m,\quad \bbV[\bmu] = \dfrac{\beta}{\nu(\alpha-1)}.
\end{equation}
The marginal distribution for $\btau$ is the Gamma distribution with
\begin{equation}\label{eqEVtau}
\bbE[\btau] = \dfrac{\alpha}{\beta}, \quad \bbV[\btau] = \dfrac{\bbE(\tau)}{\beta} = \dfrac{\alpha}{\beta^2}.
\end{equation}

By marginalizing $\bmu$ and $\btau$, one obtains the predictive distribution $p_{\rm pred}(y) = p_{\rm pred}(y;m,\nu,\alpha,\beta)$ for $\by$, which appears to be a non-standardized Student's t-distribution. Its mean and variance can be used to estimate the mean and variance of $\by$. The estimated mean $m_{\rm est}$ and variance $\Vest$ are given by
\begin{equation}\label{eqStandardCPEstimate}
m_{\rm est}:=m,\quad \Vest := \dfrac{\beta(\nu+1)}{(\alpha-1)\nu}.
\end{equation}
We refer, e.g., to~\cite{BishopBook} for further details.

\begin{figure}[t]
    \center
	   \includegraphics[width=0.5\textwidth]{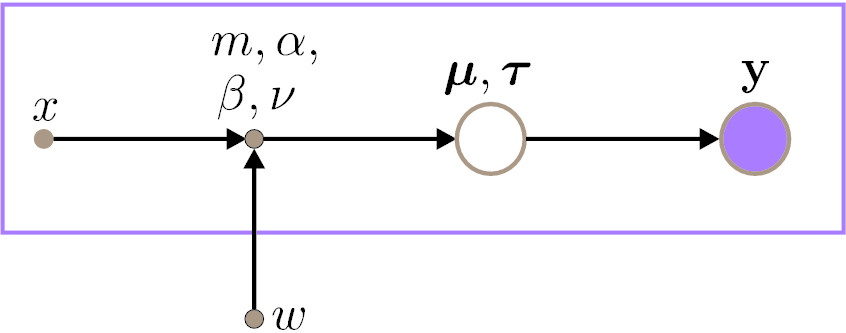}
    \caption{Deterministic parameters (input $x$, the prior's parameters $m,\alpha,\beta,\nu$, and the weights $w$) are shown by solid nodes and random variables ($\bmu,\btau$, and $\by$) by circles. The shaded circle corresponds to the observed random variable $\by$. The box encompasses the quantities depending on $x$.}\label{figGraphicalModel}
\end{figure}

Our goal is to approximate the parameters $m,\alpha,\beta,\nu$ by multi-layer neural networks, i.e., to represent them as functions of inputs and weights: $m=m(x,w_1), \alpha=\alpha(x,w_2), \beta=\beta(x,w_3), \nu=\nu(x,w_4)$,  $w_j\in\bbR^{M_j}$.  The corresponding graphical model is shown in Fig.~\ref{figGraphicalModel}.

A standard way to estimate the parameters $w$ of the neural networks is to maximize the log-likelihood
\begin{equation}\label{eqMaxLogLikelihoodPred}
\frac{1}{N}\sum\limits_{n=1}^N p_{\rm pred}(y_n;m(x_n,w_1), \alpha(x_n,w_2), \beta(x_n,w_3), \nu(x_n,w_4))\to\max\limits_{w_1,\dots,w_4}.
\end{equation}
An alternative approach is to use the conjugacy of the prior, which allows for an explicit formula for the Kullback--Leibler (KL) divergence from the posterior to the prior. This approach is explained in Sec.~\ref{secCPUpdate}. In Sec.~\ref{subsecGSPandPredictiveDistribution}, we show that both methods are actually {\em equivalent}.

\subsection{Conjugate prior update}\label{secCPUpdate}

Suppose one observes a new sample $ y $. Then, by the Bayes theorem, the conditional distribution of $(\bmu,\btau)$ under the condition that $\by= y $ (the posterior distribution denoted by $p_{\rm post}(\mu,\tau)$) appears to be normal-gamma as well~\cite{BishopBook}, namely,
\begin{equation}\label{eqPosterior}
  p_{\rm post}(\mu,\tau) = p(\mu,\tau|m',\nu',\alpha',\beta', y ),
\end{equation}
where $p(\cdot)$ is defined in~\eqref{eqNormalGamma} and the parameters are updated as follows:
\begin{equation}\label{eqCPUpdade}
\begin{gathered}
 m' = \dfrac{\nu m +  y }{\nu + 1},\quad
\nu' = \nu+1,\quad \alpha' = \alpha+\frac{1}{2},\quad
\beta' = \beta +   \frac{\nu}{\nu+1}\frac{( y - m)^2}{2}.
\end{gathered}
\end{equation}
We call~\eqref{eqCPUpdade} the {\em conjugate prior (CP) update}.

The KL divergence from a continuous distribution $p_{\rm post}$ to a continuous distribution $p$ is defined as follows:
\begin{equation}\label{eqKLDivergenceIntegral}
 {\displaystyle D_{\mathrm {KL} }(p_{\rm post}\|p):=\int p_{\rm post}(\mu,\tau)\,\ln {\frac {p_{\rm post}(\mu,\tau)}{p(\mu,\tau)}}\,d\mu\, d\tau}.
\end{equation}
We denote by $
 \Psi (x):= {\Gamma '(x)}/{\Gamma (x)}
$ the digamma function,
where $\Gamma(x)$ is the gamma function. Then for the above normal-gamma distributions~\eqref{eqNormalGamma} and~\eqref{eqPosterior} the KL divergence takes the form~\cite{Soch16}
\begin{equation}\label{eqKLDivergence}
\begin{aligned}
K( m,\nu,\alpha,\beta) & := \dfrac{1}{2}\frac{\alpha'}{\beta'}(m-m')^2\nu + \frac{1}{2}\frac{\nu}{\nu'} - \frac{1}{2}\ln\frac{\nu}{\nu'}-\frac{1}{2}
\\
& - \alpha\ln\frac{\beta}{\beta'} + \ln\frac{\Gamma(\alpha)}{\Gamma(\alpha')} - (\alpha-\alpha')\Psi(\alpha') + (\beta-\beta')\frac{\alpha'}{\beta'}.
\end{aligned}
\end{equation}

In our case, one cannot directly apply the update in~\eqref{eqCPUpdade}, but one must update the weights $w_j$ instead. The natural way to do so is to observe a sample~$ y $, to calculate the posterior distribution~\eqref{eqPosterior} and to change the weights $w$ in the direction of $-\nabla_w K$, i.e.,
\begin{equation}\label{eqGradientDescentW}
w_{\rm new} := w - \lambda \nabla_w K,\quad w\in\{w_1,\dots,w_4\},
\end{equation}
where $\lambda>0$ is a learning rate. When we compute the gradient of $K$ with respect of $w$, we keep all the prime variables in~\eqref{eqKLDivergence} fixed and do not treat them as functions of $w$, while all the nonprime variables are treated as functions of $w$. We still use the notation $\nabla_w K$ in this case. We call~\eqref{eqGradientDescentW} the {\em gradient conjugate prior (GCP) update}. In Table~\ref{tableBayesVsGCP}, we summarize our approach and highlight its difference from the Bayesian neural networks and variational inference\footnote{The latent variables are usually denoted $w$ in the Bayesian neural networks framework or $\bz$ in the variational inference framework. We use the notation $\btheta$ to make it consistent with our notation in Sec.~\ref{subsecGSPandPredictiveDistribution}.}.

As we will see below, this update induces the update for $(m,\alpha,\beta,\nu)$ that is different from the classical conjugate prior update~\eqref{eqCPUpdade} and yields a completely different dynamics.
Before we analyze this dynamics in detail, we explain why the GCP update~\eqref{eqGradientDescentW} is equivalent to maximizing the log-likelihood~\eqref{eqMaxLogLikelihoodPred} of the predictive in the general case.


\begin{table}[t]
  \centering
\resizebox{\textwidth}{!}{%
    \begin{tabular}{l|l|l}

    {} & {\bf Bayesian neural networks} & {\bf GCP networks} \\

    \hline
    \hline

    {\bf Data} & \multicolumn{2}{c}{Inputs $X=\{x_1,\dots,x_N\}$, targets $Y=\{y_1,\dots,y_N\}$}\\

    \hline

    {\bf Ground truth} & \multicolumn{2}{c}{Gaussian $q(y|x)$} \\

    \hline

    \begin{tabular}{@{}l} {\bf Captured} \\ {\bf uncertainty} \end{tabular} & {\em Epistemic, homoscedastic}  &
     {\em Aleatoric, heteroscedastic}  \\

    \hline

    {\bf Weights} & {\em Random} $\btheta$  &
     {\em Deterministic}  $w$  \\

    \hline

    \begin{tabular}{@{}l} {\bf Latent} \\ {\bf variables} \end{tabular} & Weights $\btheta\sim p(\theta)$   {\em independent} of $x\in X$  &
    \begin{tabular}{@{}l} Means and precisions $(\bmu,\btau)=\btheta\sim p(\theta|x,w)$ \\
    {\em conditioned} on $x\in X$ and $w$ \end{tabular} \\
    \hline

    {\bf Prior} & $p(\theta)$ {\em fixed} during training & $p(\theta|x,w)$ {\em evolves} during training. \\

    \hline

    {\bf Likelihood} & Gaussian $\cN(y|m(x,\theta),\tau^{-1})$ with {\em constant}  $\tau$ & \begin{tabular}{@{}l} Gaussian $\cN(y|\mu,\tau^{-1})$ with both $\mu$ and  $\tau$ \\ {\em depending on} $x\in X$ and $w$. \end{tabular}\\

    \hline

{\bf Posterior} & $p(\theta|X,Y)$ {\em intractable} and {\em fixed} during training & $p_{\rm post}(\theta)$ {\em tractable} normal-Gamma and {\em evolves} during training \\

    \hline

    {\bf Training} & \begin{tabular}{@{}l}
    Minimize  the   KL-divergence $D_{\rm KL}(p(\theta|w)\|p(\theta|X,Y))$ \\
    from $p(\theta|w)$ parametrized by deterministic $w$ \\
    to the {\em intractable fixed} posterior $p(\theta|X,Y)$ \end{tabular}
      & \begin{tabular}{@{}l} Gradient descent step to minimize the {\em reverse} \\
      KL-divergence $D_{\rm KL}(p_{\rm post}(\theta)\|p(\theta|x,w))$ w.r.t. $w$. \\
      The posterior $p_{\rm post}(\theta)$ is {\em  tractable} and {\em recalculated}\\ after each observation $(x,y)$ based on the evolving prior $p(\theta|x,w)$     \end{tabular}   \\

    \hline

    {\bf Result} & $p(\theta|w)$ {\em approximates} the posterior $p(\theta|X,Y)$  & \begin{tabular}{@{}l} $p(\theta|x,w)$ does {\em not} converge to $p_{\rm post}(\theta)$, \\ but the predictive distribution maximizes the likelihood of data \end{tabular}. \\

    \hline

    {\bf Predictive} & \begin{tabular}{@{}l} $p_{\rm pred}(y|x,w)=\int\cN(y|m(x,\theta),\tau^{-1})p(\theta|w)\,d\theta$ \\ typically evaluated by {\em sampling} \end{tabular} & \begin{tabular}{@{}l} {\em Explicit} Student's~t-distribution  $p_{\rm pred}(y|x,w) = t_{2\alpha}(y|m,\beta(\nu+1)/(\nu\alpha))$, \\ where $m,\alpha,\beta,\nu$ depend on $x$ and $w$  \end{tabular} \\

    \hline

     \begin{tabular}{@{}l} {\bf Output outliers}\\
     {\bf in the training set} \end{tabular} & \begin{tabular}{@{}l} {\em Distorted} means and {\em overestimated} variances \end{tabular} & {\em Robust} means and variances via the correction formula~\eqref{eqEstMeanVarNetwork} \\

  \end{tabular}
  }
  \caption{Comparison of Bayesian neural networks with variational inference and GCP networks}\label{tableBayesVsGCP}
\end{table}

\subsection{Maximization of the predictive distribution and the GCP update}\label{subsecGSPandPredictiveDistribution}

Suppose we want to learn a ground truth probability distribution $q(y)$ of a random variable~$\by$ (a normal distribution in our particular case). Since the ground truth distribution is a priori unknown, we conjecture that it belongs to a family of distributions $L(y|\theta)$ parametrized by~$\theta$ (in our case $\theta=(\mu,\tau)$ and $L(y|\theta)$ is a normal distribution with mean $\mu$ and precision~$\tau$). Since $\theta$ is a priori unknown, we assume it is a random variable $\btheta$ with a prior distribution from a family $p(\theta|w)$ parametrized by $w$ (in our case, $p(\theta|w)$ is the normal-gamma distribution and $w$ are the weights of neural networks approximating $m,\alpha,\beta,\nu$). We denote the predictive distribution by
\begin{equation}\label{eqPpred}
p_{\rm pred}(y|w) := \int L(y|\theta)p(\theta|w)\,d\theta
\end{equation}
(non-standardized Student's t-distribution in our case).
Given an observation $\by=y$, the Bayes rule determines the posterior distribution of $\btheta$:
\begin{equation}\label{eqPpost}
p_{\rm post}(\theta|w,y):=\frac{L(y|\theta) p(\theta|w)}{p_{\rm pred}(y|w)}.
\end{equation}
In our case, $p_{\rm post}(\theta|w,y)$ is normal-gamma again, but we emphasize that, in general, it need not be from the same family as the prior $p(\theta|w)$ is.

Now we compute the gradient of the KL divergence
\begin{equation}\label{eqKLGeneral}
K(w,y):= {\displaystyle D_{\mathrm {KL} }(p_{\rm post}\|p) =\int p_{\rm post}(\theta|w,y)\,\ln {\frac {p_{\rm post}(\theta|w,y)}{p(\theta|w)}}\ d\theta}
\end{equation}
(cf.~\eqref{eqKLDivergenceIntegral}) with respect to $w$, assuming that $w$ in the posterior distribution is {\em freezed}, and we do not differentiate it. Denoting such a gradient by $\nabla_w K(w,y)$, we obtain the following lemma.

\begin{lemma}\label{lNablawKy}
  $
\nabla_w K(w,y) =-\nabla_w \ln p_{\rm pred}(y|w).
  $
\end{lemma}
\proof
Freezing $p_{\rm post}(\theta|w,y)$ in~\eqref{eqKLGeneral}, we have
$$
\nabla_w K(w,y) = - \int \frac{p_{\rm post}(\theta|w,y) \nabla_w p(\theta|w)}{p(\theta|w)}\ d\theta.
$$
Plugging in $p_{\rm post}(\theta|w,y)$ from~\eqref{eqPpost} and using~\eqref{eqPpred} yields
$$
\nabla_w K(w,y) = - \int \frac{L(y|\theta) \nabla_w p(\theta|w)}{p_{\rm pred}(y|w)}\ d\theta = -\frac{\nabla_w p_{\rm pred}(y|w)}{p_{\rm pred}(y|w)}
=-\nabla_w \ln p_{\rm pred}(y|w).
$$
\endproof

Lemma~\ref{lNablawKy} shows that the GCP update~\eqref{eqGradientDescentW} is the gradient ascent step in the direction of maximizing the log-likelihood of the predictive distribution $p_{\rm pred}(y|w)$ given a new observation $\by=y$. Furthermore, using Lemma~\ref{lNablawKy}, we see that given observations $y_1,\dots,y_N$,
the averaged  GCP update of the parameters $w$ is given by (cf.~\eqref{eqGradientDescentW})
\begin{equation}\label{eqGradientDescentWAveragedFinite}
w_{\rm new} := w - \lambda\frac{1}{N}\sum_{n=1}^{N}[\nabla_w K(w,y_n)] = w +\lambda \nabla_w \left(\frac{1}{N}\sum_{n=1}^{N}
\ln p_{\rm pred}(y_n|w)\right).
\end{equation}

Further, if the observations are sampled from the ground truth distribution $q(y)$ and their number tends to infinity, then the GCP update~\eqref{eqGradientDescentWAveragedFinite} assumes the form
\begin{equation}\label{eqGradientDescentWAveraged}
\begin{aligned}
w_{\rm new} &:= w - \lambda \cdot\bbE_{\by \sim q(y)}[\nabla_w K(w,\by)]
 =
w-\lambda \nabla_w \int \ln p_{\rm pred}(y|w) q(y)\,dy \\
&= w- \lambda \nabla_w \int q(y) \ln\frac{q(y)}{p_{\rm pred}(y|w)} \,dy = w-\lambda\nabla_w D_{\rm KL}(q \| p_{\rm pred}(\cdot|w)).
\end{aligned}
\end{equation}

\begin{remark}\label{rGCPupdatePredictiveDistr}
\begin{enumerate}

 \item\label{rGCPupdatePredictiveDistr0} Formula~\eqref{eqGradientDescentWAveragedFinite} shows that the GCP update  {\em maximizes the likelihood of the predictive distribution $p_{\rm pred}(y|w)$ for the observations $y_1,\dots,y_N$}.

  \item\label{rGCPupdatePredictiveDistr1} Formula~\eqref{eqGradientDescentWAveraged} shows that the GCP update is equivalent to the gradient descent step for the {\em minimization of the KL divergence from the ground truth distribution $q(y)$ to the predictive distribution $p_{\rm pred}(y|w)$}. If the ground truth distribution $q(y)$ belongs to the family $p_{\rm pred}(y|w)$, then the minimum equals zero and is achieved for some (not necessarily unique) $w_*$ such that $p_{\rm pred}(y;w_*)=q(y)$; otherwise the minimum is positive and provides the best possible approximation of the ground truth in the sense of the KL divergence.

  \item\label{rGCPupdatePredictiveDistr2} In our case, $q(y)$ is a normal distribution and $p_{\rm pred}(y|w)$ are Student's t-distributions. In accordance with item~\ref{rGCPupdatePredictiveDistr1}, we will see below that the GCP update forces the number of degrees of freedom of $p_{\rm pred}(y|w)$ to tend to infinity. However, due to the overparametrization of the predictive distribution (four parameters $m,\alpha,\beta,\nu$ instead of three), the learned variance of $q(y)$ will be represented by a curve in the space $(\beta,\nu)$. The limit point $\beta_*,\nu_*$ to which $\beta,\nu$ will converge during the GCP update, will depend on the initial condition. Interestingly, $\beta_*,\nu_*$ will always be different from the limit  obtained by the classical CP update~\eqref{eqCPUpdade} (cf. Remark~\ref{remVestAlphaInfinity}).
\end{enumerate}
\end{remark}

\subsection{Practical approaches}

Based on Remark~\ref{rGCPupdatePredictiveDistr} (items~\ref{rGCPupdatePredictiveDistr0} and \ref{rGCPupdatePredictiveDistr1}), we suggest the following general practical approach.

\begin{approach}\label{apprGeneral}
\begin{enumerate}
  \item One approximates the parameters of the prior by neural networks:
  \begin{equation}\label{eqNetworksMAlphaBetaNu}
  m=m(x,w_1),\ \alpha=\alpha(x,w_2),\ \beta=\beta(x,w_3),\ \nu=\nu(x,w_4).
  \end{equation}
  We call them the {\em GCP neural networks}.
  \item One trains these four networks by the GCP update~\eqref{eqGradientDescentW} until convergence of  $m,\alpha,\beta,\nu$.
  \item The resulting predictive distribution is the non-standardized Student's t-distribution $t_{2\alpha}(y|m,\beta(\nu+1)/(\nu\alpha))$. The estimated mean $m_{\rm est}$ and variance $\Vest$ (for $\alpha>1$) are given by
  \begin{equation}\label{eqEstMeanVarNetwork'}
    m_{\rm est}:=m,\quad \Vest := \dfrac{\beta(\nu+1)}{(\alpha - 1)\nu}.
  \end{equation}
  
  \item Student's t variance $\Vest$ {\em overestimates} the ground truth variance of the normal distribution $q(y)$. However, one can still recover the correct variance of $q(y)$ by appropriately correcting $\Vest$. We show that the correction is given by
  \begin{equation}\label{eqEstMeanVarNetwork}
     \tVest := \dfrac{\beta(\nu+1)}{(\alpha-A(\alpha))\nu}
  \end{equation}
with $A(\alpha_*)$ from Definition~\ref{defA}. We call it a {\em correction formula} for the variance. The interplay between $\Vest$ and $\tVest$ is illustrated in sections~\ref{subsecSensitivityOutliers}, \ref{subsecSynthetic}, and~\ref{subsecRealWorldOutliers}. Our experiments show that this correction allows for reconstructing the ground truth variance even in the presence of outliers in the training set. 
\end{enumerate}
\end{approach}

In the rest of the paper, we rigorously justify the above approach,  study the dynamics of $m,\alpha,\beta,\nu$ under this update, and analyze how one should correct the variance for a fixed~$\alpha$.

\section{Dynamics of $m,\alpha,\beta,\nu$}\label{secDynamics4}

\subsection{Dynamical system  for $m,\alpha,\beta,\nu$}

The GCP update~\eqref{eqGradientDescentW}  induces the update for $(m,\alpha,\beta,\nu)$ as follows:
\begin{equation}\label{eqGradientDescentBeta}
m_{\rm new} := m(w-\lambda\cdot\nabla_w K)= m\left(w-\lambda\frac{\partial K}{\partial m}\nabla_w m(w)\right) \approx m(w) - \lambda(\nabla_w m)^T(\nabla_w m)\frac{\partial K}{\partial m},
\end{equation}
where $w=w_1$, and similarly for $\alpha,\beta,\nu$ and $w_2,w_3,w_4$, respectively.

Obviously, the new parameters $m_{\rm new},\alpha_{\rm new},\beta_{\rm new},\nu_{\rm new}$ are different from $m',\alpha',\beta',\nu'$ given by the classical conjugate prior update~\eqref{eqCPUpdade}. From now on, we replace $\lambda(\nabla_w m)^T(\nabla_w m)$, etc. by  new learning rates and analyze {\it how the parameters will change and to which values they will converge under the updates of the form
\begin{equation}\label{eqGradientDescent}
m_{\rm new} := m - \lambda_1\frac{\partial K}{\partial m}, \quad \alpha_{\rm new} := \alpha - \lambda_2\frac{\partial K}{\partial \alpha}, \quad
\beta_{\rm new} := \beta - \lambda_3\frac{\partial K}{\partial \beta}, \quad \nu_{\rm new} := \nu - \lambda_4\frac{\partial K}{\partial \nu},
\end{equation}
where $\lambda_j>0$ are the learning rates}.
As before, when we compute the derivatives of $K$, we keep all the prime-variables in~\eqref{eqKLDivergence} fixed and do not treat them as functions of $m,\nu,\alpha,\beta$. In other words, {\em we first compute the derivatives of $K$ with respect to $m,\nu,\alpha,\beta$ and then substitute $m',\nu',\alpha',\beta'$ from~\eqref{eqCPUpdade}}. For brevity, we will simply write $\partial K/\partial m$, etc. We call~\eqref{eqGradientDescent} the {\em GCP update} as well.

Setting
\begin{equation}\label{eqSigmaBetaNu}
\sigma:= \frac{\beta(\nu+1)}{\nu},
\end{equation}
we have
\begin{align}
\frac{\partial K}{\partial m} &=
\frac{\alpha+1/2}{\sigma  +   \frac{(m- y )^2}{2}}(m- y ),\label{eqDKmPure}\\
{}
\frac{\partial K}{\partial \alpha} &= \ln\left(1+ \frac{(m- y )^2}{2\sigma}\right) + \Psi(\alpha)-\Psi\left(\alpha+\frac{1}{2}\right),\label{eqDKalphaPure}\\
{}
\frac{\partial K}{\partial \beta} & = \frac{1}{\beta}\left(\frac{\alpha+1/2}{1+\frac{(m- y )^2}{2\sigma}} - \alpha\right),\label{eqDKbetaPure}\\
{}
\frac{\partial K}{\partial \nu} &=
\frac{1}{2\nu(\nu+1)}\left(\frac{\alpha+1/2}{\sigma +\frac{(m- y )^2}{2}}(m- y )^2 - 1\right).\label{eqDKnuPure}
\end{align}

In this section and in the next one, we will treat the parameters $m,\alpha,\beta,\nu$ as functions of time $t>0$ and study a dynamical system that approximates the GCP update~\eqref{eqGradientDescent} when the number of observations is large.
We will concentrate on the prototype situation, where all new learning rates are the same.
\begin{condition}\label{condLambdaEqual}
In the GCP update~\eqref{eqGradientDescent}, we have
  $\lambda_1=\lambda_2=\lambda_3=\lambda_4$.
\end{condition}

\begin{figure}[t]
    \center
	   \includegraphics[width=0.3\textwidth]{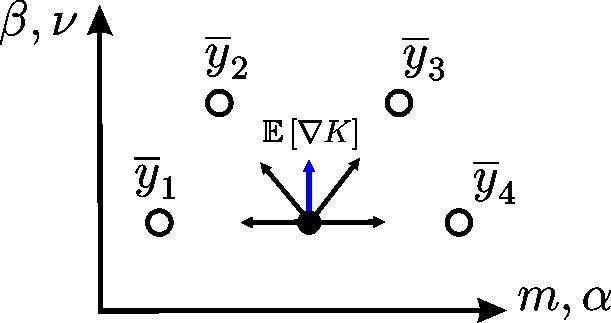}
    \caption{The black circle indicates the prior probability distribution~\eqref{eqNormalGamma} in the space of the parameters $m,\alpha,\beta,\nu$. The white circles indicate the posterior probability distributions~\eqref{eqPosterior} corresponding to different observations $y_1,y_2,\dots$. The black vectors are the gradients with respect to the nonprime variables of the corresponding KL divergences. The blue vector is the averaged gradient. An equilibrium of system~\eqref{eqODE4} would correspond to the case where the blue vector vanishes. Theorem~\ref{thStrip} shows that this actually never happens. However, Theorem~\ref{thDynamicsFixedAlpha} shows that if one keeps $\alpha$ fixed, but updates $m,\beta,\nu$, then one obtains a whole curve of equilibria.}\label{figKLDivergence}
\end{figure}

Under Condition~\ref{condLambdaEqual}, the approximating dynamical system takes the form
\begin{equation}\label{eqODE4}
\dot m = -\bbE\left[\frac{\partial K}{\partial m}\right],\quad \dot\alpha = -\bbE\left[\frac{\partial K}{\partial \alpha}\right],\quad
\dot\beta = -\bbE\left[\frac{\partial K}{\partial \beta}\right],\quad \dot\nu = -\bbE\left[\frac{\partial K}{\partial \nu}\right];
\end{equation}
hereinafter the expectations are taken with respect to the true distribution $q(y)$ of $\by$ which is treated as a normally distributed random variable with mean $\bbE[\by]$ and variance $V:=\bbV[\by]$, see Fig.~\ref{figKLDivergence}.

\begin{remark}\label{rGradientFlow4Eq}
  Due to~\eqref{eqGradientDescentWAveraged}, system~\eqref{eqODE4} defines a gradient flow with the potential  $D_{\rm KL}(q \| p_{\rm pred}(\cdot|w))$, where $p_{\rm pred}(y|w))$ is the  Student's t-distribution $t_{2\alpha}(y|m,\beta(\nu+1)/(\nu\alpha))$.
\end{remark}

\begin{remark}
  If Condition~\ref{condLambdaEqual} does not hold, then the respective factors $\lambda_j$ will appear in the right-hand sides in~\eqref{eqODE4}. The modifications one has to make in the arguments below are straightforward.
\end{remark}

\subsection{Estimation of the mean $m$}

Using~\eqref{eqDKmPure}, we obtain the formula for the expectation
\begin{equation}\label{eqDKmu}
\bbE\left[\frac{\partial K}{\partial m}\right] =\frac{\alpha+1/2}{(2\pi V)^{1/2}} \int\limits_{-\infty}^\infty \frac{m-y}{\sigma  + \frac{(m-y)^2}{2}} e^{-\frac{(\bbE[\by]-y)^2}{2V}}  dy.
\end{equation}

\begin{theorem}\label{thDKm}
  The first equation in~\eqref{eqODE4} has a unique equilibrium $m=\bbE[\by]$. It is stable in the sense that, for any $\alpha,\beta,\nu$, we have
  $$
  \dot m < 0\ \text{if } m>\bbE[\by],\quad \dot m > 0\ \text{if } m<\bbE[\by].
  $$
\end{theorem}
\proof
Without loss of generality, assume that $\bbE[\by]=0$ and $V=1$ (otherwise, make a change of variables $z=(y-\bbE[\by])/V^{1/2}$ in the integral in~\eqref{eqDKmu}). Then we obtain from~\eqref{eqDKmu}
\begin{equation}\label{eqDKmu1}
\bbE\left[\frac{\partial K}{\partial m}\right] =C_1  \int\limits_{-\infty}^\infty \frac{m-y}{C_2 + (m-y)^2} e^{-\frac{y^2}{2}}  dy
= - C_1  \int\limits_{-\infty}^\infty \frac{z}{C_2 + z^2} e^{-\frac{(m+z)^2}{2}}  dz,
\end{equation}
where $C_1,C_2>0$ do not depend on $m$.

Obviously, the right-hand side in~\eqref{eqDKmu1} vanishes at $m=0$. Furthermore, due to~\eqref{eqDKmu1}, for $m>0$,
$$
\bbE\left[\frac{\partial K}{\partial m}\right] = C_1 \int\limits_{0}^\infty \frac{z}{C_2 + z^2} \left(e^{-\frac{(m-z)^2}{2}}-e^{-\frac{(m+z)^2}{2}}\right)  dz > 0
$$
because $-(m-z)^2>-(m+z)^2$ for $m,z>0$. Similarly, $\bbE\left[\frac{\partial K}{\partial m}\right]<0$ for $m<0$.
\endproof

\subsection{Estimation of the variance. The unbounded absorbing set}

From now on, taking into account Theorem~\ref{thDKm}, we assume the following.
\begin{condition}\label{condmEqualsE}
  $m=\bbE[\by]$.
\end{condition}
Under Condition~\ref{condmEqualsE}, we study the other three equations in~\eqref{eqODE4}, namely,
\begin{equation}\label{eqODE3}
\dot\alpha = -\bbE\left[\frac{\partial K}{\partial \alpha}\right],\quad
\dot\beta = -\bbE\left[\frac{\partial K}{\partial \beta}\right],\quad \dot\nu = -\bbE\left[\frac{\partial K}{\partial \nu}\right],
\end{equation}
where (due to Condition~\ref{condmEqualsE})
\begin{align}
\bbE\left[\frac{\partial K}{\partial \alpha}\right] &= \frac{1}{(2\pi)^{1/2}}\int_{-\infty}^{\infty} \ln\left(1+
\frac{V}{\sigma}\frac{z^2}{2}\right)e^{-\frac{z^2}{2}}dz + \Psi(\alpha)-\Psi\left(\alpha+\frac{1}{2}\right),\label{eqDKalphaExp}\\
{}
\bbE\left[\frac{\partial K}{\partial \beta}\right] &=
\frac{1}{\beta}\left(\frac{\alpha+1/2}{(2\pi)^{1/2}}\int_{-\infty}^{\infty} \frac{1}{1 + \frac{V}{\sigma}\frac{z^2}{2}} e^{-\frac{z^2}{2}}dz  - \alpha\right),\label{eqDKbetaExp} \\
{}
\bbE\left[\frac{\partial K}{\partial \nu}\right] & =
\frac{1}{2\nu(\nu+1)}\left(\frac{2\alpha+1}{(2\pi)^{1/2}}\int_{-\infty}^{\infty}\frac{z^2 }{\frac{2\sigma}{V}+z^2} e^{-\frac{z^2}{2}} dz - 1\right).\label{eqDKnuExp}
\end{align}

\subsubsection{The functions $A(\alpha)$ and $\sigma_\varkappa(\alpha)$}

To formulate the main theorem of this section, we introduce a function $A(\alpha)$, which plays the central role throughout the paper.
\begin{definition}\label{defA}
For each $\alpha>0$,  $A=A(\alpha)$ is defined as a unique root of the equation
\begin{equation}\label{eqDKbetanu02}
F(\alpha-A) = \frac{\alpha}{(2\alpha+1)(\alpha-A)}
\end{equation}
with respect to $A$,
where
\begin{equation}\label{eqFAalpha}
F(x) := \frac{1}{(2\pi)^{1/2}}\int_{-\infty}^{\infty} \frac{1}{2x+z^2}e^{-\frac{z^2}{2}}dz \ \left(=
 \frac{\pi^{1/2}}{2}\frac{ e^{x} {\rm erfc}(x^{1/2})}{x^{1/2}}\right),\quad x>0,
\end{equation}
and ${\rm erfc}$ is the complementary error function.
\end{definition}

\begin{figure}[t]
\center
	   \includegraphics[width=0.4\textwidth]{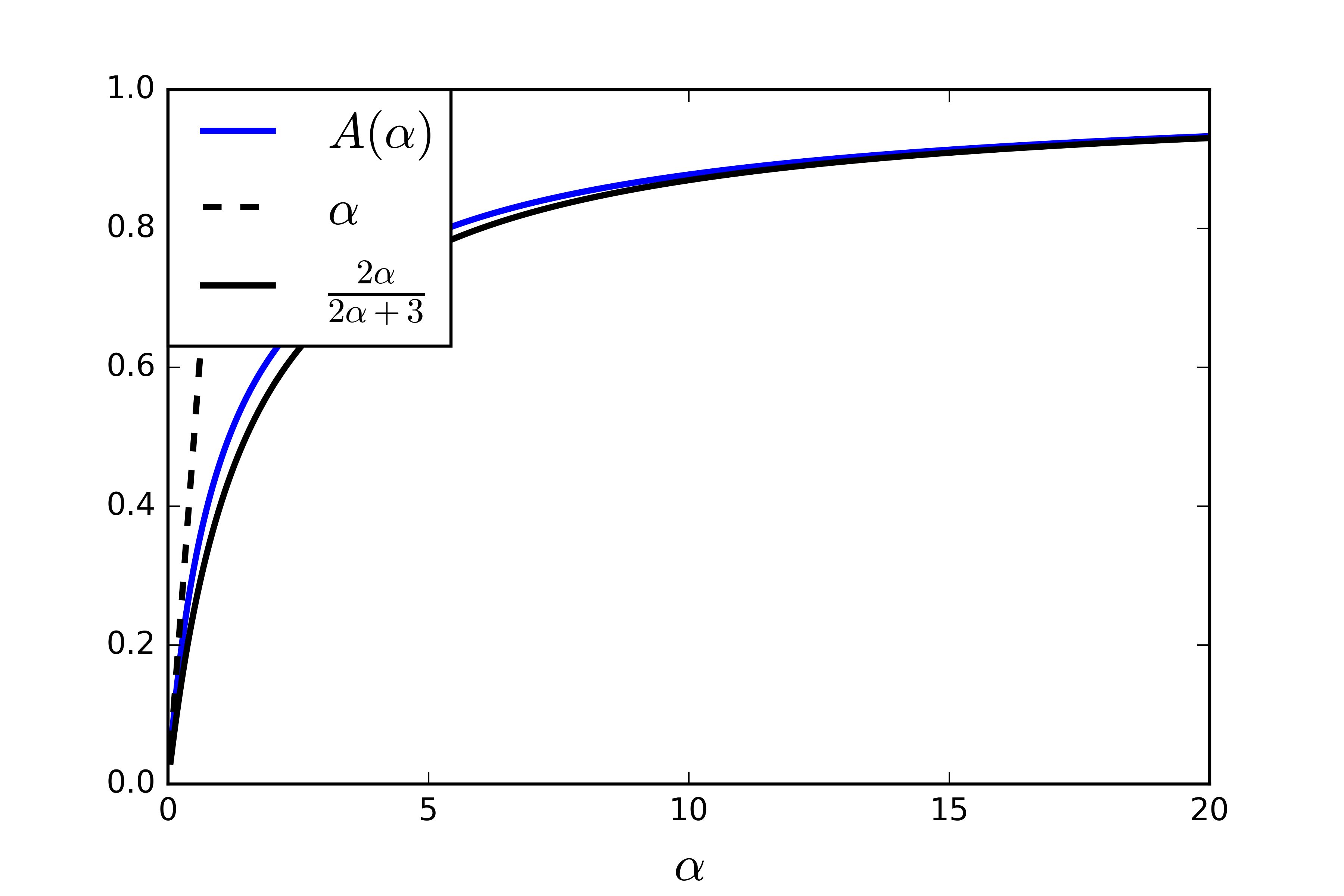}
\caption{The function $A(\alpha)$ from Definition~\ref{defA}.}\label{figAGalpha}
\end{figure}

The main properties of $A(\alpha)$ are given in the following lemma (see Fig.~\ref{figAGalpha}).
\begin{lemma}\label{lAProperties}
\begin{enumerate}
  \item\label{lAProperties1} Equation~\eqref{eqDKbetanu02} has a unique root $A(\alpha)$,
  \item\label{lAProperties2} $A(\alpha)$ is monotonically increasing,
  \item\label{lAProperties3} $\frac{2\alpha}{2\alpha+3}<A(\alpha)<\min(\alpha,1)$,
  \item\label{lAProperties4} $A(\alpha)$ satisfies the differential equation
    \begin{equation}\label{eqAODE}
  A'(\alpha) = 1 - \frac{2(\alpha-A(\alpha))}{(2\alpha+1)A(\alpha)},
  \end{equation}
  \item\label{lAProperties5} $A(\alpha)$ has the following asymptotics:
    \begin{equation}\label{eqPropertiesA}
  \begin{aligned}
  A(\alpha)=\alpha-k_0 \alpha^2 + o(\alpha^2)\ \text{as }\alpha\to 0,\\
   A(\alpha)=1-\frac{k_1}{\alpha} + o\left(\frac{1}{\alpha}\right) \ \text{as }\alpha\to \infty,
  \end{aligned}
  \end{equation}
  where $k_0=4/\pi$, $k_1=3/2$.
\end{enumerate}
\end{lemma}
\proof These properties are proved in Lemmas~\ref{lEalpha1}--\ref{lAmonotone}.
\endproof

\begin{definition}
For each $\varkappa\ge 0$, we define the functions (see Fig.~\ref{figSkappa_figCurvesCalpha}, left)
\begin{equation}\label{eqSigmaFunctions}
  \sigma_\varkappa(\alpha):= \left(1-\frac{\varkappa}{\alpha}\right)(\alpha-A(\alpha))V,\quad \alpha>0.
\end{equation}
We remind that $V=\bbV[\by]$.
\end{definition}

\subsubsection{Estimation of the variance}

\begin{figure}[t]
     \begin{minipage}{0.35\textwidth}
	   \includegraphics[width=\textwidth]{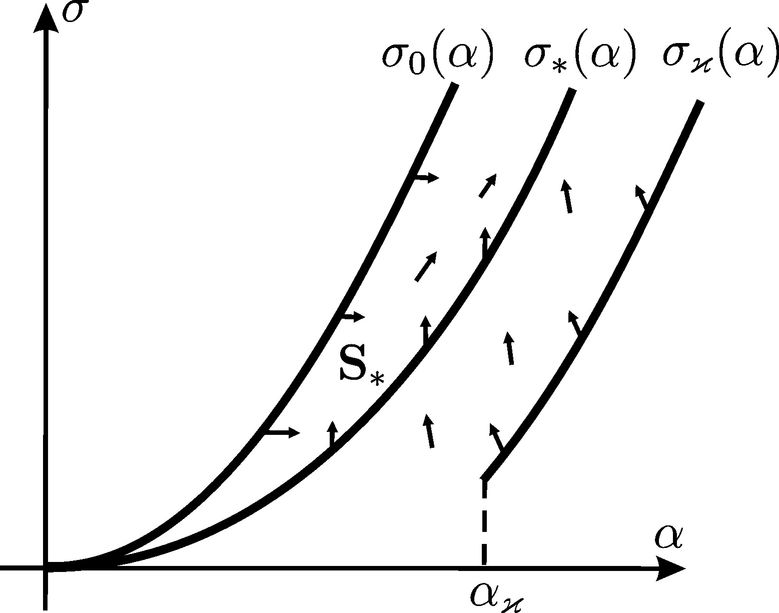}
    \end{minipage}
\hfill
    \begin{minipage}{0.5\textwidth}
    \includegraphics[width=\textwidth]{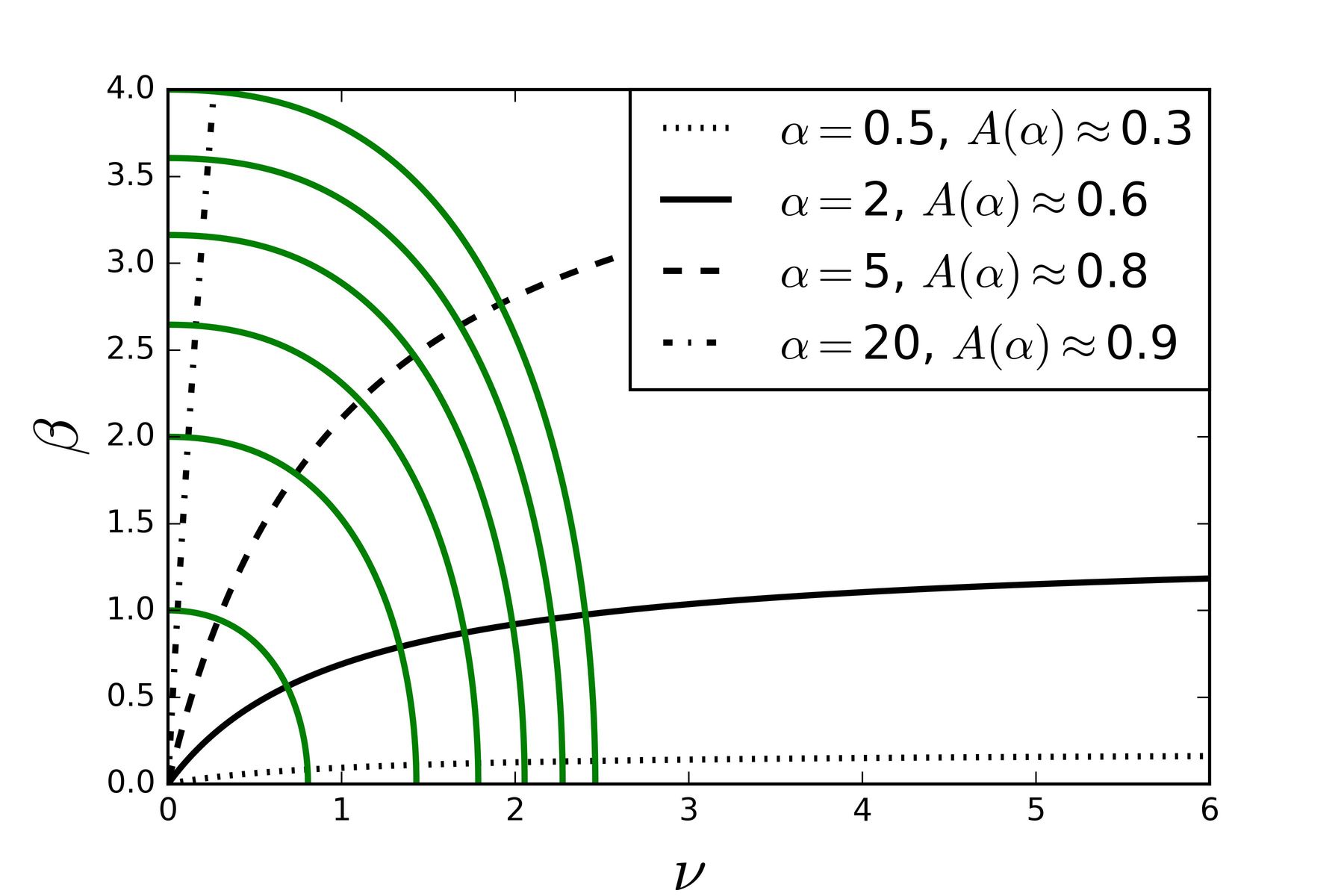}
    \end{minipage}
\caption{Left: The curves $\sigma_0(\alpha), \sigma_\varkappa(\alpha)$ given by~\eqref{eqSigmaFunctions}, the curve $\sigma_*(\alpha)$ from Theorem~\ref{thStrip}, item~\ref{thStrip1},  and the region~$\bS_*$ given by~\eqref{eqSStar}. The arrows indicate the directions of the vector field. Right:  Green lines are the curves given by~\eqref{eqCurvesNuBeta} for $C=1,4,7,10,13,16$. Black lines are the curves $\bC_{\alpha,V}$ given by~$\frac{\beta(\nu+1)}{\nu}=(\alpha-A(\alpha))V$ with $V=1$. The black and green curves are orthogonal to each other.}\label{figSkappa_figCurvesCalpha}
\end{figure}

\begin{figure}[t]
	\begin{minipage}{0.35\textwidth}
	   \includegraphics[width=\textwidth]{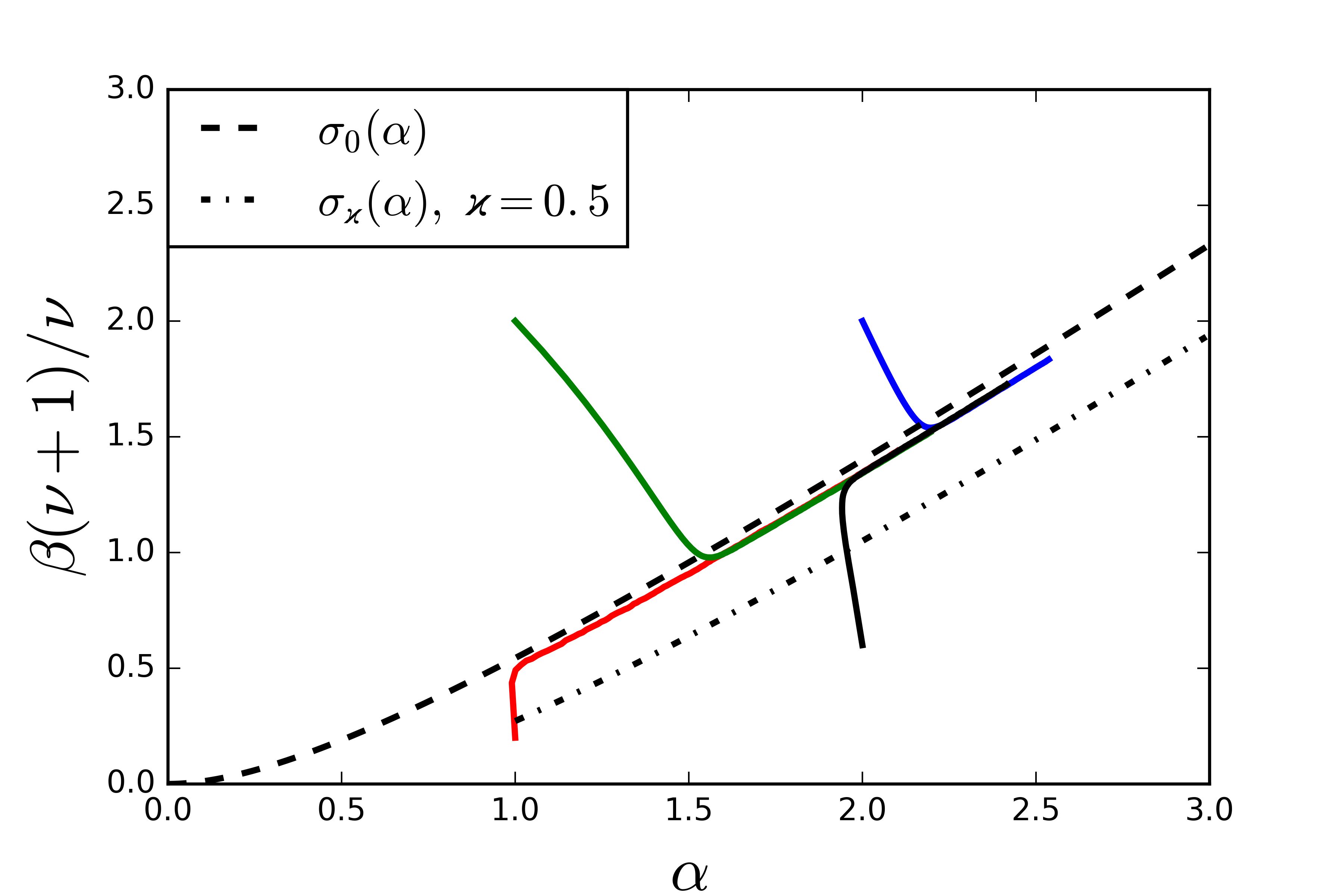}
	\end{minipage}
\hfill
	\begin{minipage}{0.31\textwidth}
       \includegraphics[width=\textwidth]{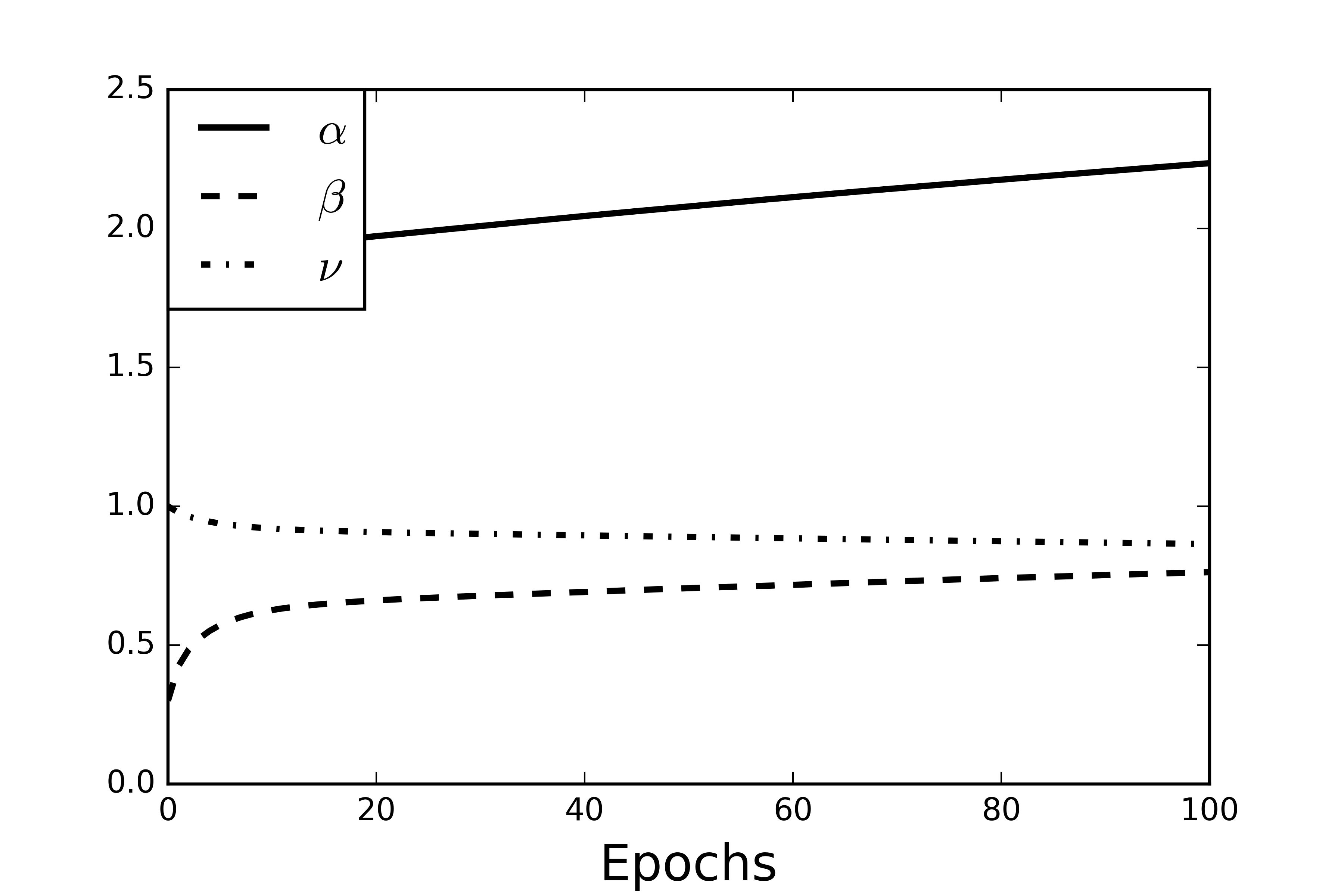}
	\end{minipage}
\hfill
	\begin{minipage}{0.31\textwidth}
       \includegraphics[width=\textwidth]{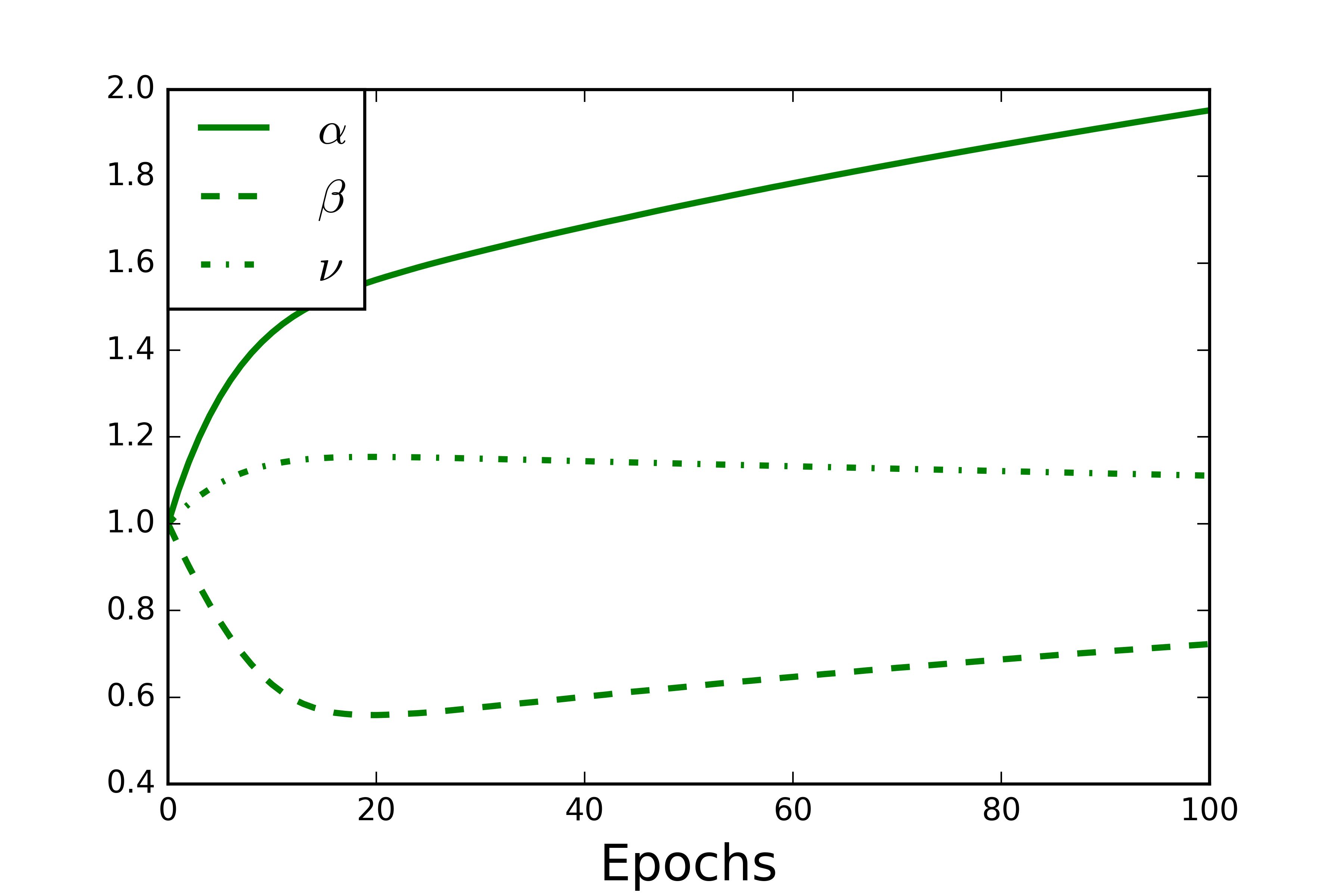}
	\end{minipage}
\caption{Left: Several trajectories obtained via iterating~\eqref{eqGradientDescent} for 2000 samples drawn from the normal distribution with mean $0$ and variance $1$. Middle/Right: Graphs of $\alpha,\beta,\nu$ plotted versus the number of epochs, corresponding to the lower-right/upper-left trajectory in the left figure.}\label{figTrajectories3D}
\end{figure}

The main result of this section (illustrated by Figures~\ref{figSkappa_figCurvesCalpha} and~\ref{figTrajectories3D}) is as follows.

\begin{theorem}\label{thStrip}
\begin{enumerate}
  \item\label{thStrip1} There is a smooth increasing function $\sigma_*(\alpha)$, $\alpha>0$, such that
  \begin{enumerate}
	\item $ \dot \alpha = 0 $ on the curve $ (\alpha, \sigma_* (\alpha)) $,
    \item $\lim\limits_{\alpha\to 0}\sigma_*(\alpha) = 0$ and $\lim\limits_{\alpha\to \infty}\sigma_*(\alpha) = \infty$,

    \item $\sigma_*(\alpha)<\sigma_0(\alpha)$ for all $\alpha>0$,

    \item   for any $\varkappa>0$, there exists $\alpha_\varkappa>0$ such that
  $$
   \sigma_*(\alpha)>\sigma_\varkappa(\alpha)\quad \text{for all}\ \alpha>\alpha_\varkappa,
  $$
    \item the region
\begin{equation}\label{eqSStar}
\bS_*:=\left\{(\alpha,\beta,\nu)\in\bbR^3:\alpha,\beta,\nu>0\ \text{and } \sigma_*(\alpha)<\frac{\beta(\nu+1)}{\nu} < \sigma_0(\alpha)\right\}
\end{equation}
 is forward invariant for system~\eqref{eqODE3}.
  \end{enumerate}

  \item\label{thStrip2} For any $\alpha(0),\beta(0),\nu(0)>0$, there exists a time moment $t_0$ depending on the initial condition such that for all $t>t_0$, $(\alpha(t),\beta(t),\nu(t))\in \bS_*$,  $\dot\alpha(t),\dot\beta(t)>0$, $\dot\nu(t)<0$.

  \item\label{thStrip4} For any $\alpha(0),\beta(0),\nu(0)>0$, there is $C>0$ depending on the initial conditions such that the points $(\nu(t),\beta(t))$ for all $t\ge 0$ lie on the integral curve
  \begin{equation}\label{eqCurvesNuBeta}
  {\beta^2} +  {\nu^2} + \frac{2\nu^3}{3} = C
  \end{equation}
      of the equation
    \begin{equation}\label{eqCurvesNuBetaODE}
    \frac{d\beta}{d\nu} = -\frac{\nu(\nu+1)}{\beta}.
  \end{equation}

  \item\label{thStrip3} For any $\alpha(0),\beta(0),\nu(0)>0$, we have
  $$
  \alpha(t)\to\infty,\ A(\alpha(t))\to 1,\ \nu(t)\to 0,\ \beta(t)\to\beta_*\quad\text{as } t\to\infty,
  $$
  where $\beta_*:= \left({\beta^2(0)} + {\nu^2(0)} + \frac{2\nu^3(0)}{3}\right)^{1/2}$.

\end{enumerate}
\end{theorem}

Theorem~\ref{thStrip} immediately implies the following corollary about the asymptotics of the variance $\Vest$ in~\eqref{eqStandardCPEstimate} for the predictive Student's t-distribution.
\begin{corollary}\label{corVestAlphaInfintiy}
  For any $\alpha(0),\beta(0),\nu(0)>0$, we have
\[
\left(1-\frac{\varkappa}{\alpha}\right)\frac{\alpha-A(\alpha)}{\alpha-1} < \frac{\Vest}{V} < \frac{\alpha-A(\alpha)}{\alpha-1} \quad\text{for all large enough}\ t.
\]
In particular,
$$
\Vest \to V\quad\text{as } t\to\infty.
$$
\end{corollary}
\begin{proof}
From Theorem \ref{thStrip}, item \ref{thStrip1}, we have by definition of $ S_* $
\[
\left(1-\frac{\varkappa}{\alpha}\right)(\alpha-A(\alpha))  V < \frac{\beta(\nu+1)}{\nu} < (\alpha-A(\alpha)) V \quad\text{for all sufficiently large}\ t.
\]
Deviding these inequalities by $\alpha-1$ and recalling that $ \alpha (t) \to \infty $ as $ t \to \infty $  and $A(\alpha)\to 1$ as $\alpha\to\infty$, we obtain the desired result.
\end{proof}

\begin{remark}
  One can show that $\varkappa/\alpha$ in the definition of the function $\sigma_\varkappa(\alpha)$ can be replaced by $\tilde\varkappa/\alpha^2$ with a sufficiently large $\tilde\varkappa$. In particular, the asymptotics in Corollary~\ref{corVestAlphaInfintiy} will assume the form
\[
\left|\frac{\Vest}{V} - 1\right| = O(\alpha^{-2})\quad\text{as }\ t\to\infty.
\]
   The proof   would require obtaining an extra term in the asymptotics of $A(\alpha)$ as $\alpha\to\infty$. However, we will not elaborate on these details.
\end{remark}

\begin{remark}\label{remVestAlphaInfinity}
  Suppose the number of observations tends to infinity. Then in the standard conjugate prior update~\eqref{eqCPUpdade}, the parameters $\alpha,\beta,\nu$ tend to infinity and the estimated mean and variance given by~\eqref{eqStandardCPEstimate} converge to the ground truth mean $\bbE[\by]$ and variance $V=\bbV[\by]$, while
  $$
  \bbE[\bmu]\to\bbE[\by],\quad \bbV[\bmu]\to 0,\quad \bbE[\btau]\to\frac{1}{\bbV[\by]},\quad \bbV[\btau]\to 0.
  $$

  The situation is quite different in Theorem~\ref{thStrip}. Although the parameter $\alpha$ tends to infinity,  $\beta$ converges to a finite positive value and $\nu$ converges to zero. Nevertheless, the estimated  variance $\Vest$ in Corollary~\ref{corVestAlphaInfintiy} converges to the ground truth variance $V=\bbV[\by]$, while (due to \eqref{eqEVmu}, \eqref{eqEVtau}, and~\eqref{eqStandardCPEstimate})
  $$
  \bbE[\bmu]\to\bbE[\by],\quad \bbV[\bmu]\to \bbV[\by],\quad \bbE[\btau]\to\infty,\quad \bbV[\btau]\to \infty.
  $$
\end{remark}

\subsection{Dynamics of $\alpha,\beta,\nu$. Proof of Theorem~\ref{thStrip}}

First, we show that $\bbE\left[\frac{\partial K}{\partial \beta}\right]$ and $\bbE\left[\frac{\partial K}{\partial \nu}\right]$ simultaneously vanish on the two-dimensional manifold
\begin{equation}\label{eqTwoDimSurface}
 \left\{(\alpha,\beta,\nu)\in\bbR^3:\alpha,\beta,\nu>0\ \text{and } \frac{\beta (\nu + 1)}{\nu} = \sigma_0(\alpha)\right\},
\end{equation}
where $\sigma_0(\alpha)$ is defined in~\eqref{eqSigmaFunctions}
Note that this manifold corresponds to the curve $\sigma=\sigma_0(\alpha)$ in Fig.~\ref{figSkappa_figCurvesCalpha}, left. We will also see that  $\dot\sigma,\dot\beta>0$ and $\dot\nu<0$ in $\bS_*$.

\begin{lemma}\label{lDKbetanu}
 We have
  \begin{equation}\label{eqDKbetanu0}
    \dot\beta=\dot\nu =0 \quad \text{if } \sigma=\sigma_0(\alpha),
  \end{equation}
    \begin{equation}\label{eqDKbetanu1}
    \begin{aligned}
    &\dot\sigma, \dot\beta <0,\ \dot\nu>0 & & \text{if } \sigma>\sigma_0(\alpha),\\
    &\dot\sigma, \dot\beta >0,\ \dot\nu<0  & & \text{if } \sigma<\sigma_0(\alpha).
    \end{aligned}
  \end{equation}
 \end{lemma}
\proof This lemma is proved in Appendix~\ref{appendixLemmaDKbetanu}.
\endproof

Now we show that the trajectories $(\nu(t),\beta(t))$ lie on curves  that do not depend on $\alpha$ or $V$, see the green lines in Fig.~\ref{figSkappa_figCurvesCalpha} (right).

\begin{lemma}\label{lCurvesNuBeta}
Let $\beta(t),\nu(t)$ ($t>0$) satisfy the last two equations in~\eqref{eqODE3} (for an arbitrary $\alpha(t)>0$). Then there is $C>0$ such that all the points $(\nu(t),\beta(t))$ belong to the integral curve~\eqref{eqCurvesNuBeta} of the equation~\eqref{eqCurvesNuBetaODE}.
\end{lemma}
\proof This lemma is proved in Appendix~\ref{appendixLemmaCurvesNuBeta}.
\endproof

Now we show that $\bbE\left[\frac{\partial K}{\partial \alpha}\right]$ is strictly negative on the manifold~\eqref{eqTwoDimSurface}, and, hence, neither system~\eqref{eqODE4}, nor system~\eqref{eqODE3} possesses an equilibrium.

\begin{lemma}\label{lDKalpha}
 We have
  \begin{equation}\label{eqDKalpha}
    \dot\alpha>0\quad \text{if } \sigma=\sigma_0(\alpha),\ \alpha>0.
  \end{equation}
  Moreover, for any $\varkappa>0$, there exists $\alpha_\varkappa>0$ such that
  \begin{equation}\label{eqDKalphaKappa}
    \dot\alpha<0 \quad\text{if } \sigma=\sigma_\varkappa(\alpha),\
     \alpha>\alpha_\varkappa.
  \end{equation}
\end{lemma}
\proof This lemma is proved in Appendix~\ref{appendixLemmaDKalpha}.
\endproof

\proof[Proof of Theorem~$\ref{thStrip}$] The arguments below are illustrated by Fig.~\ref{figSkappa_figCurvesCalpha}.

\textit{Item~\ref{thStrip1}.} Note that, for each fixed $ \sigma $, the function $\bbE\left[\frac{\partial K}{\partial \alpha}\right]$ is monotonically decreasing in~$ \alpha $. Furthermore by Lemma \ref{lDKalpha}, we have $ \dot \alpha > 0 $ on the curve $ \sigma=\sigma_0 (\alpha) $. On the other hand, $ \dot \alpha < 0 $ for large $ \alpha $, since $\Psi(\alpha)-\Psi(\alpha+1/2)\to 0$ as $\alpha\to\infty$. Thus, for each fixed $ \sigma $ there exists a unique value $ \alpha_* (\sigma) $ such that $ \dot \alpha = 0$. Moreover, since $ \bbE\left[\frac{\partial K}{\partial \alpha}\right] $ depends monotonically on $ \sigma $ and $ \alpha $, the function $ \alpha_*(\sigma) $ is smooth and can be inverted to a function $ \sigma_*(\alpha) $ by the inverse function theorem. By construction, $\sigma_*(\alpha)$ satisfies all the properties in Theorem~\ref{thStrip}, item~\ref{thStrip1}.

\textit{Item \ref{thStrip2}.} We argue by contradiction. Suppose $(\alpha(t),\beta(t),\nu(t))\notin\bS_*$ for all $t\ge0$. Then $ \alpha(t) $ either decreases for all $t\ge0$ or increases for all $t\ge0$, since the trajectory cannot cross the manifold $ \sigma=\sigma_*(\alpha) $. Suppose that it decreases. Since $ \sigma_* (\alpha) $  increases in $ \alpha $, it follows that $ \sigma(t) $ remains bounded. Furthermore, by Lemma~\ref{lDKbetanu}, $ \sigma(t) $ increases. Hence, there exists $ (\hat \alpha, \hat \sigma) $ such that
\begin{equation}\label{eqLimAlphaSigma0}
 \lim_{t \to \infty} (\alpha (t), \sigma(t)) = (\hat \alpha, \hat \sigma),\qquad \lim_{t \to \infty} (\dot\alpha (t), \dot\sigma(t)) = (0,0).
\end{equation}
In particular, $\bbE\left[\frac{\partial K}{\partial \alpha}\right]$ has to vanish at $(\hat \alpha, \hat \sigma)$. This may happen only if $(\hat \alpha, \hat \sigma)$ belongs to the curve $ \sigma=\sigma_*(\alpha) $. However, due to Lemma~\ref{lDKbetanu},  $\lim_{t\to\infty}\dot\sigma(t)$ cannot vanish in this case, which is a contradiction with~\eqref{eqLimAlphaSigma0}.  A similar argument applies if $\alpha(t)$ increases.

\textit{Item \ref{thStrip4}} follows from Lemma~\ref{lCurvesNuBeta}.

\textit{Item \ref{thStrip3}.} Due to item~\ref{thStrip2}, we can assume that $t_0=0$, so that $(\alpha(0),\beta(0),\nu(0))\in\bS_*$. By Lemma~\ref{lDKbetanu}, $\alpha(t)$ and $\beta(t)$ are increasing, while $\nu(t)$ is decreasing. Furthermore, by Lemma~\ref{lCurvesNuBeta}, $\beta(t)$ is bounded for all $t$.  Let us show that $\nu(t)\to 0$ as $t\to\infty$. Since the right-hand side in~\eqref{eqDKnuExp} has a singularity only for $\nu=0$, it remains to exclude the following two cases.

{\bf Case 1:} $\nu(t)\to \tilde\nu$ as $t\to\infty$ for some $\tilde\nu>0$. In this case, $\beta(t)\to\tilde\beta$ for some finite $\tilde\beta>0$, and hence $\alpha(t)\to\tilde\alpha$ for some finite $\tilde\alpha>0$ since the trajectory must stay in $\bS_*$.
Therefore, $(\tilde\alpha,\tilde\beta,\tilde\nu)$ must be an equilibrium of system~\eqref{eqODE3}. This contradicts Lemmas~\ref{lDKbetanu} and~\ref{lDKalpha}.

{\bf Case 2:} $\nu(t)\to 0$ as $t\to T$ for some finite $T>0$. In this case, $\beta(t)\to\tilde\beta$ as $t\to T$ for $\tilde\beta>0$ and hence $\alpha(t)\to\infty$ as $t\to T$ since the trajectory must stay in $\bS_*$. But this is possible only if $\bbE\left[\frac{\partial K}{\partial \alpha}\right]$ is unbounded as $\alpha,\sigma\to\infty$, which is not the case due to~\eqref{eqDKalphaExp}.

By having excluded Cases 1 and 2, we see that $\nu(t)\to 0$ as $t\to\infty$. Then $\beta(t)\to\tilde\beta$, where $\tilde\beta:= \left({\beta^2(0)} + {\nu^2(0)} + \frac{2\nu^3(0)}{3}\right)^{1/2}$ by to Lemma~\ref{lCurvesNuBeta}. Hence, $\sigma(t)=\frac{\beta(t)(\nu(t)+1)}{\nu(t)}\to\infty$. As  the trajectory must stay in $\bS_*$, it follows that $\alpha(t)\to\infty$. Finally, by Lemma~\ref{lAProperties}, $A(\alpha(t))\to 1$.
\endproof

\section{Dynamics of $m,\beta,\nu$ for a fixed $\alpha$}\label{secDynamicsFixedAlpha}

\subsection{Estimation of the variance. The curves of equilibria}

According to Theorem~\ref{thStrip},   neither system~\eqref{eqODE4}, nor system~\eqref{eqODE3} possesses an equilibrium. 
However, this is not the case any more if the training dataset contains output outliers. Experiments show that $\alpha$ then converges to a finite value. In this section, we analyze the dynamics $m,\beta,\nu$ in~\eqref{eqGradientDescent} for a fixed $\alpha$. We consider the update
\begin{equation}\label{eqGradientDescentFixedAlpha}
m_{\rm new} := m - \lambda\frac{\partial K}{\partial m}, \quad
\beta_{\rm new} := \beta - \lambda\frac{\partial K}{\partial \beta}, \quad \nu_{\rm new} := \nu - \lambda\frac{\partial K}{\partial \nu}.
\end{equation}

As in Sec.~\ref{secDynamics4}, taking into account Theorem~\ref{thDKm}, we assume that the mean $m$ has already been learned: $m=\bbE[\by]$ (Condition~\ref{condmEqualsE}). Then
the corresponding approximating dynamical system is given by the two equations for $\beta,\nu$ from~\eqref{eqODE4}:
\begin{equation}\label{eqODE3FixedAlpha}
\dot\beta = -\bbE\left[\frac{\partial K}{\partial \beta}\right],\quad \dot\nu = -\bbE\left[\frac{\partial K}{\partial \nu}\right],
\end{equation}
where the right-hand sides are explicitly given by~\eqref{eqDKbetaExp} and~\eqref{eqDKnuExp}.
We consider this system on the quadrant $\{(\nu,\beta)\in\bbR^2:\nu,\beta>0\}$. Due to~\eqref{eqDKbetanu1}, this quadrant is forward invariant.

\begin{remark}\label{rGradientFlow2Eq}
As in Remark~\ref{rGradientFlow4Eq}, we conclude from~\eqref{eqGradientDescentWAveraged} that system~\eqref{eqODE3FixedAlpha} defines the gradient flow with the potential  $D_{\rm KL}(q \| p_{\rm pred}(\cdot|w))$, where $p_{\rm pred}(\cdot|w))$ is the probability distribution of the Student's t-distribution $t_{2\alpha}(y|\bbE[\by],\beta(\nu+1)/(\nu\alpha))$.
\end{remark}

\begin{theorem}\label{thDynamicsFixedAlpha}
Let $\alpha>0$ be fixed. Then the following hold.
\begin{enumerate}
  \item\label{thDynamicsFixedAlpha1} Dynamical system~\eqref{eqODE3FixedAlpha} possesses a globally attracting family of equilibria lying on the curve
\begin{equation}\label{eqAttractingCurveEquil}
  \bC_{\alpha,V}:=\left\{(\nu,\beta)
  \in\bbR^2: \nu>0,\ \frac{\beta(\nu+1)}{\nu}=\sigma_0(\alpha)\right\},
\end{equation}
where $\sigma_0(\alpha)$ is defined in~\eqref{eqSigmaFunctions}.

  \item\label{thDynamicsFixedAlpha2} Each trajectory $(\nu(t),\beta(t))$ lies on one of the integral curves~\eqref{eqCurvesNuBeta}. If $(\nu(0),\beta(0))$ lies below the curve $\bC_{\alpha,V}$, then $\nu(t)$ decreases and converges to $\nu_*$ and $\beta(t)$ increases and converges to $\beta_*$. If $(\nu(0),\beta(0))$ lies above the curve $\bC_{\alpha,V}$, then $\nu(t)$ increases and converges to $\nu_*$ and $\beta(t)$ decreases and converges to $\beta_*$. In both cases, $(\nu_*,\beta_*)$ is the point of intersection of the corresponding integral curve and the curve of equilibria $\bC_{\alpha,V}$, see Fig.~\ref{figSkappa_figCurvesCalpha}.

  \item\label{thDynamicsFixedAlpha3} The family of integral curves~\eqref{eqCurvesNuBeta} is orthogonal to the family of the curves of equilibria $\{\bC_{\alpha,V}\}_{\alpha,V>0}$.
\end{enumerate}
\end{theorem}
\proof
Lemmas~\ref{lDKbetanu} and~\ref{lCurvesNuBeta} imply items~\ref{thDynamicsFixedAlpha1} and~\ref{thDynamicsFixedAlpha2}. Let us prove item~\ref{thDynamicsFixedAlpha3}. Assume that $(\nu,\beta)$ is a point of intersection of the curves
\begin{equation}\label{eqOrthogonalCurves}
\begin{aligned}
\beta & = f(\nu):=(\alpha-A)V \frac{\nu}{\nu+1} & &  \text{(a curve of equilibria)},\\
\beta & = g(\nu):= \left(C - {\nu^2} - \frac{2\nu^3}{3}\right)^{1/2} & & \text{(an integral curve)}.
\end{aligned}
\end{equation}
Then $f'(\nu) = \frac{(\alpha-A)V}{(\nu+1)^2} = \frac{\beta}{\nu(\nu+1)}$. On the other hand, by Lemma~\ref{lCurvesNuBeta}, $g'(\nu)=-\frac{\nu(\nu+1)}{\beta}$. Thus, $f'(\nu)g'(\nu)=-1$, which implies the orthogonality of the two curves in~\eqref{eqOrthogonalCurves}.
\endproof

Theorem~\ref{thDynamicsFixedAlpha} immediately implies the following corollary.
\begin{corollary}\label{corVestAlphaFixed}
Let $\alpha>0$ be fixed. Then for any $\beta(0),\nu(0)>0$, we have
\begin{equation}\label{eqVarianceEst}
  \tVest(t):=\frac{\beta(t)(\nu(t)+1)}{(\alpha-A(\alpha))\nu(t)}\to V\quad\text{as } t\to\infty.
\end{equation}
\end{corollary}

\begin{figure}[t]
	\begin{minipage}{0.32\textwidth}
	   \includegraphics[width=\textwidth]{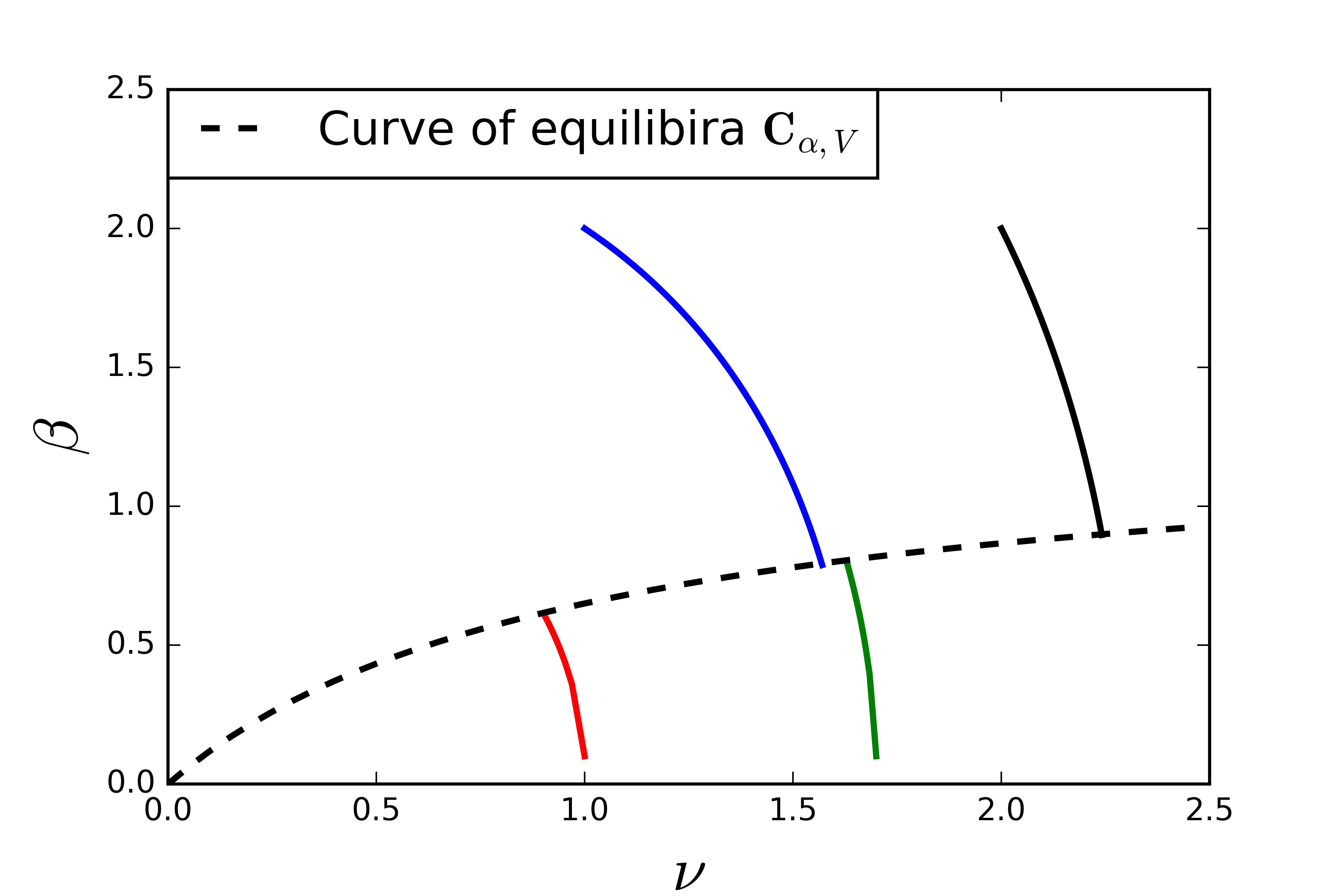}
	\end{minipage}
\hfill
	\begin{minipage}{0.32\textwidth}
       \includegraphics[width=\textwidth]{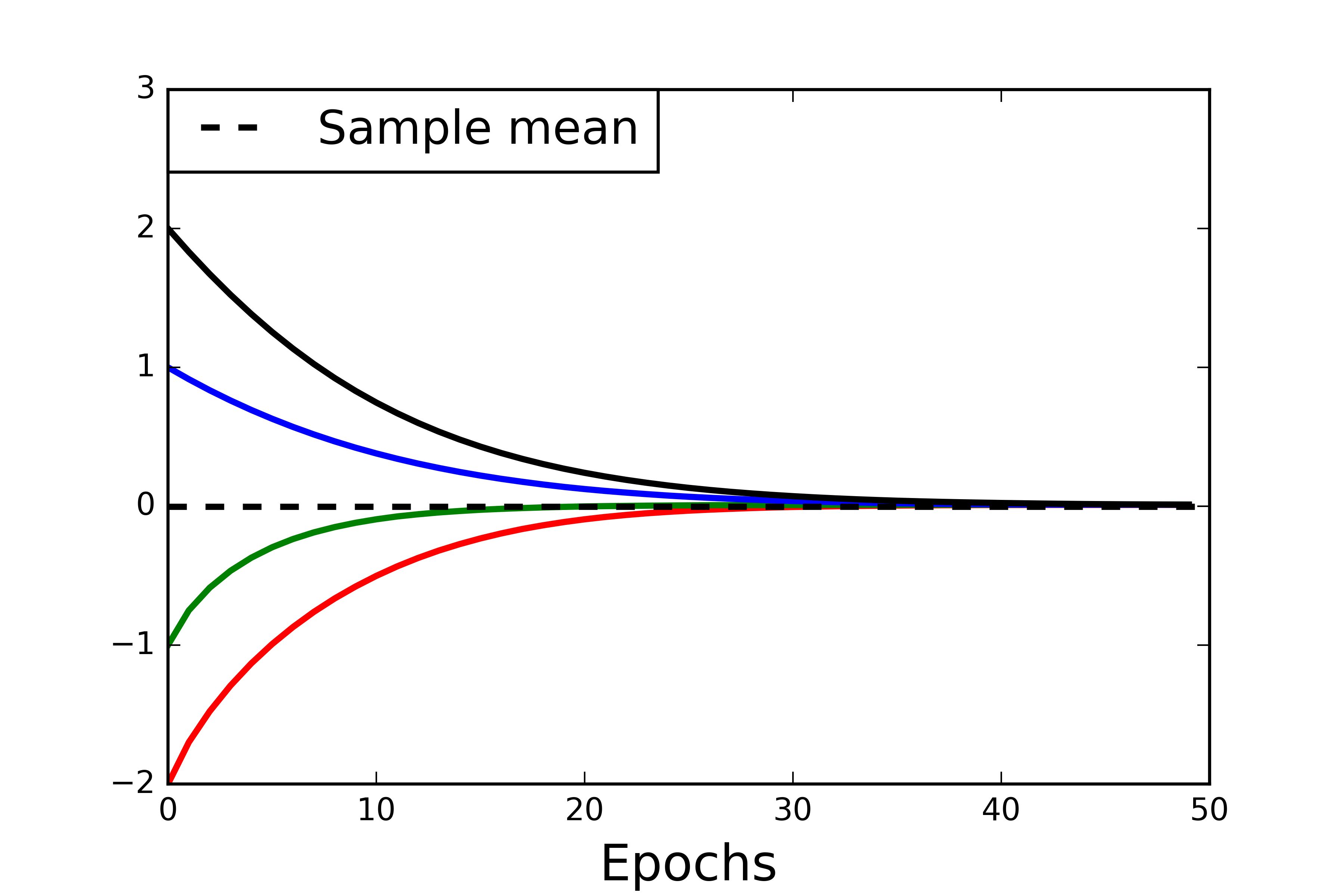}
	\end{minipage}
\hfill
	\begin{minipage}{0.32\textwidth}
       \includegraphics[width=\textwidth]{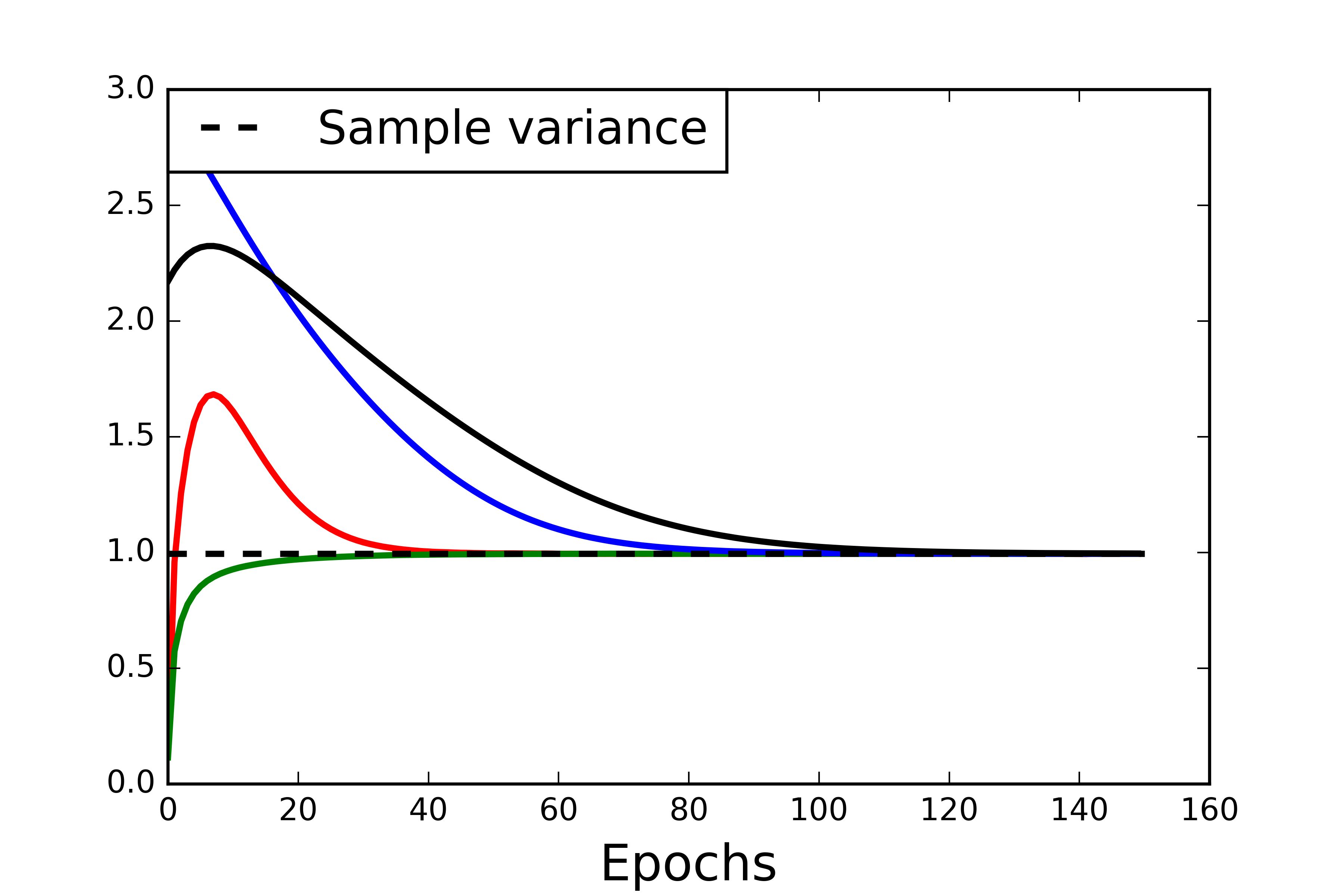}
	\end{minipage}
\caption{Several trajectories obtained via iterating~\eqref{eqGradientDescentFixedAlpha} with 2000 samples drawn from the normal distribution with mean $0$ and variance $1$. The parameter $\alpha$ is fixed: $\alpha=2$, $A(\alpha)\approx 0.619$. Left: plane $(\nu,\beta)$, with the initial condition $m(0)=0$ for all trajectories. Middle and Right: means $m$ and variances $\tVest$ plotted versus the number of epochs.  The initial conditions for $\beta$ and $\nu$ are the same as for the trajectories of the respective colors in the left plot, while $m(0)=-2,-1,1,2$.}\label{figTrajectories2D}
\end{figure}

Figure~\ref{figSkappa_figCurvesCalpha} (right) shows the mutual configuration of the integral curves~\eqref{eqCurvesNuBeta}  and the curves of equilibria $\bC_{\alpha,V}$ corresponding to different $\alpha$. Figure~\ref{figTrajectories2D} (left) shows several trajectories in the $(\nu,\beta)$ plane converging to the curve of equilibria $\bC_{\alpha,V}$. Figure~\ref{figTrajectories2D} (middle and right) shows that taking initial conditions with $m(0)\ne\bbE[\by]$ still yields the proper convergence of the estimated mean $m_{\rm est}:=m$ and the estimated variance $\tVest$ in~\eqref{eqVarianceEst}.

\begin{remark}
  Due to Theorem~\ref{thDynamicsFixedAlpha}, each trajectory $(\nu(t),\beta(t))$ of system~\eqref{eqODE3FixedAlpha} can be obtained by solving the scalar differential equation
  \begin{equation}\label{eqNuODE}
  \dot\nu = -\bbE\left[\frac{\partial K}{\partial \nu}\right]\bigg|_{\beta=g(\nu)},
  \end{equation}
  where $g(\nu)$ is defined in~\eqref{eqOrthogonalCurves} with a fixed $C>0$ (uniquely determined by $\nu(0),\beta(0)$).

  Furthermore, one can use other functions $\tilde g(\nu)$ in~\eqref{eqNuODE} instead of $g(\nu)$. Due to~\eqref{eqDKbetanu0}, the resulting ODE would still have an equilibrium $\nu_*$ such that $(\nu_*,\tilde g(\nu_*))\in\bC_{\alpha,V}$, and at this equilibrium, we would have
  $$
  V=\frac{\beta_*(\nu_*+1)}{(\alpha-A(\alpha))\nu_*},\quad \beta_*:=\tilde g(\nu_*).
  $$
  One can also show that this equilibrium is globally stable for a broad class of functions $\tilde g(\nu)$.

  However, the function $g(\nu)$ from~\eqref{eqOrthogonalCurves}, corresponding to the integral curve~\eqref{eqCurvesNuBeta}, is {\em optimal} in the following sense, see Fig.~\ref{figGradientDescent}.
 \begin{figure}[t]
    \center
	   \includegraphics[width=0.5\textwidth]{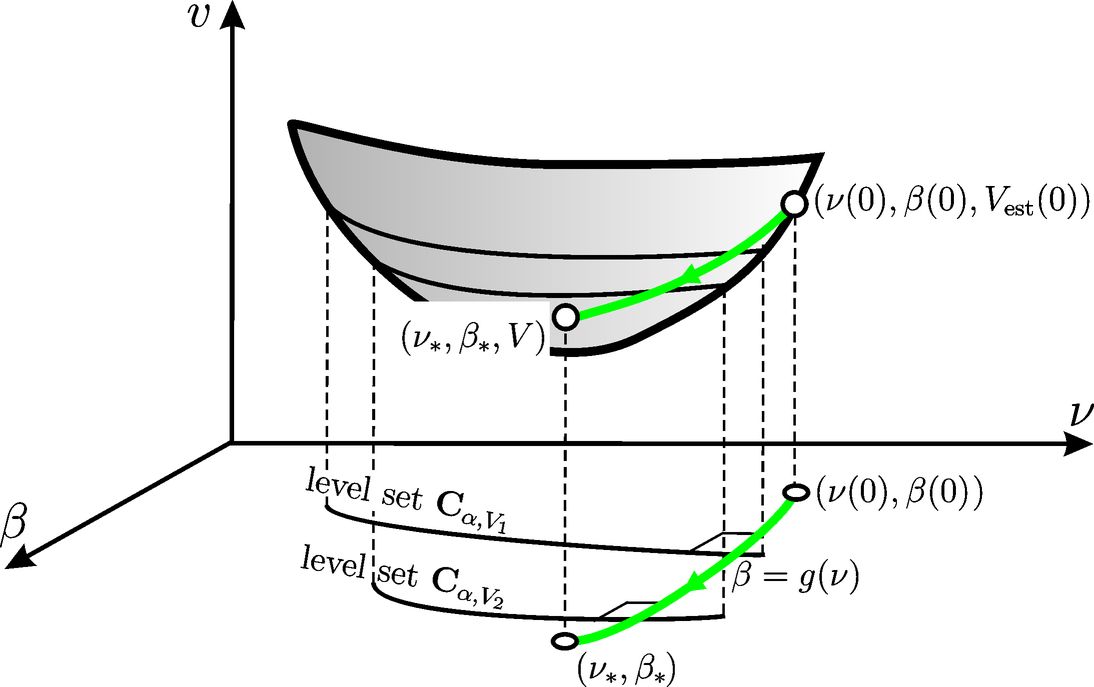}
    \caption{Schematic surface~\eqref{eqSurfaceNuBetaV}, its level sets $\bC_{\alpha,V_1}$ and $\bC_{\alpha,V_2}$, and the trajectory (in green) connecting the initial point $(\nu(0),\beta(0),\Vest(0))$ and the target point $(\nu_*,\beta_*,V)$. The projection $\beta=g(\nu)$ of the trajectory is orthogonal to the level sets, i.e., the trajectory follows the gradient descent on the surface~\eqref{eqSurfaceNuBetaV}.}\label{figGradientDescent}
\end{figure}
  Consider the two-dimensional surface in $\bbR^3$
  \begin{equation}\label{eqSurfaceNuBetaV}
  \left\{(\nu,\beta,v)\in\bbR^3: v = \frac{\beta(\nu+1)}{(\alpha-A(\alpha))\nu}\right\}
  \end{equation}
  (with $\alpha$ fixed). Then the initial point $(\nu(0),\beta(0),\tVest(0))$ (where $\tVest(t)$ is defined in~\eqref{eqVarianceEst}) and the target point $(\nu_*,\beta_*,V)$ (where $\nu_*,\beta_*$ are defined in Theorem~\ref{thDynamicsFixedAlpha}, item~\ref{thDynamicsFixedAlpha2}) both lie on this surface. On the other hand, the curves $\{\bC_{\alpha,V_1}\}_{V_1>0}$ are the level sets of this surface. Hence, due to Theorem~\ref{thDynamicsFixedAlpha}, item~\ref{thDynamicsFixedAlpha3}, the curve $\beta=g(\nu)$ corresponds to the path of the {\em gradient descent} (or {\em ascent}) connecting the initial point $(\nu(0),\beta(0),\Vest(0))$ and the target point $(\nu_*,\beta_*,V)$.
\end{remark}

\begin{remark}
  It follows from Remark~\ref{rGradientFlow2Eq} and item~\ref{thDynamicsFixedAlpha3} in Theorem~\ref{thDynamicsFixedAlpha} that, for any fixed $\alpha$ and $V=\bbV[\by]$, the curve $\bC_{\alpha,V}$ is the set of minima of the potential $D_{\rm KL}(q \| p_{\rm pred}(\cdot|w))$, while all the other curves $\bC_{\alpha_1,V_1}$, $\alpha_1,V_1>0$, are the level sets of this potential.
\end{remark}

\begin{remark}
    The situation in Theorem~\ref{thDynamicsFixedAlpha} is different both from the standard CP update~\eqref{eqCPUpdade} and from the GCP update~\eqref{eqGradientDescent} (cf. Remark~\ref{remVestAlphaInfinity}). First, the parameter $\alpha$ is now fixed. Furthermore, each trajectory of system~\eqref{eqODE3FixedAlpha} (approximating the GCP update~\eqref{eqGradientDescentFixedAlpha}) converges to a finite equilibrium $(\nu_*,\beta_*)$, where $\nu_*,\beta_*>0$. Nevertheless, the estimated  variance $\tVest$ given by~\eqref{eqVarianceEst} again converges to the ground truth variance $V=\bbV[\by]$.
\end{remark}

\section{Role of a fixed $\alpha$}\label{secRoleFixedAlpha}

\subsection{Sensitivity to outliers}\label{subsecSensitivityOutliers}

\begin{figure}[t]
	\begin{minipage}{0.32\textwidth}
	   \includegraphics[width=\textwidth]{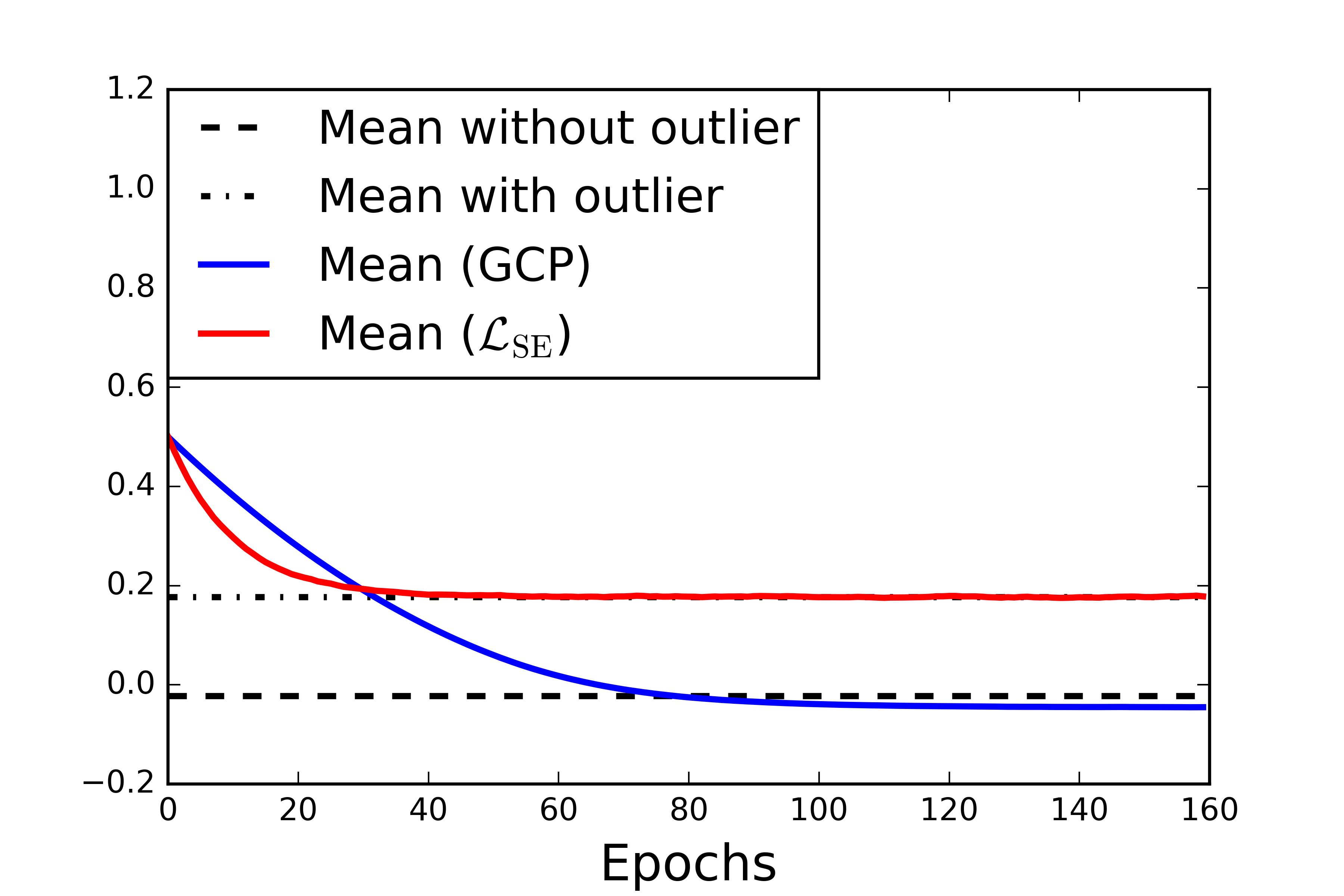}
	\end{minipage}
\hfill
	\begin{minipage}{0.32\textwidth}
       \includegraphics[width=\textwidth]{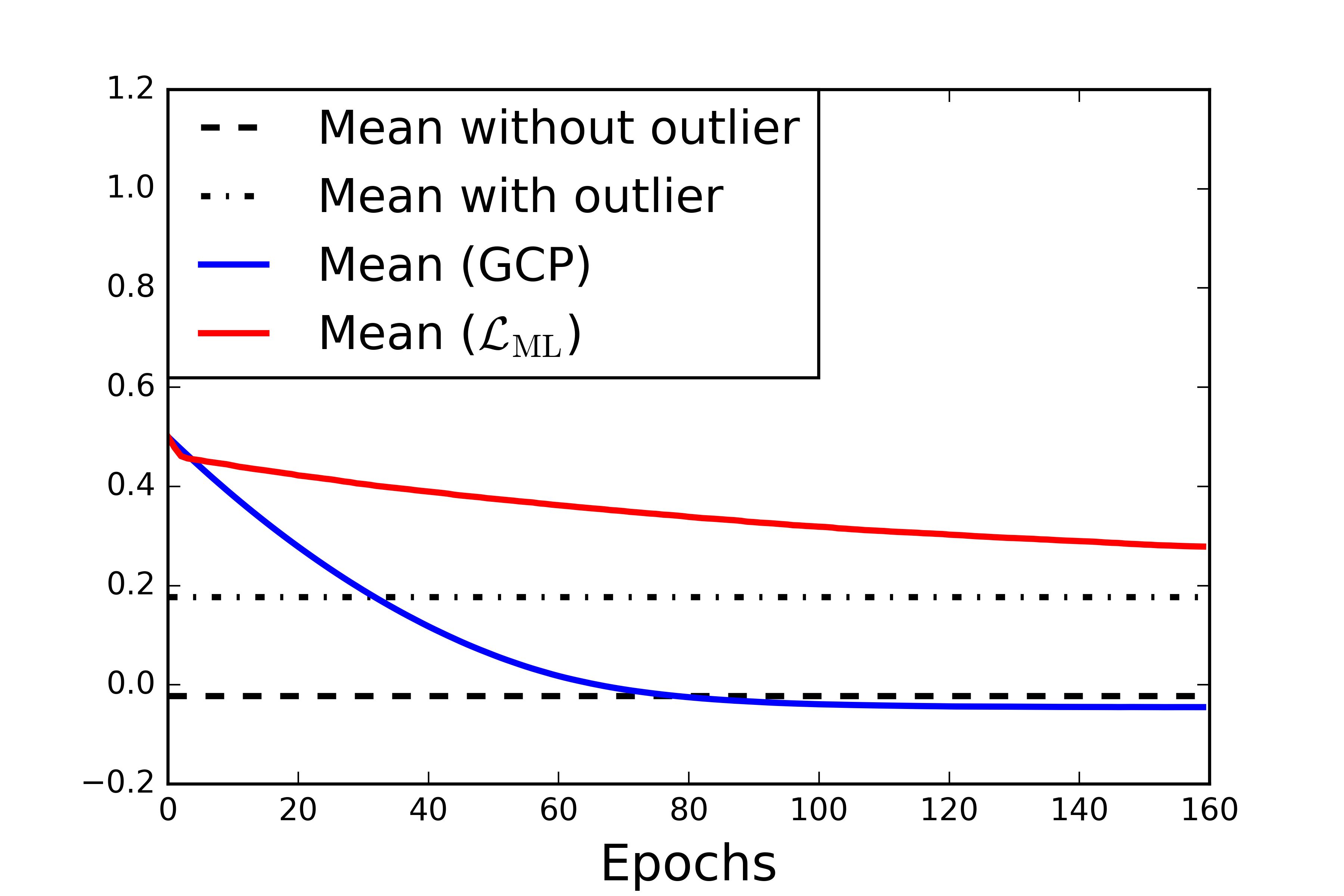}
	\end{minipage}
\hfill
	\begin{minipage}{0.32\textwidth}
       \includegraphics[width=\textwidth]{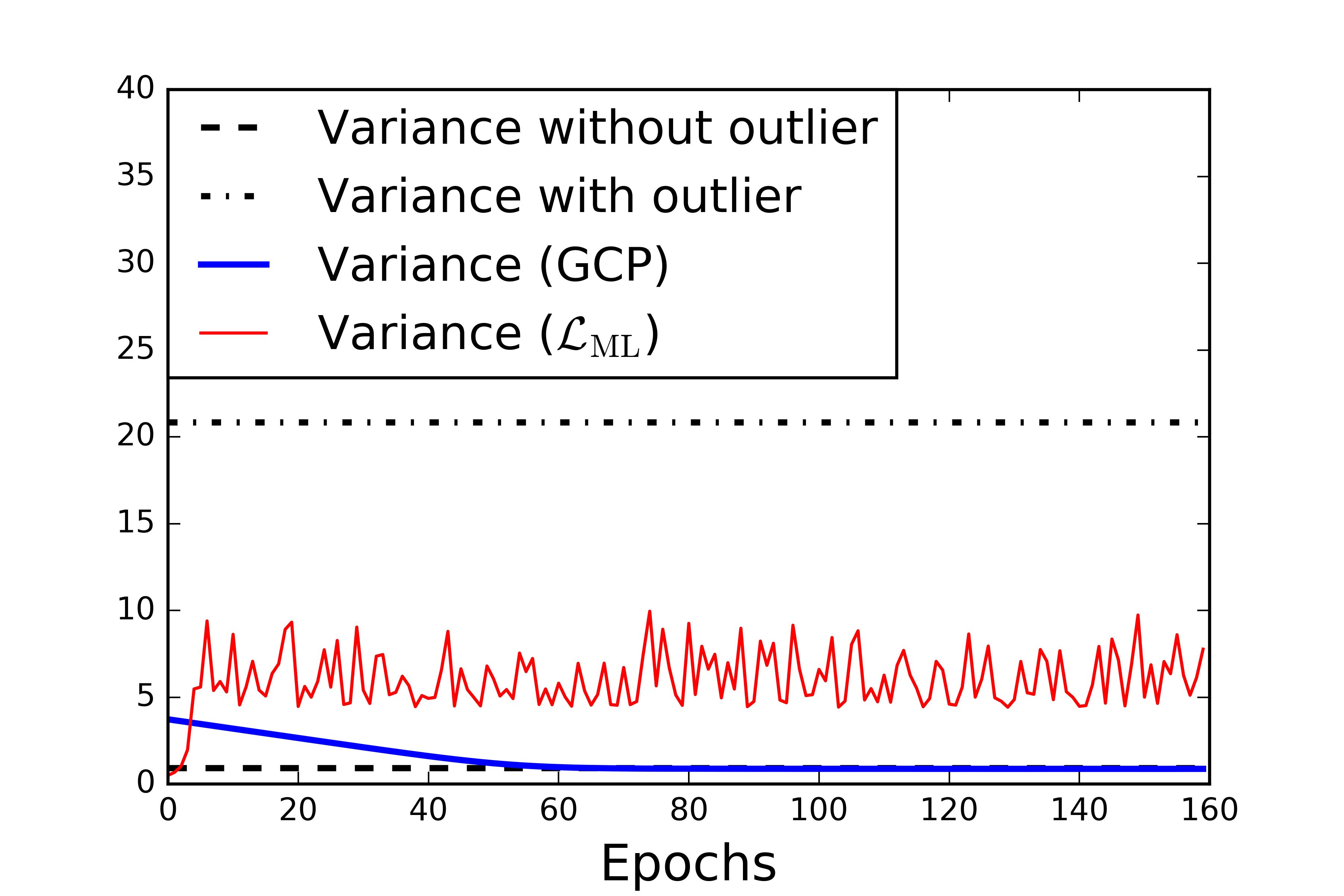}
	\end{minipage}
\caption{Fitting 500 samples drawn from the normal distribution with mean $0$ and variance~$1$ and supplemented by an outlier $ y =100$. Left: Means  fitted with the GCP update~\eqref{eqGradientDescentFixedAlpha} and, respectively, with the standard squared error loss~$\cL_{\rm SE}$. Middle/Right: Means/variances fitted with the GCP update~\eqref{eqGradientDescentFixedAlpha} and, respectively, via maximizing the likelihood, i.e., minimizing~$\cL_{\rm ML}$. For the GCP update, the parameter $\alpha$ is fixed: $\alpha=1$, $A(\alpha)\approx 0.46$, and the variance is estimated by $\tVest$ in~\eqref{eqVarianceEst}.}\label{figOutliers}
\end{figure}

It is well known that outliers essentially influence the estimate of the mean $m$ if one uses the standard squared error loss
$$
\cL_{\rm SE}( y ,m) = (m- y )^2.
$$
The same is true when one estimates both mean $m$ and precision $p$  via maximizing the log-likelihood of a normal distribution, or, equivalently, minimizing the loss
$$
\cL_{\rm ML}( y ,m,p) = p(m- y )^2 - \ln p.
$$
The reason is that, in both cases, the derivatives of the loss functions~$\cL_{\rm SE}$ and~$\cL_{\rm ML}$ with respect to $m$ are proportional to $m- y $, while the derivative of $\cL_{\rm ML}( y ,m,p)$ with respect to $p$  contains even $(m- y )^2$. It turns out that the GCP update~\eqref{eqGradientDescentFixedAlpha} is much less sensitive to outliers, see Fig.~\ref{figOutliers}. This can be explained by the fact that the derivatives of the KL divergence with respect to $m$, $\beta$ and $\nu$ are bounded with respect to $m- y $, see~\eqref{eqDKmPure}, \eqref{eqDKbetaPure}, and~\eqref{eqDKnuPure}. Moreover, $\frac{\partial K}{\partial m}$ even vanishes as $m- y \to\infty$. Another explanation is that the GCP update is equivalent to maximizing the likelihood of the Student's t-distribution (item~\ref{rGCPupdatePredictiveDistr0} in Remark~\ref{rGCPupdatePredictiveDistr}). It is known~\cite{Nadarajah08} that the optimal value of $m$ is different from the sample mean due to downweighting the outlying observations. In sections~\ref{subsecSynthetic} and~\ref{subsecRealWorldOutliers}, we further analyze the performance of the GCP neural networks on contaminated data sets in comparison with other neural network methods.

\subsection{Learning speed in clean and noisy regions}\label{subsecRoleAlpha}

\subsubsection{Observations}

When one approximates the parameters $m,\beta,\nu$ by neural networks, one represents these parameters as functions of an input variable $x\in  \bbR^m$ and of a set of weights $w\in\bbR^M$. Since neural networks have finite capacity ($M$ is finite), they cannot perfectly approximate $m,\beta,\nu$ for all $x$ simultaneously. Therefore, it is important to understand in which regions of the input space $\bbR^m$ the parameters are approximated better and in which worse, cf.~\cite{GurHannesLU}. This is directly related to the values of the gradients in the GCP update~\eqref{eqGradientDescentFixedAlpha}, which determine the learning speed. The faster the learning in a certain region occurs, the more influential this region is. In particular, we are interested in the learning speed in so called clean regions (where $V$ is small) compared with noisy regions (where $V$ is large).

Below, we will concentrate on the regime where the learning process starts and the initial conditions for $\beta$ and $\nu$ satisfy
\begin{equation}\label{eqInitialBetaNu}
\beta(0)\approx\nu(0)\approx 1.
\end{equation}
 This is often the case if $\beta$ and $\nu$ are approximated by neural networks with the softplus output, e.g.,
\begin{equation}\label{eqBetaSoftplus}
  \beta = \ln\left(1+e^w\right),
\end{equation}
where $w\in\bbR$ is the input of the softplus output.

In the observations below, we denote the learning speed of the mean and the variance by $\LerSp(m)$ and $\LerSp(\Var)$, respectively.

\begin{observation}\label{obsSmallalpha}
Let $\alpha$ be small.
\begin{enumerate}
  \item\label{obsSmallalpha1}  {In clean regions $($small $V)$$:$}   $\LerSp(m)$ is of order $1$ and $\LerSp(\Var)$ is of order $1$.
  \item\label{obsSmallalpha2}  {In noisy regions $($large $V)$$:$} $\LerSp(m)$ is of order $1/V$ and   $\LerSp(\Var)$ is of order $\alpha$.
\end{enumerate}
\end{observation}

\begin{observation}\label{obsLargealpha}
Let $\alpha$ be large.
\begin{enumerate}
  \item\label{obsLargealpha1}  {In clean regions $($small $V)$$:$} $\LerSp(m)$ is of order $\alpha$ and $\LerSp(\Var)$ is of order $1$.
  \item\label{obsLargealpha2}  {In noisy regions $($large $V)$$:$} $\LerSp(m)$ is of order $1/V$ and $\LerSp(\Var)$ is of order $\alpha$.
\end{enumerate}
\end{observation}

Observations~\ref{obsSmallalpha} and~\ref{obsLargealpha} are summarized in Table~\ref{table}. In particular, we see that the mean is always learned faster in clean regions. Taking large $\alpha$ further increases the learning speed of the mean in clean regions, but simultaneously increases the learning speed of variance in noisy regions compared with clean regions.
\begin{table}[t]
  \centering
    \begin{tabular}{c|c|c}

    $\alpha$ & $\LerSp(m)$ & $\LerSp(\Var)$ \\

    \hline

    Small & \text{cl.} $>$ \text{noisy} & \text{cl.} $>$ \text{noisy} \\

    \hline

    Large & \text{cl.} $\gg$ \text{noisy} & \text{cl.} $<$ \text{noisy}

  \end{tabular}
  \caption{Relative learning speed of estimated mean $m$ and variance $\tVest$ for clean (cl.) and noisy regions. Notation   ``$<$'' and ``$>$'' stands for a ``lower'' and a ``higher'' speed, and  ``$\gg$'' for a ``much higher'' speed.}\label{table}
\end{table}

\begin{observation}\label{obsValpha}
The values of $\beta_*$ and $\nu_*$ to which the trajectory of~\eqref{eqODE3FixedAlpha} will converge are determined by the value $(\alpha-A(\alpha))V$.
\begin{enumerate}
  \item\label{obsValpha1}  If $(\alpha-A(\alpha))V \ll 1$, then $\beta_*\approx 0$ and $\nu_*\approx 1$.
  \item\label{obsValpha2}  If $(\alpha-A(\alpha))V \gg 1$, then $\beta_*\approx 1$ and $\nu_*\approx 0$.
\end{enumerate}
Small values of $\beta$ and $\nu$ will lead to large gradients $\frac{\partial K}{\partial \beta}$ and $\frac{\partial K}{\partial \nu}$, respectively, which may cause large oscillations of $\tVest$.
\end{observation}

\begin{figure}[t]
	\begin{minipage}{0.32\textwidth}
	   \includegraphics[width=\textwidth]{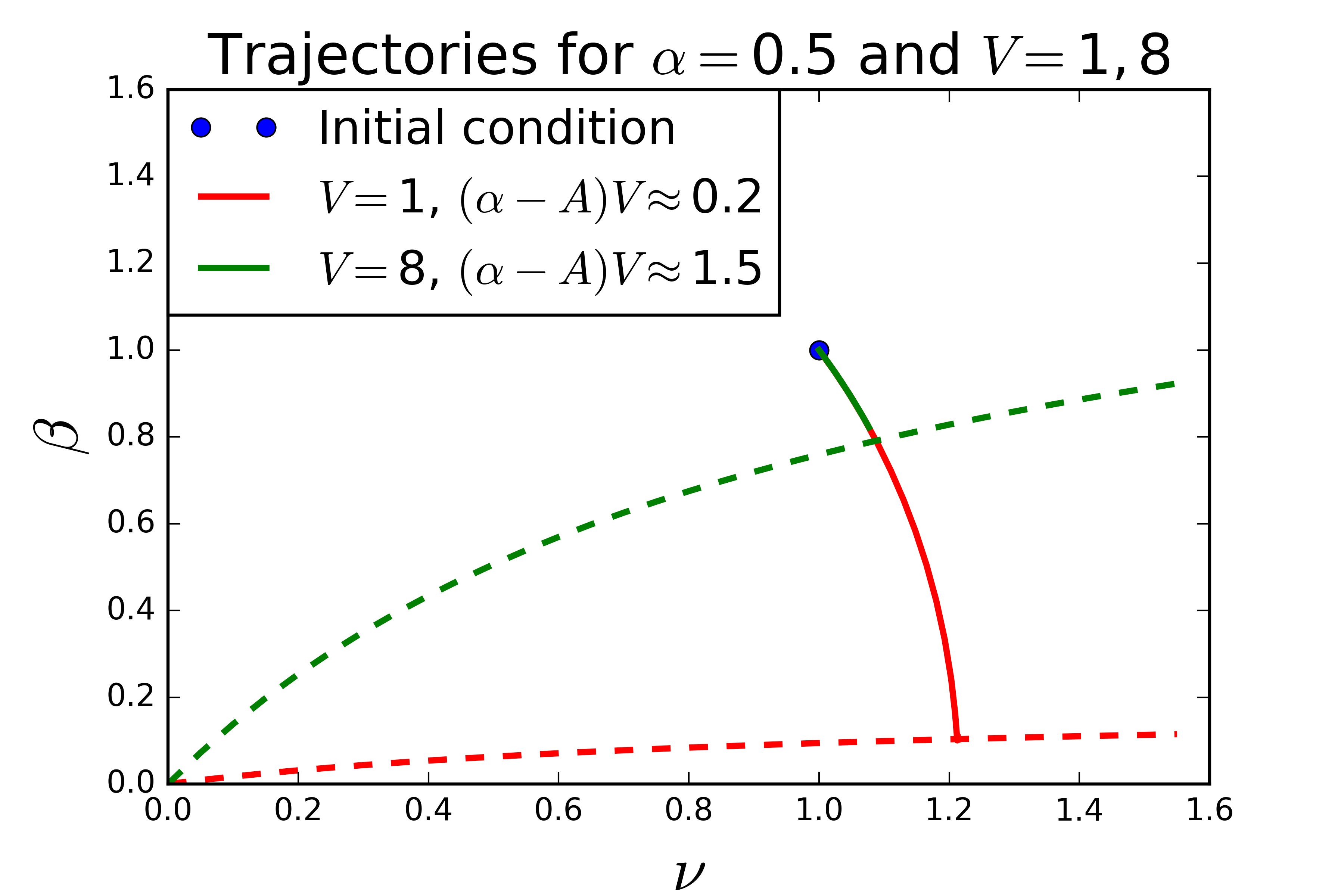}
	\end{minipage}
\hfill
	\begin{minipage}{0.32\textwidth}
       \includegraphics[width=\textwidth]{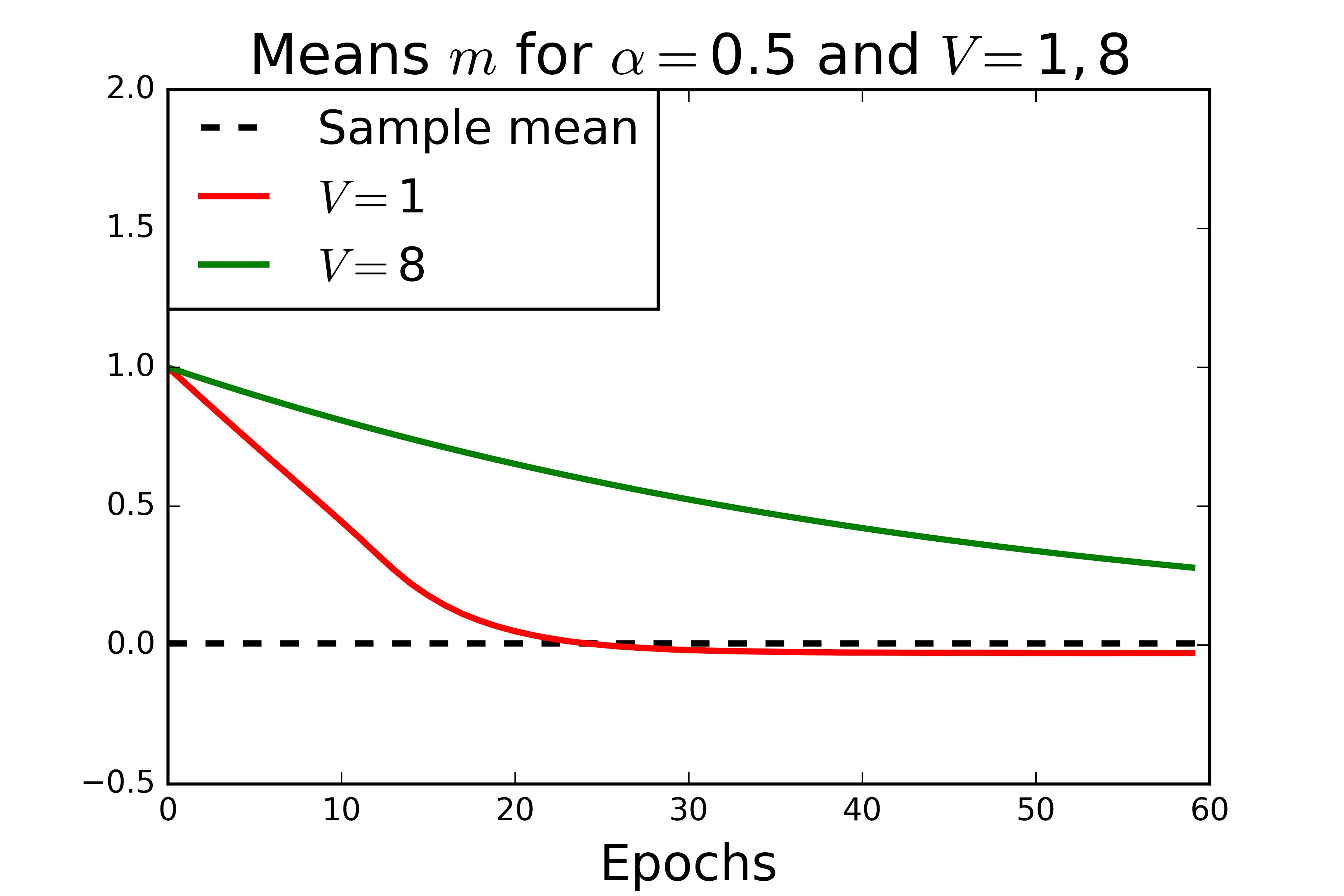}
	\end{minipage}
\hfill
	\begin{minipage}{0.32\textwidth}
       \includegraphics[width=\textwidth]{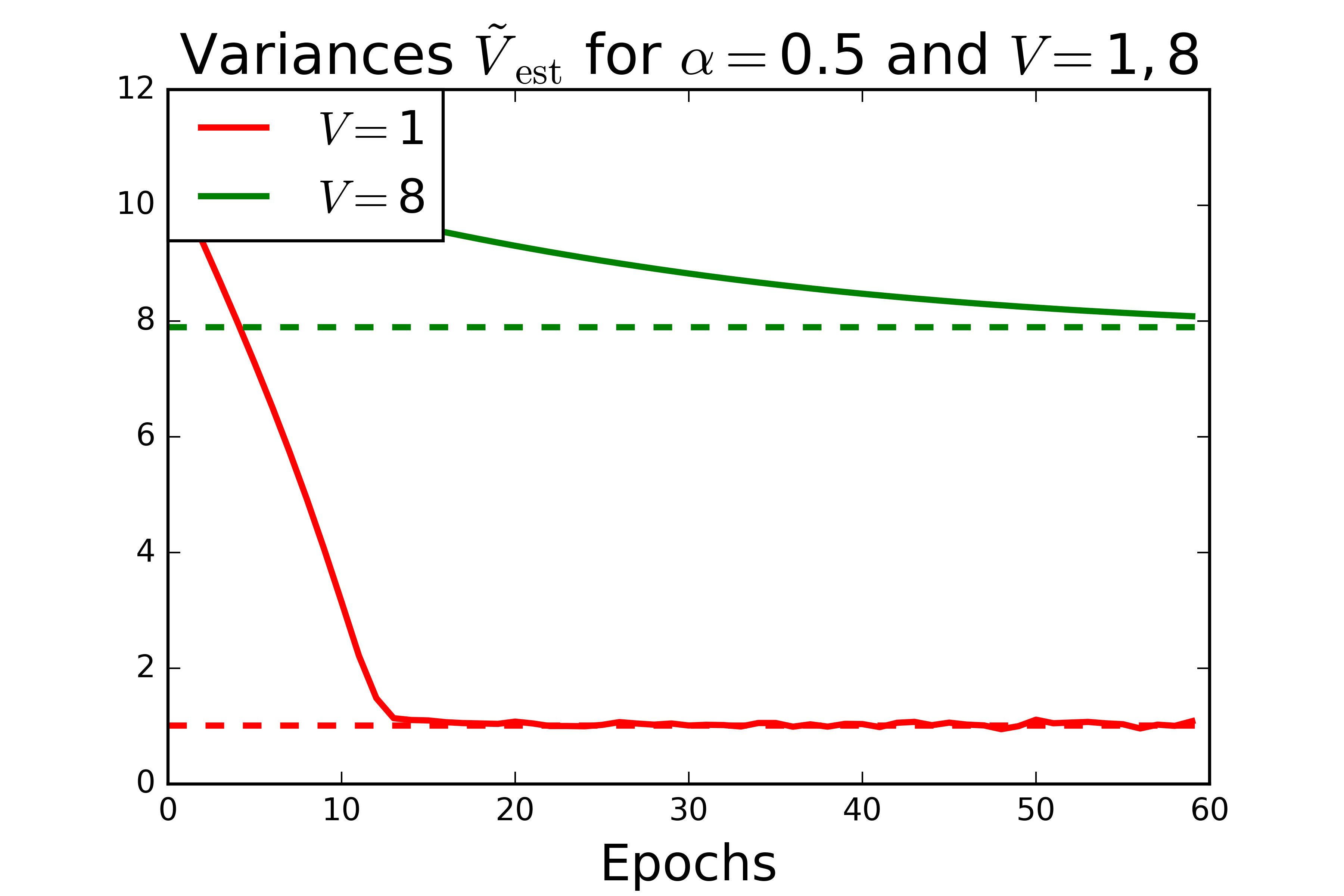}
	\end{minipage}

	\begin{minipage}{0.32\textwidth}
	   \includegraphics[width=\textwidth]{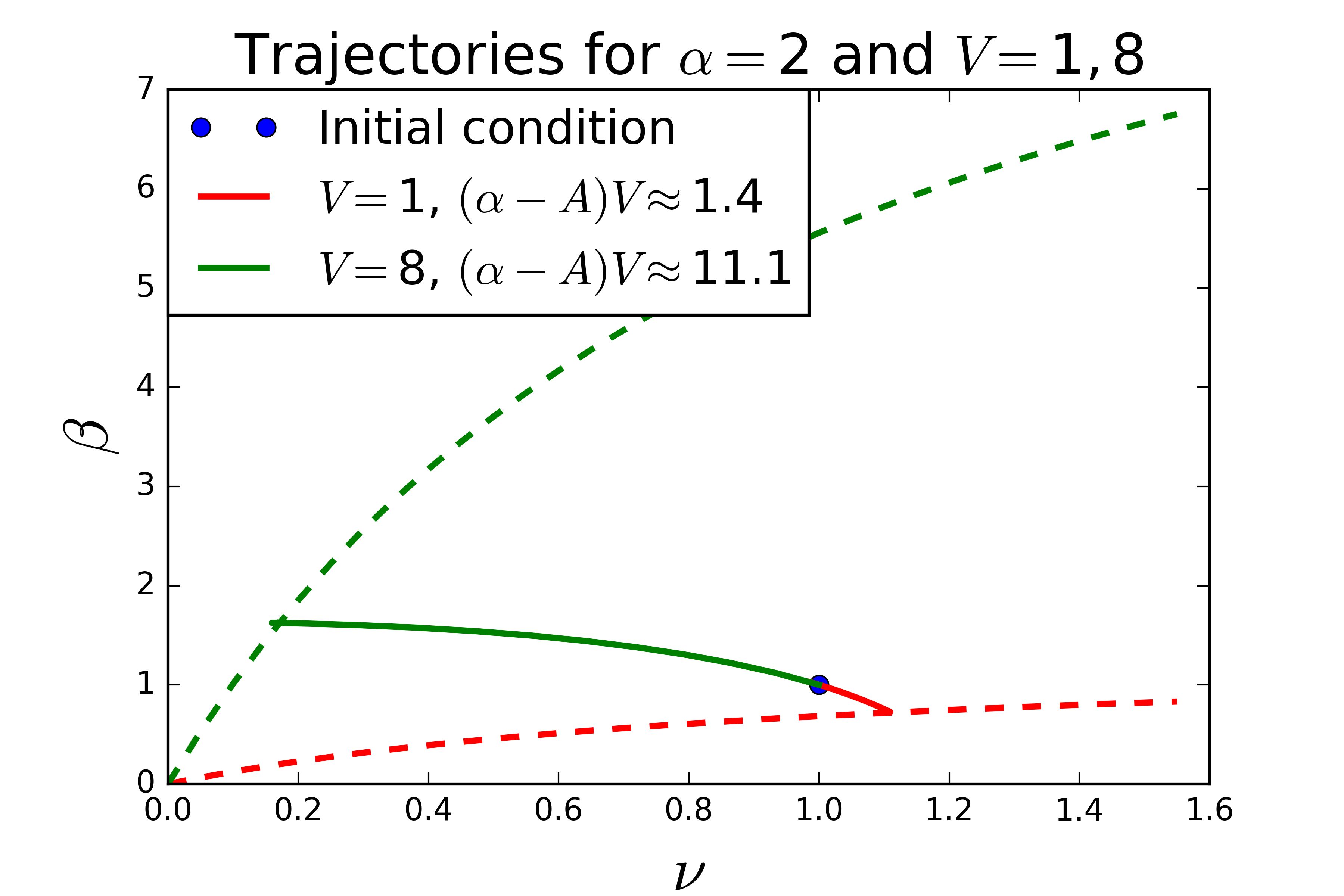}
	\end{minipage}
\hfill
	\begin{minipage}{0.32\textwidth}
       \includegraphics[width=\textwidth]{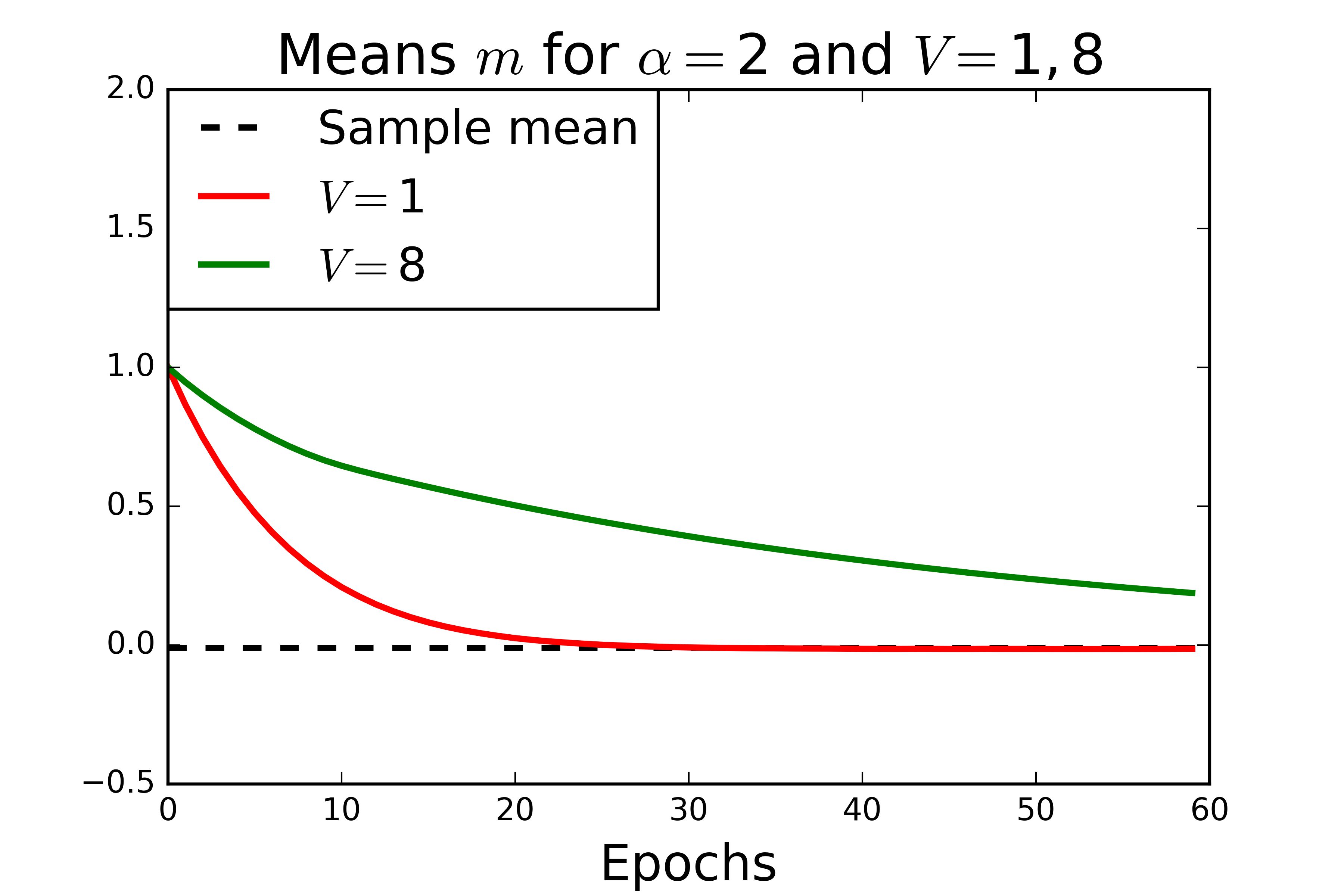}
	\end{minipage}
\hfill
	\begin{minipage}{0.32\textwidth}
       \includegraphics[width=\textwidth]{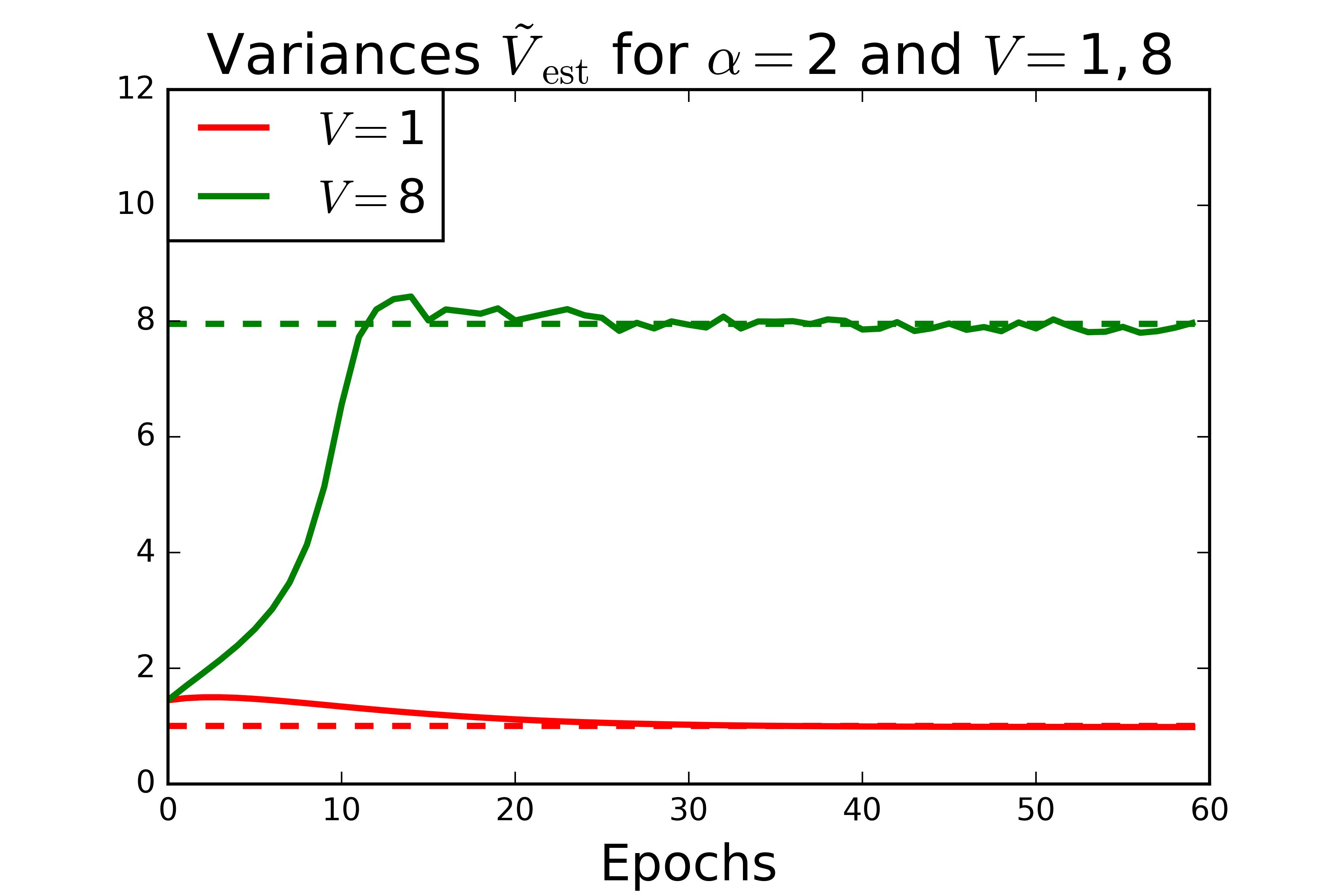}
	\end{minipage}

	\begin{minipage}{0.32\textwidth}
	   \includegraphics[width=\textwidth]{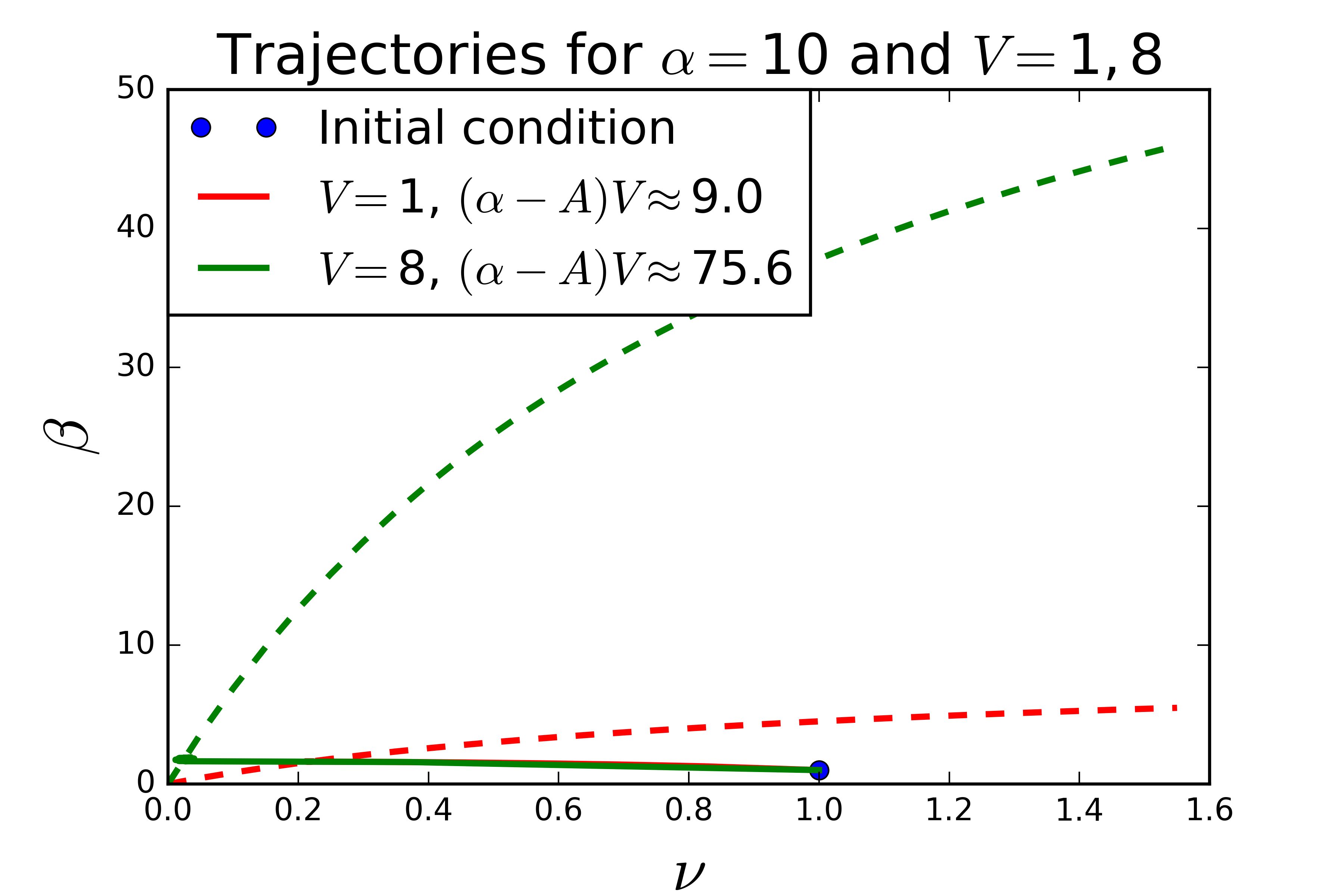}
	\end{minipage}
\hfill
	\begin{minipage}{0.32\textwidth}
       \includegraphics[width=\textwidth]{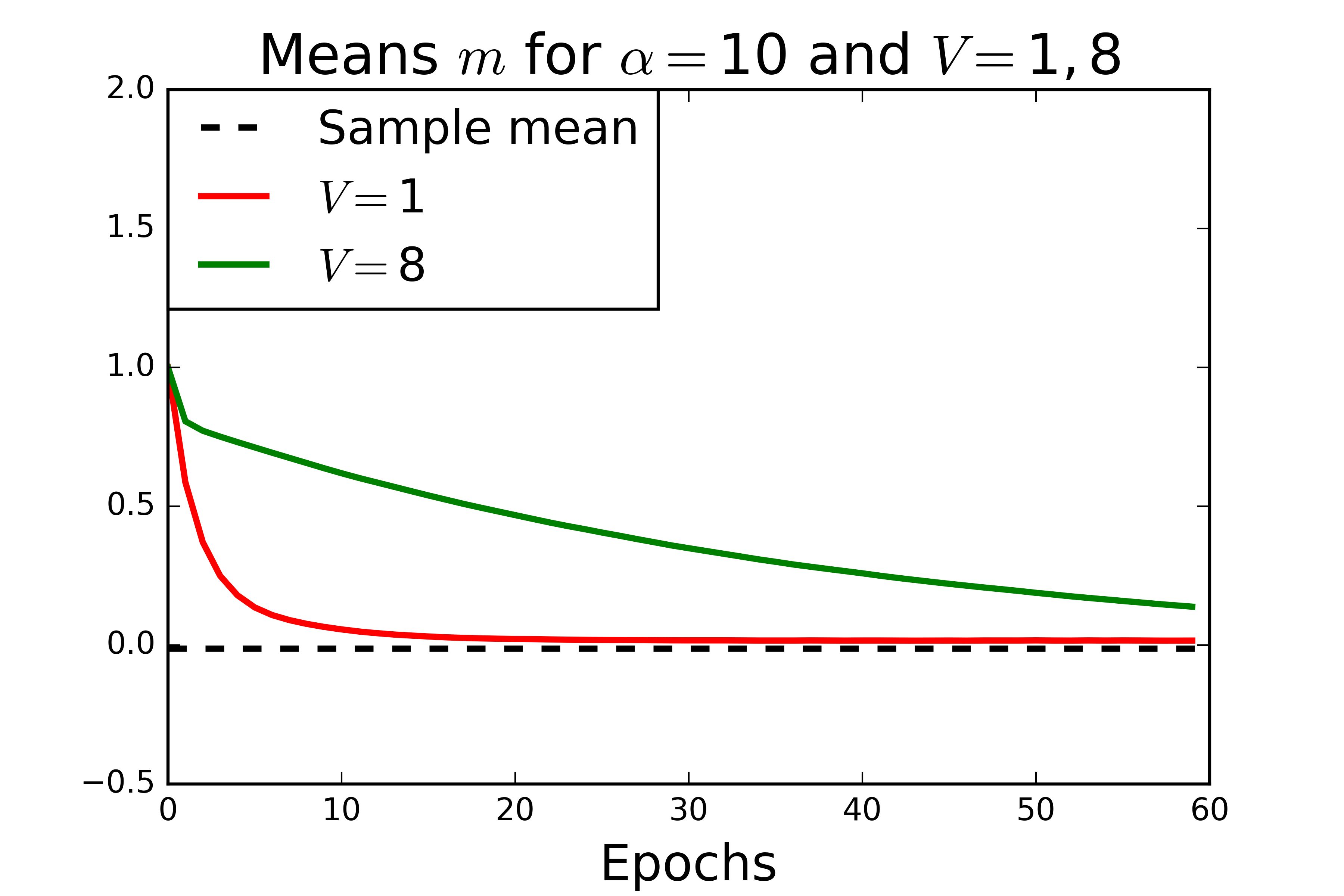}
	\end{minipage}
\hfill
	\begin{minipage}{0.32\textwidth}
       \includegraphics[width=\textwidth]{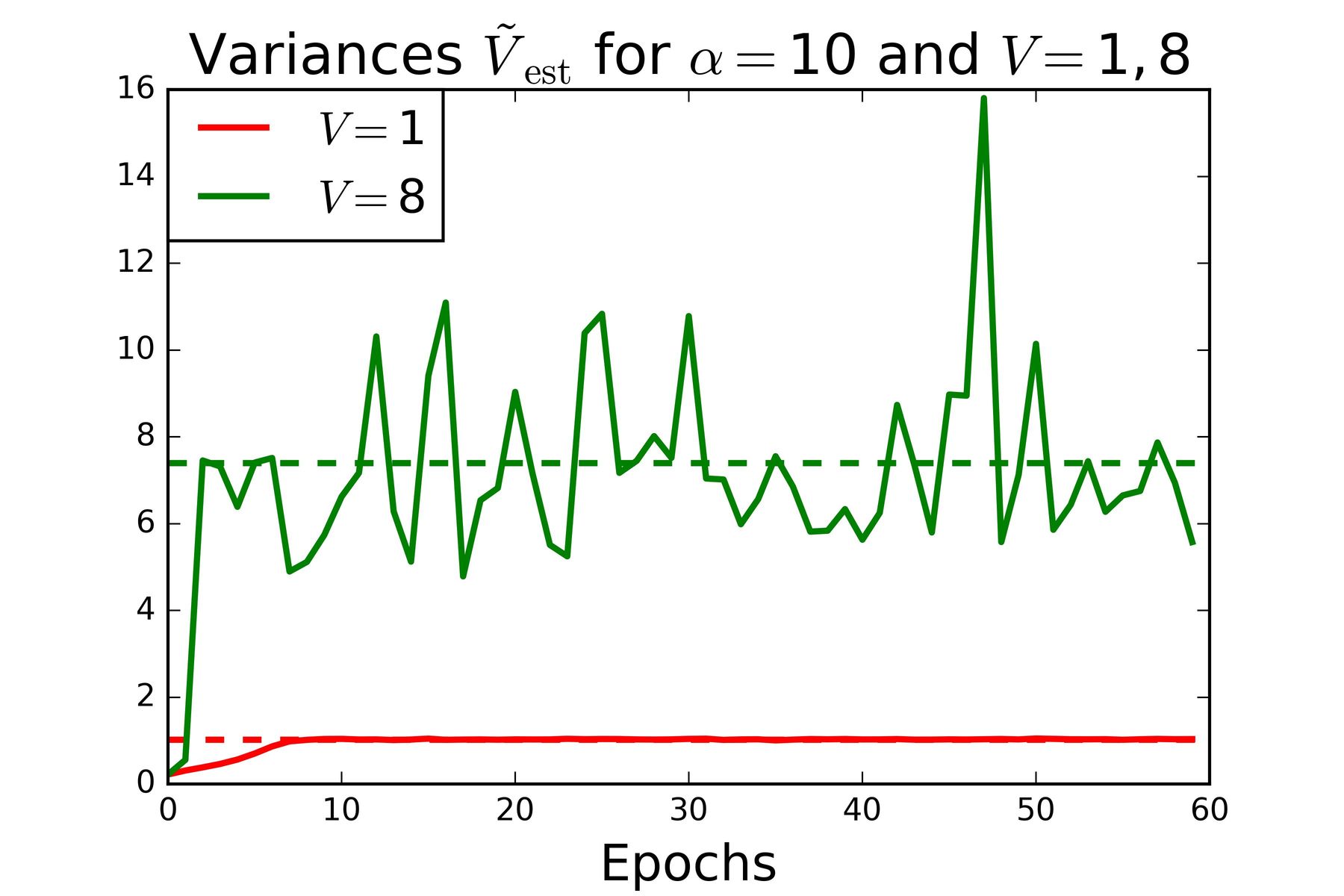}
	\end{minipage}
\caption{Trajectories (left) and graphs of means (middle) and variances (right) versus the number of epochs for different values of $\alpha$ and $V$, based on the GCP update~\eqref{eqGradientDescentFixedAlpha} for 2000 samples drawn from the normal distribution with mean $0$ and variance~$V$. The initial conditions are $m(0)=\beta(0)=\nu(0)=1$. The dashed lines in the left-hand column indicate the corresponding curves of equilibria $\bC_{\alpha,V}$. The dashed lines in the right-hand column indicate the corresponding sample variances.}\label{figCleanNoisy}
\end{figure}

Observations~\ref{obsSmallalpha}--\ref{obsValpha} are illustrated in Fig.~\ref{figCleanNoisy} and explained  in detail below.

\subsubsection{Justification of the observations}

{\bf 1.}
First, we analyze $\LerSp(\Var)$. Consider the limit $V\to 0$ (clean regions). Due to~\eqref{eqAttractingCurveEquil}, the point $(\nu,\beta)$ lies above the line of equilibria $\bC_{\alpha,V}$ at a vertical distance of order $1$ from it. Using~\eqref{eqDKbetaExp} and~\eqref{eqDKnuExp}, we see that
$$
\lim\limits_{V\to 0}\bbE\left[\frac{\partial K}{\partial \beta}\right] = \frac{1}{2\beta},\quad \lim\limits_{V\to 0}\bbE\left[\frac{\partial K}{\partial \nu}\right] = -\frac{1}{2\nu(\nu+1)},
$$
i.e. $(\nu,\beta)$ approaches $\bC_{\alpha,V}$ with speed of order $1$.

Now consider the limit $V\to\infty$. Due to~\eqref{eqAttractingCurveEquil}, the point $(\nu,\beta)$ lies below and to the right from the curve of equilibria $\bC_{\alpha,V}$ at a horizontal distance of order $1$ from it. Using~\eqref{eqDKbetaExp} and~\eqref{eqDKnuExp}, we see that
$$
\lim\limits_{V\to \infty}\bbE\left[\frac{\partial K}{\partial \beta}\right] = -\frac{\alpha}{\beta},\quad \lim\limits_{V\to \infty}\bbE\left[\frac{\partial K}{\partial \nu}\right] = \frac{\alpha}{\nu(\nu+1)},
$$
i.e. $(\nu,\beta)$ approaches $\bC_{\alpha,V}$ with speed of order $\alpha$.

These arguments justify the assertions about $\LerSp(\Var)$ in Observations~\ref{obsSmallalpha} and~\ref{obsLargealpha}.

Further, recall that the trajectory $(\nu,\beta)$ lies on one of the curves~\eqref{eqCurvesNuBeta}.
Thus, if $(\alpha-A(\alpha))V\ll 1$ and $(\nu,\beta)$ approaches $\bC_{\alpha,V}$, the value of $\beta$ will approach $0$, while  $\nu$ will stay of order $1$. On the other hand, if $(\alpha-A(\alpha))V\gg 1$ and $(\nu,\beta)$ approaches $\bC_{\alpha,V}$, the value of $\beta$ will stay of order~$1$ and  $\nu$ will approach $0$. This is illustrated in Fig.~\ref{figCleanNoisy} (left-hand column).
This justifies Observation~\ref{obsValpha}.

{\bf 2.}
Now we analyze $\LerSp(m)$. Here we assume that the {\em variance $V$ has already been estimated approximately}. We express this fact by assuming that the parameters $\beta$ and $\nu$ are such that
\begin{equation}\label{eqKappaV}
\tVest =  c V
\end{equation}
for some $ c\in( c_1, c_2)$, where $ c_2> c_1>0$ do not depend on $V$.
Without loss of generality, assume that $\bbE[\by]=0$ and $m>0$.  Then, due to~\eqref{eqDKmu},
\begin{equation}\label{eqDkmuVAsymp}
\begin{aligned}
\bbE\left[\frac{\partial K}{\partial m}\right] &= (\alpha+1/2)\frac{1}{(2\pi)^{1/2}} \int\limits_{-\infty}^\infty \frac{m-V^{1/2}z}{ c(\alpha-A)V + \frac{\left(m-V^{1/2} z\right)^2}{2}} e^{-\frac{z^2}{2}}  dz \\
& =
(\alpha+1/2)\frac{1}{(2\pi V)^{1/2}} \int\limits_{-\infty}^\infty \frac{\frac{m}{V^{1/2}}-z}{ c(\alpha-A) + \frac{\left(\frac{m}{V^{1/2}}-z\right)^2}{2}} e^{-\frac{z^2}{2}}  dz.
\end{aligned}
\end{equation}
Hence,
\begin{equation}\label{eqDkmuVAsymp1}
\begin{aligned}
& \bbE\left[\frac{\partial K}{\partial m}\right] = \frac{2\alpha+1}{m} + o(V)\quad\text{as } V\to 0,\\
& \bbE\left[\frac{\partial K}{\partial m}\right] = \frac{k(\alpha)}{V} + o\left(\frac{1}{V}\right)\quad\text{as } V\to \infty,
\end{aligned}
\end{equation}
where
$$
k(\alpha)= \frac{(\alpha+1/2)m}{(2\pi)^{1/2}}\int_{-\infty}^{\infty}\frac{ c(\alpha-A)-z^2/2}{( c(\alpha-A)+z^2/2)^2}e^{-\frac{z^2}{2}}dz
.$$
The constant $k(\alpha)$ can be obtained by dividing the integral in the right-hand side of~\eqref{eqDkmuVAsymp} by $V^{-1/2}$ and applying  L'Hospital's rule. Using the properties of $A(\alpha)$ in~\eqref{eqPropertiesA}, one can show that $k(\alpha)$ is positive, bounded, and bounded away from zero for all $\alpha>0$.

These arguments justify the assertions about $\LerSp(m)$ in Observations~\ref{obsSmallalpha} and~\ref{obsLargealpha}.

\section{GCP neural networks: experiments}\label{secNN}

\subsection{Methods}\label{subsecMethods}

We compare the  following methods:
\begin{enumerate}

\item the GCP method with and without the correction formula   for $\tVest$ in~\eqref{eqEstMeanVarNetwork}. Whenever we apply the correction formula, we indicate this by writing GCP$_{\rm corr}$;

  \item the maximum likelihood method (ML), in which one maximizes the likelihood of the normal distribution with the input-dependent mean and precision,

  \item density power divergence method (DPD)~\cite{Basu1998,Gosh2016}, in which one minimizes the density power (instead of the KL) divergence from the ground truth distribution to the approximating normal distribution; as the GCP, this method is known to be robust against outliers,

  \item Bayesian maximum likelihood method (ML$_{\rm Bayes}$)~\cite{KendallGal17}, in which one maximizes the likelihood of the normal distribution with input-dependent precision, using dropout for both training and prediction to approximate the posterior distribution of the weights.


  \item Stein variational gradient descent (SVGD)~\cite{LiuWang2016}; Bayesian method, in which we use a Gaussian likelihood with input-dependent mean and input-independent variance and a particle approximation of the posterior distribution of the weights,

  \item the probabilistic back propagation (PBP)~\cite{HernandezLobato15}; Bayesian method, in which the posterior of weights is approximated with assumed density filtering~\cite{Opper98} and expectation-propagation~\cite{Minka01} methods. The variance of the Gaussian likelihood is assumed input-independent.

\end{enumerate}

Note that the GCP (GCP$_{\rm corr}$), ML, and DPD estimate aleatoric heteroscedastic uncertainty; the Bayesian methods SVGD and PBP estimate aleatoric homoscedastic and epistemic heteroscedastic uncertainty; and the ML$_{\rm Bayes}$ estimates both uncertainties as heteroscedastic.

\subsection{Synthetic data set: aleatoric vs. epistemic uncertainty}\label{subsecSyntheticAleatoricEpistemic}

\begin{figure}[!t]
	\begin{minipage}{0.32\textwidth}
	   \includegraphics[width=\textwidth]{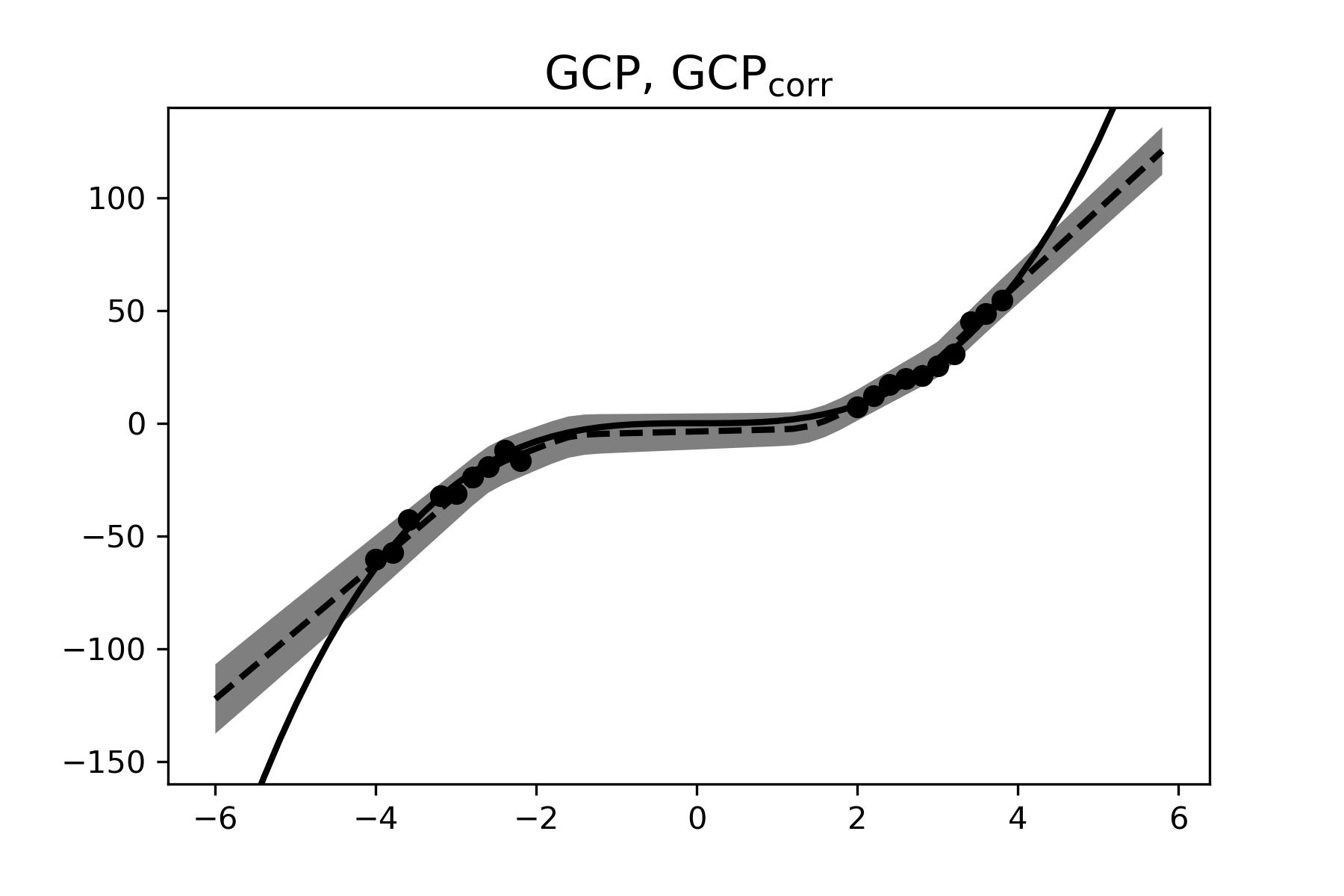}
	\end{minipage}
\hfill
	\begin{minipage}{0.32\textwidth}
       \includegraphics[width=\textwidth]{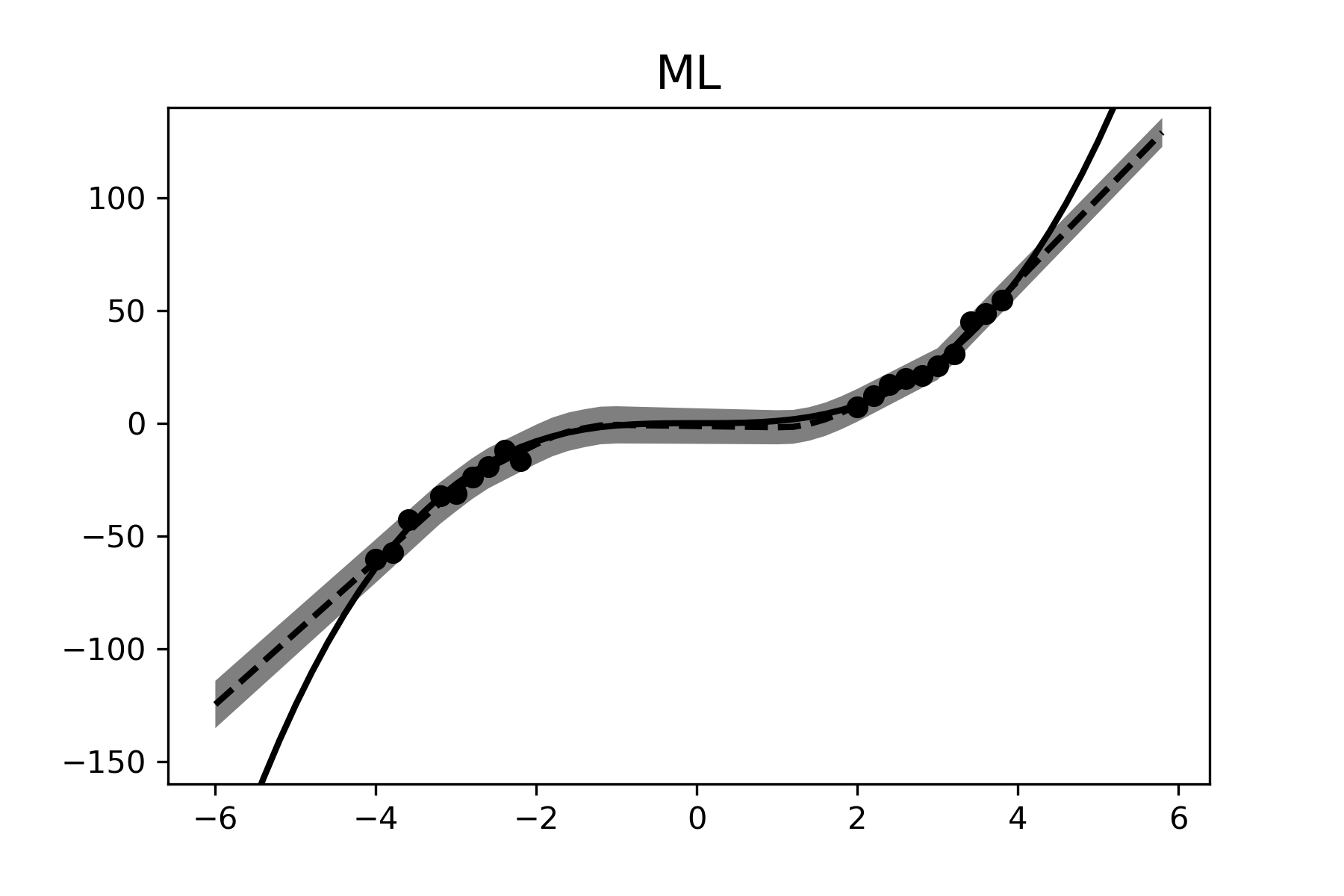}
	\end{minipage}
\hfill
	\begin{minipage}{0.32\textwidth}
       \includegraphics[width=\textwidth]{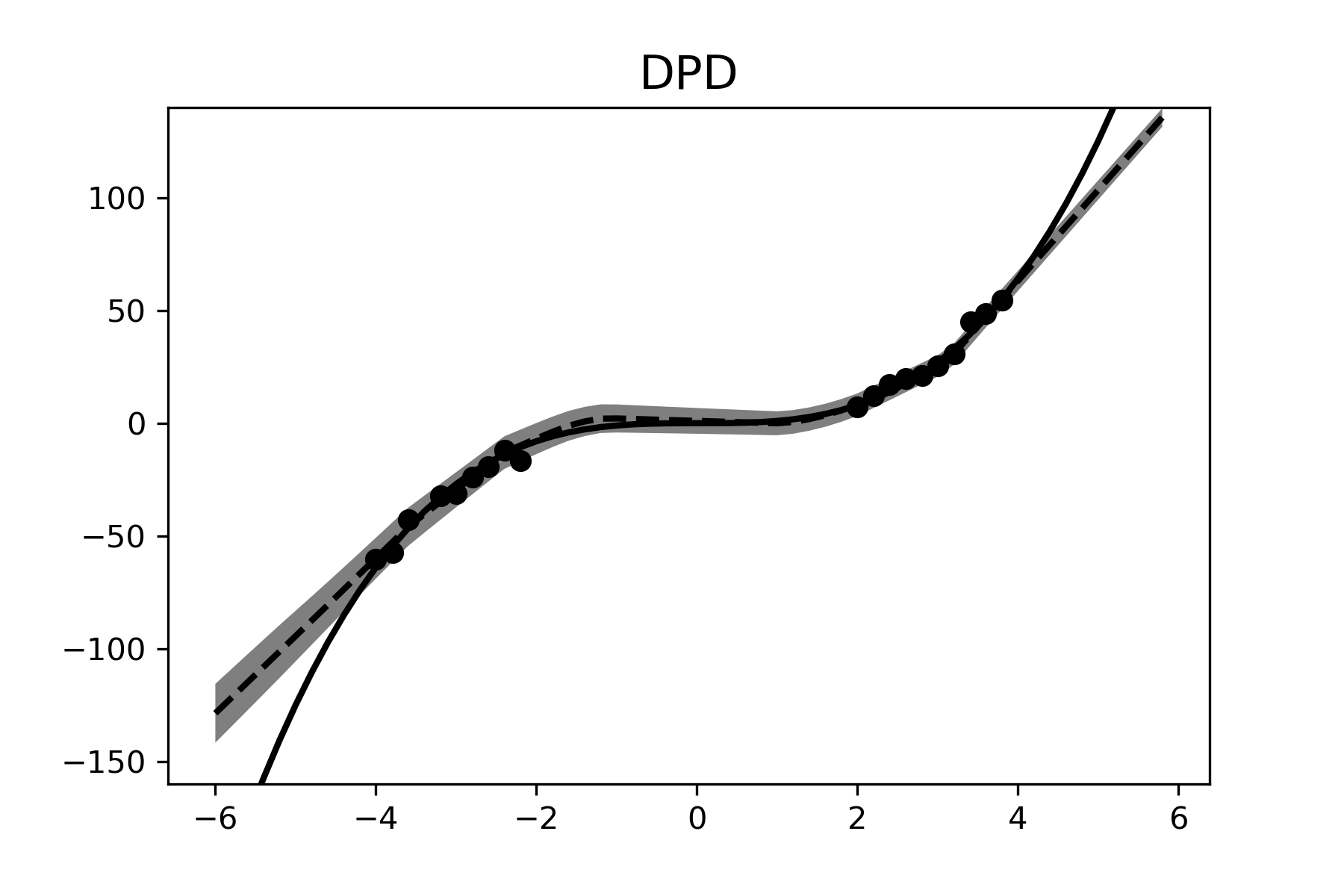}
	\end{minipage}

	\begin{minipage}{0.32\textwidth}
       \includegraphics[width=\textwidth]{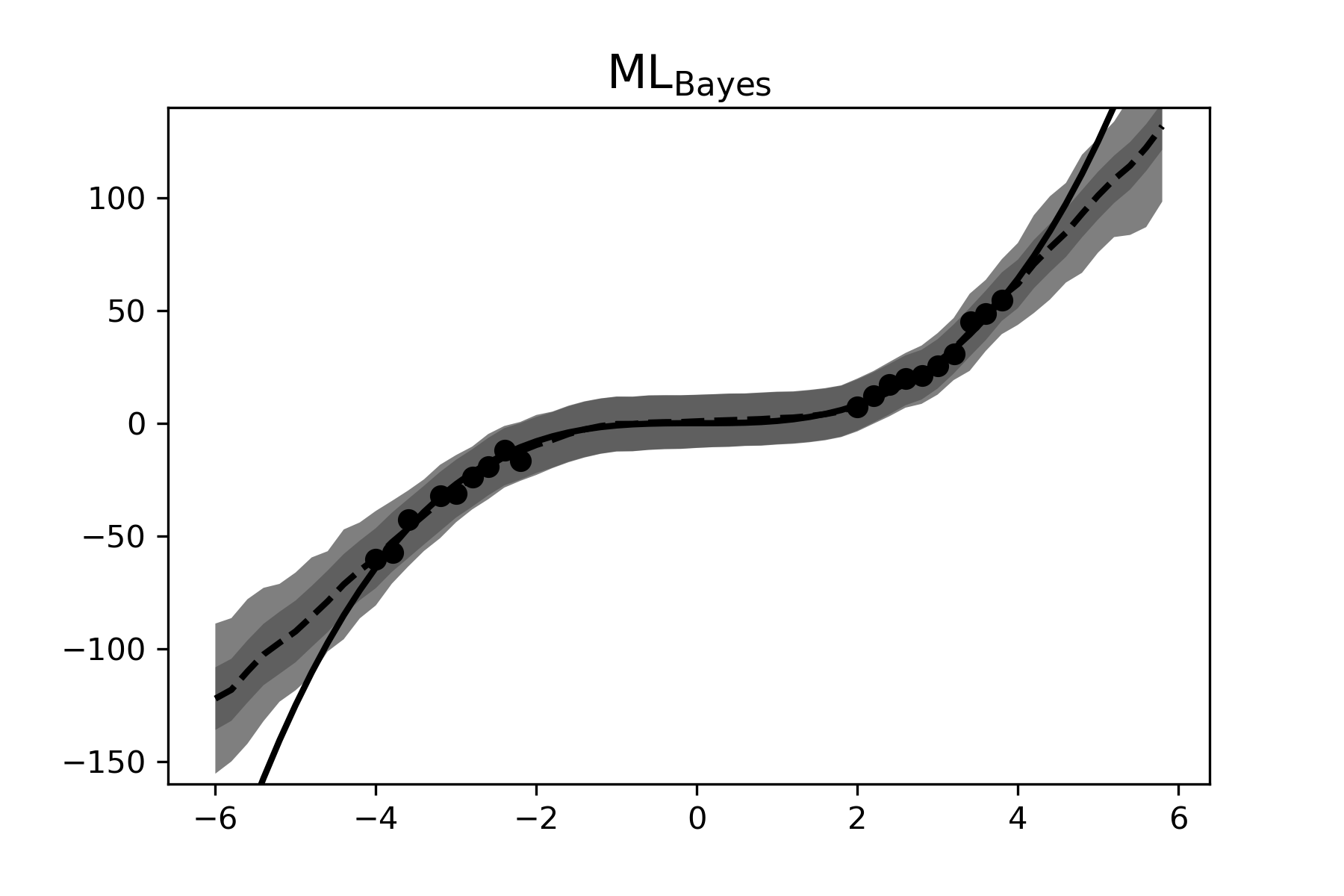}
	\end{minipage}
\hfill
	\begin{minipage}{0.32\textwidth}
       \includegraphics[width=\textwidth]{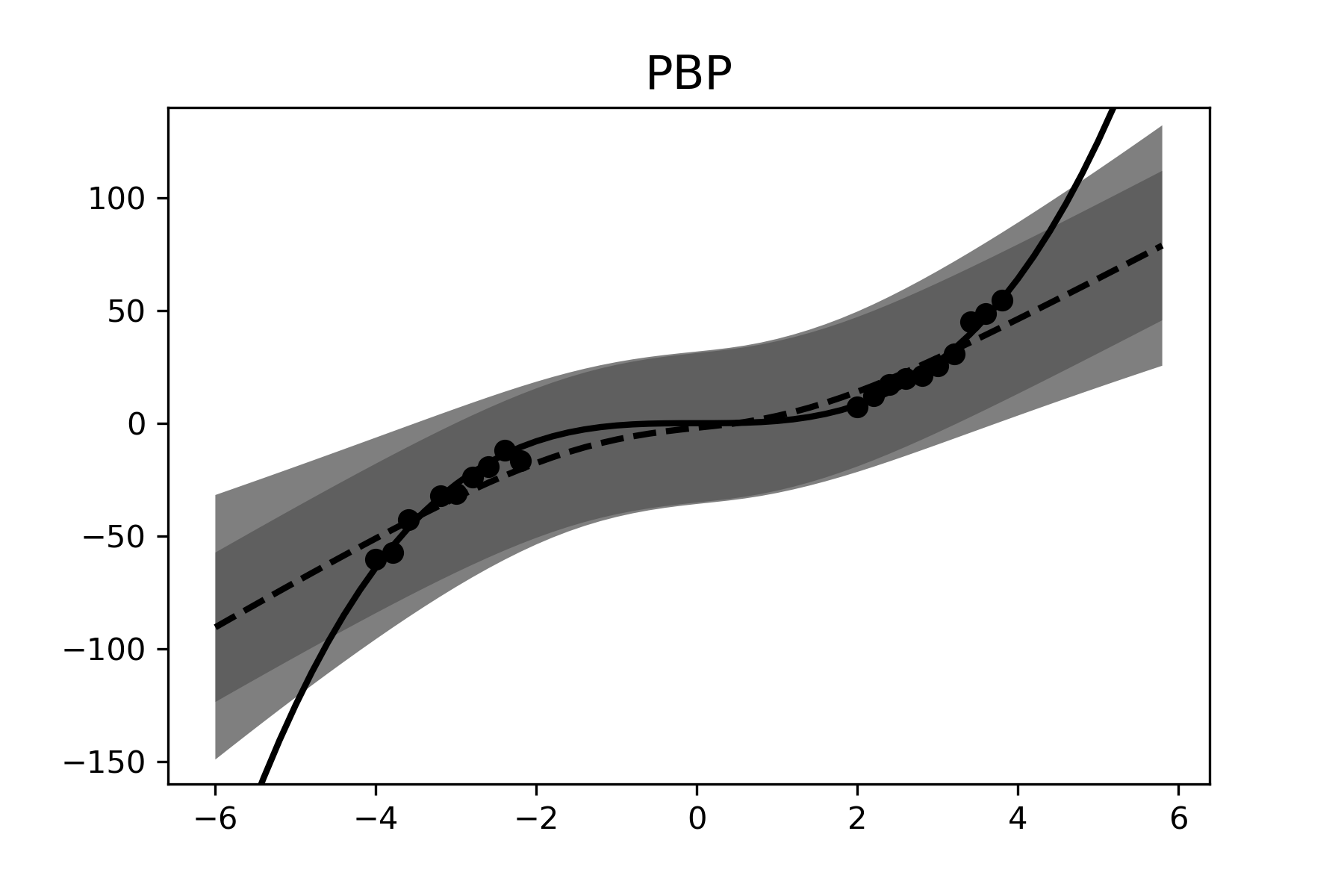}
	\end{minipage}
\hfill
	\begin{minipage}{0.32\textwidth}
       \includegraphics[width=\textwidth]{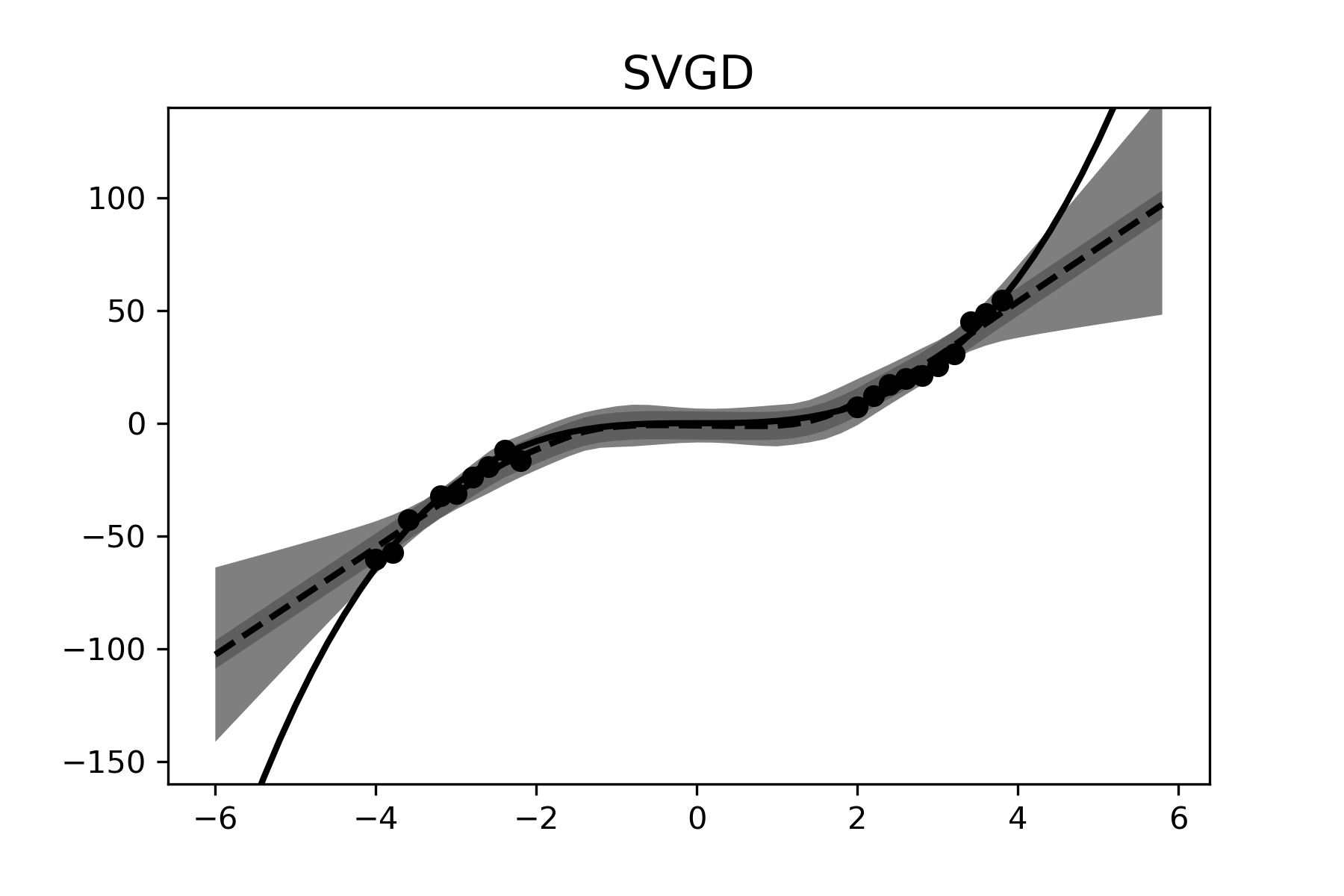}
	\end{minipage}
\caption{Predictions made by different methods on the synthetic data set $\by\sim\cN(x^3,3^2)$. The samples are shown as filled circles, the ground truth as a solid line, the predicted mean as a dashed line. For the Bayesian methods ML$_{\rm Bayes}$, PBP, and SVGD, the darker shaded areas correspond to $\pm 3$ aleatoric standard deviations and the lighter shaded areas correspond to $\pm 3$ standard deviations of the predictive distribution. For the non-Bayesian methods GCP (GCP$_{\rm corr}$), ML, and DPD, the aleatoric standard deviation coincides with the predictive one (both indicated as the lighter shaded areas). In the absence of outliers, GCP and GCP$_{\rm corr}$ predict almost identical variance due to large values of $\alpha$.}\label{figSyntheticDataAleatoricEpistemic}
\end{figure}

We evaluate the predictive distributions in the above methods for a synthetic data set. The input data set $X$ consists of 20 points uniformly distributed on $[-4,-2]\cup [2,4]$. For each $x\in X$, we sample $y$ from the normal distribution with mean $x^3$ and standard deviation~$3$. In this and next subsections, we use one-hidden layer networks with 100 hidden units. Figure~\ref{figSyntheticDataAleatoricEpistemic} confirms that the GCP (GCP$_{\rm corr}$), ML, and DPD  capture only aleatoric uncertainty, while the Bayesian methods ML$_{\rm Bayes}$, PBP, and SVGD also capture epistemic uncertainty due to the lack of data. To visualize the two types of uncertainty, we follow~\cite{Depeweg2018} and use the law of total variance $\bbV[\by|x]=\bbE_w[\bbV[\by|x,w]] + \bbV_w[\bbE[\by|x,w]]$, where the first term corresponds to aleatoric variance and the second term to epistemic.

\subsection{Synthetic data set: robustness to outliers in the training set}\label{subsecSynthetic}

\begin{figure}[!t]
	\begin{minipage}{0.3\textwidth}
	   \includegraphics[width=\textwidth]{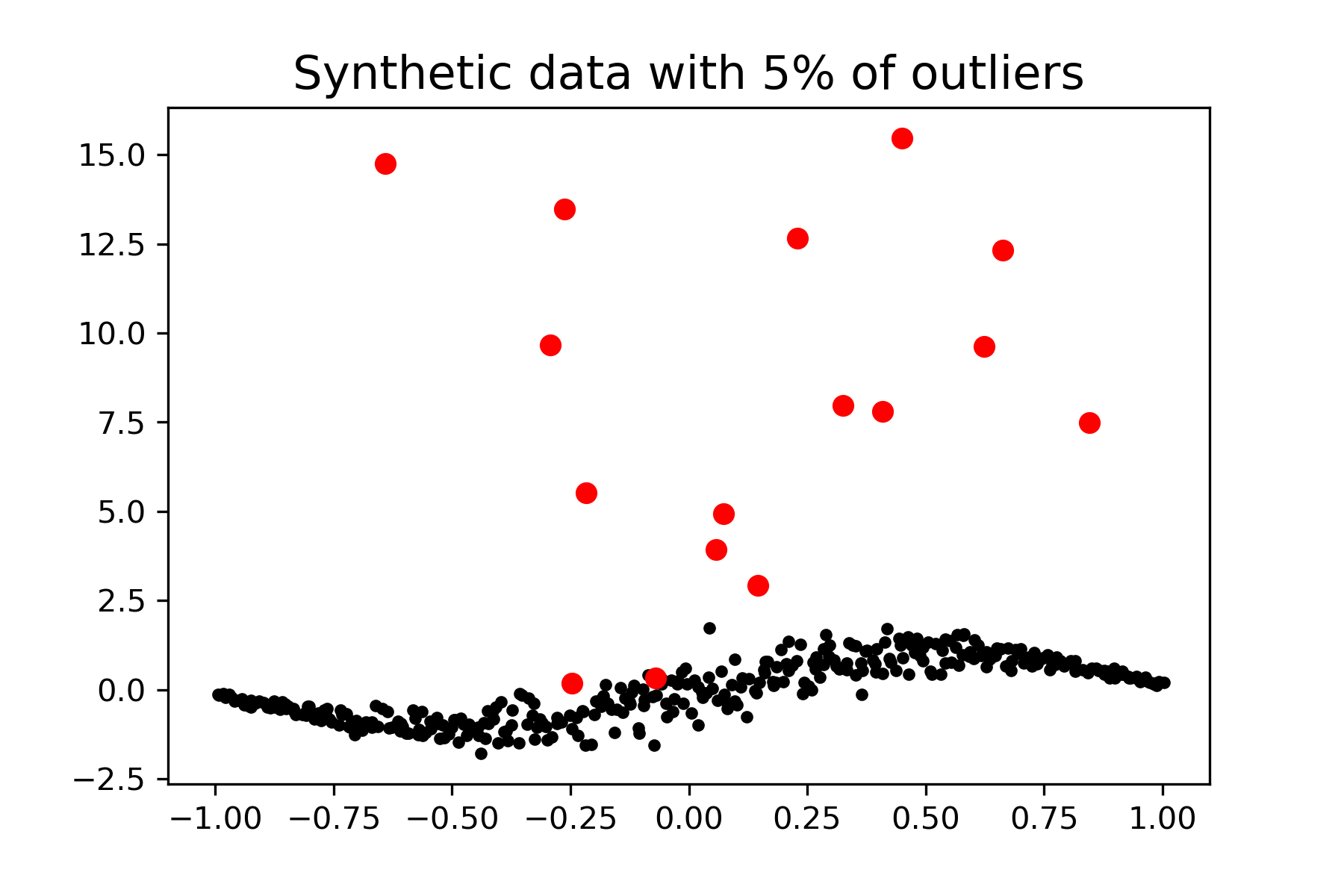}
	\end{minipage}
\hfill
	\begin{minipage}{0.3\textwidth}
       \includegraphics[width=\textwidth]{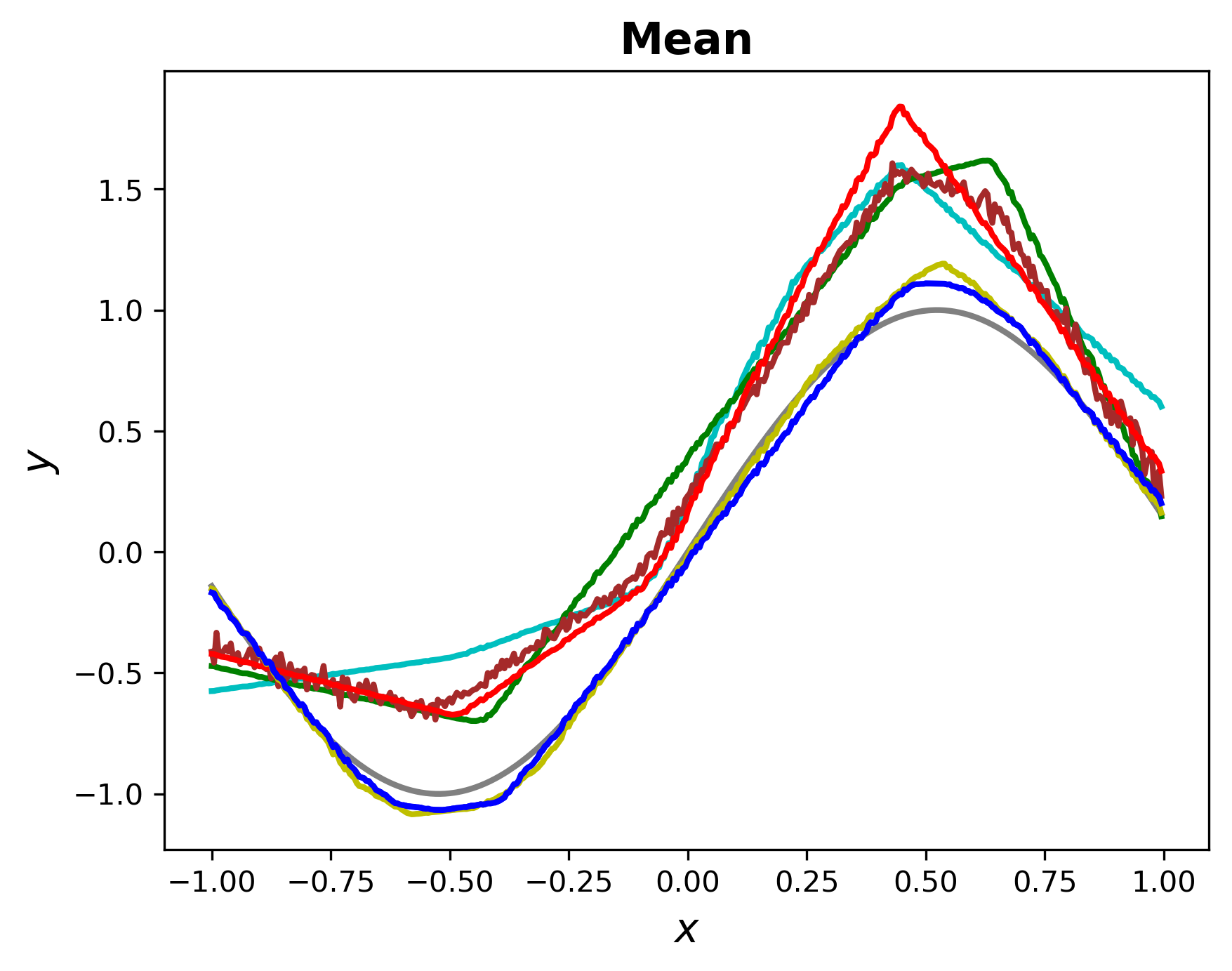}
	\end{minipage}
\hfill
	\begin{minipage}{0.38\textwidth}
       \includegraphics[width=\textwidth]{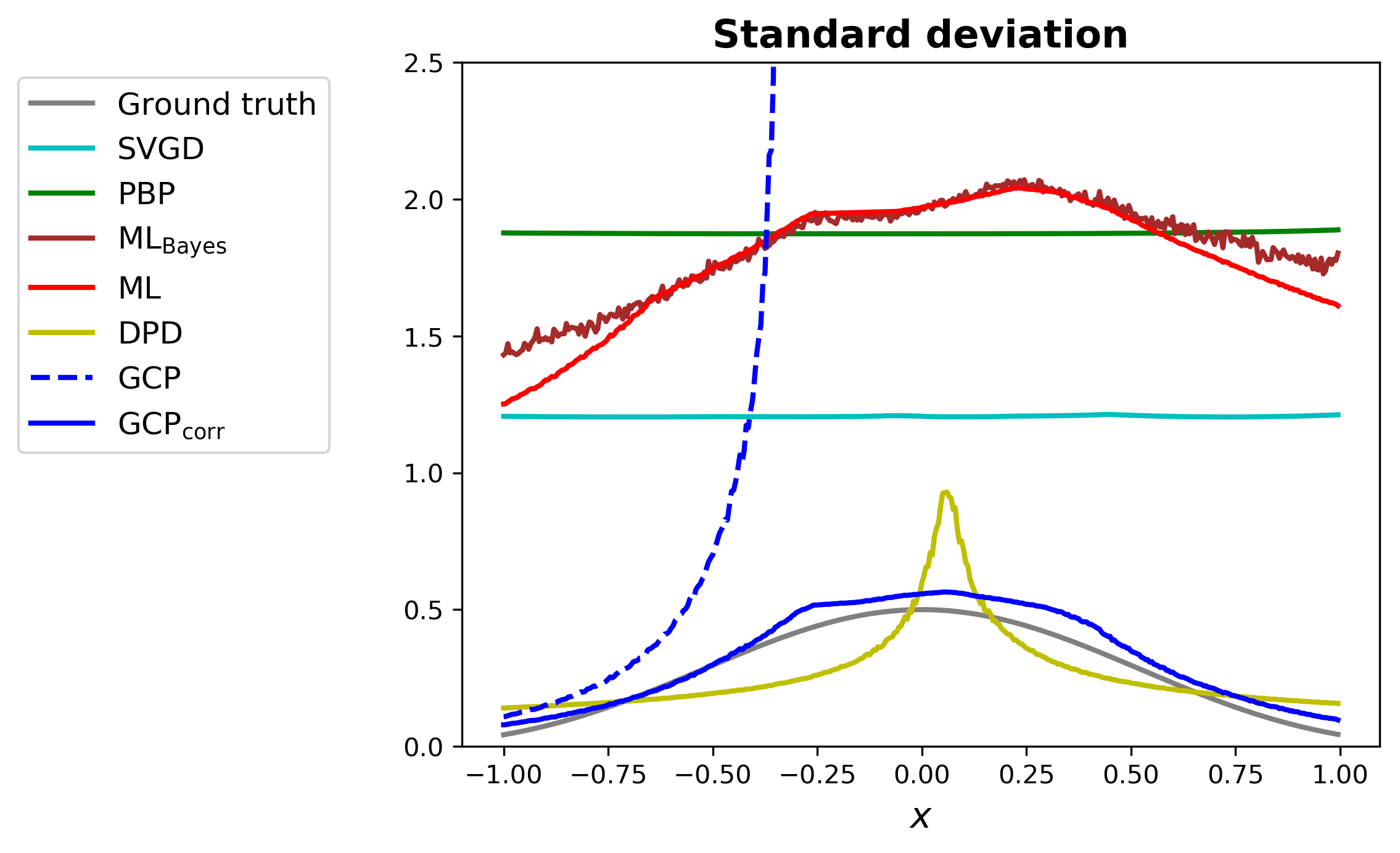}
	\end{minipage}
\caption{Left: Synthetic data with mean $\sin(3x)$ and standard deviation $0.5\cos^4 x$ (black dots), complemented by 5\% of outliers sampled from a uniform distribution on the interval $(-4,16)$ (red disks). Middle: The ground truth mean $\sin(3x)$ and the means predicted by the different methods. Right: The ground truth standard deviation $0.5\cos^4 x$ and the standard deviations predicted by the different methods.}\label{figSyntheticDataWithOutliers}
\end{figure}

We generate a synthetic data set containing 5\% of outliers. The set $X$ consists of 400 points uniformly distributed on the interval $(-1,1)$. For each $x\in X$, with probability 0.95 we sample $y$ from the normal distribution with mean $\sin(3x)$ and standard deviation $0.5\cos^4 x$, and with probability 0.05 we sample $y$ from a uniform distribution on the interval $(-4,16)$. Figure~\ref{figSyntheticDataWithOutliers} shows the data and the fits of different methods. 
The means predicted by the GCP and DPD
  are significantly less affected by the outliers compared with the other methods. However, the standard deviations\footnote{For Bayesian methods, we plot the standard deviations of the predictive distributions. The epistemic standard deviations are negligible in this example due to a relatively large amount of data.} predicted by all the methods, except for the GCP$_{\rm corr}$, are significantly distorted. Although the DPD is known to be robust against the outliers, it does not manage to properly capture the $x$-dependence of the variance. The ML, ML$_{\rm Bayes}$, SVGD, PBP, and the GCP without the correction overestimate the variance. However, using the correction formula for $\tVest$ in~\eqref{eqEstMeanVarNetwork} allows GCP$_{\rm corr}$ to reconstruct the ground truth variance. Furthermore, the knowledge of $\alpha$ in the GCP provides additional information, namely, small values of $\alpha$ indicate that the corresponding samples belong to a (less trust-worthy) region in which the training set contained outliers.



\subsection{Real world data sets}\label{subsecRealWorld}

\textbf{Architectures.} We use one-hidden layer networks for the parameters of the prior~\eqref{eqNetworksMAlphaBetaNu}  with ReLU
nonlinearities. Each network contains 50 hidden units for all the data sets below, except for the largest MSD set. For the latter, we use 100 hidden units. For regularization, we use a dropout layer between the hidden layer and the output unit. Our approach is directly applicable to neural networks of any depth and structure, however we kept one hidden layer for the compatibility of our validation with~\cite{HernandezLobato15,Lakshminarayanan16,Lakshminarayanan17,GurHannesLU}.

\textbf{Measures.} We use two measures to estimate the quality of the fit.
\begin{enumerate}
  \item The overall root mean squared error ({\em RMSE}).
  \item The area under the following curve ({\em AUC}), measuring the trade-off between properly learning the mean and the variance. Assume the test set contains $N$ samples. We order them with respect to their predicted variance. For each $n=0,\dots,N-1$, we remove $n$  samples with the highest variance and calculate the RMSE for the remaining $N-n$ samples (with the lowest variance). We denote it by ${\rm RMSE}(n)$ and plot it versus~$n$ as a continuous piecewise linear curve. The second measure is the area under this curve normalized by $N-1$:
      $$
      {\rm AUC} := \frac{1}{N-1}\sum\limits_{n=0}^{N-2}\frac{{\rm RMSE}(n)+{\rm RMSE}(n+1)}{2}.
      $$

\end{enumerate}

\textbf{Data sets.}
We analyze the following publicly available data sets: Boston House Prices~\cite{Harrison78} ($506$ samples, 13 features), Concrete Compressive Strength~\cite{ICheng98} ($1030$ samples, 8 features), Combined Cycle Power Plant~\cite{Tufekci14,Kaya12} ($9568$ samples, 4 features), Yacht Hydrodynamics~\cite{Gerritsma81,Ortigosa07} (308 samples, 6 features), Kinematics of an 8 Link Robot Arm Kin8Nm\footnote{http://mldata.org/repository/data/viewslug/regression-datasets-kin8nm/} (8192 samples, 8 feature), and Year Prediction MSD~\cite{Lichman13} ($515345$ samples, 90 features). For each data set, a one-dimensional target variable is  predicted.
Each data set, except for the year prediction MSD, is randomly split into 50 train-test folds with 95\% of samples in each train subset. All the measure values reported below are the averages of the respective measure values over 50 folds. For the year prediction MSD, we used a single split recommended in~\cite{Lichman13}.

\textbf{Results.}
Table~\ref{tableSqErrorBostonConcrete}
shows the measure values of the different methods. In each column,
we mark a method in bold if it is significantly (due to
the two-tailed paired difference test with p = 0.05) better or
indistinguishable from all the other methods. We do not present the GCP$_{\rm corr}$ results in this table because, in the absence of outliers, it yielded AUC very close to that of GCP. We see that the GCP achieves the best AUC values on all the data sets (except MSD), which indicates the best trade-off between properly learning the mean and the variance.

\begin{table}
\resizebox{\textwidth}{!}{%
\begin{tabular}{lcc}
             &         {\bf Boston } &                    \\
\midrule
{} &                  RMSE &                  AUC \\
\midrule
SVGD               &  {\bf  2.93$\pm$1.02} &        2.21$\pm$0.61 \\
PBP                &  {\bf  2.99$\pm$0.97} &        2.10$\pm$0.47 \\
ML$_{\rm Bayes}$   &         3.55$\pm$1.41 &        1.92$\pm$0.39 \\
ML                 &         3.63$\pm$1.49 &        2.00$\pm$0.48 \\
DPD                &         3.91$\pm$1.74 &        2.21$\pm$0.63 \\
GCP                &         3.74$\pm$1.59 &  {\bf 1.79$\pm$0.40} \\
\bottomrule
\end{tabular}
\
\begin{tabular}{lcc}
      {\bf Concrete } &                    \\
\midrule
                 RMSE &                  AUC \\
\midrule
 {\bf  5.15$\pm$0.65} &        4.21$\pm$0.72 \\
        5.37$\pm$0.60 &        4.45$\pm$0.56 \\
        5.45$\pm$0.63 &        3.75$\pm$0.52 \\
 {\bf  5.29$\pm$0.81} &  {\bf 3.51$\pm$0.59} \\
        6.05$\pm$0.86 &        4.23$\pm$0.83 \\
        5.47$\pm$0.65 &  {\bf 3.53$\pm$0.65} \\
\bottomrule
\end{tabular}
\
\begin{tabular}{lcc}
      {\bf Power } &                    \\
\midrule
                 RMSE &                  AUC \\
\midrule
 {\bf  4.13$\pm$0.32} &        4.10$\pm$0.28 \\
 {\bf  4.11$\pm$0.30} &        3.85$\pm$0.23 \\
        4.19$\pm$0.30 &        3.72$\pm$0.21 \\
        4.15$\pm$0.30 &        3.70$\pm$0.24 \\
 {\bf  4.12$\pm$0.32} &        3.66$\pm$0.36 \\
        4.13$\pm$0.31 &  {\bf 3.54$\pm$0.33} \\
\bottomrule
\end{tabular}
\
\begin{tabular}{lcc}
      {\bf Yacht } &                    \\
\midrule
                 RMSE &                  AUC \\
\midrule
        0.88$\pm$0.40 &        0.37$\pm$0.18 \\
        1.04$\pm$0.39 &        0.57$\pm$0.13 \\
 {\bf  0.78$\pm$0.39} &  {\bf 0.24$\pm$0.07} \\
 {\bf  0.82$\pm$0.41} &  {\bf 0.25$\pm$0.09} \\
        2.48$\pm$1.24 &        0.36$\pm$0.16 \\
        0.96$\pm$0.49 &  {\bf 0.23$\pm$0.09} \\
\bottomrule
\end{tabular}
\
\begin{tabular}{lcc}
      {\bf Kin8nm } &                    \\
\midrule
                 RMSE &                  AUC \\
\midrule
 {\bf 0.09$\pm$0.01} &         0.08$\pm$0.01 \\
       0.10$\pm$0.00 &         0.08$\pm$0.00 \\
       0.11$\pm$0.01 &         0.07$\pm$0.00 \\
       0.10$\pm$0.01 &  {\bf  0.06$\pm$0.00} \\
       0.17$\pm$0.02 &         0.07$\pm$0.00 \\
 {\bf 0.09$\pm$0.01} &  {\bf  0.06$\pm$0.00} \\
\bottomrule
\end{tabular}
\
\begin{tabular}{lcc}
      {\bf MSD } &                    \\
\midrule
                 RMSE &                  AUC \\
\midrule
 {   8.93$\pm$NA} &  {  8.25$\pm$NA} \\
 { \bf  8.88$\pm$NA} &  {  6.79$\pm$NA} \\
 {   8.90$\pm$NA} &  {  \bf 5.19$\pm$NA} \\
 {   8.91$\pm$NA} &  {  5.22$\pm$NA} \\
 {   9.91$\pm$NA} &  {  6.06$\pm$NA} \\
 {   9.23$\pm$NA} &  {  5.24$\pm$NA} \\
\bottomrule
\end{tabular}

}
\caption{RMSE and AUC for the different data sets.}\label{tableSqErrorBostonConcrete}
\end{table}

\begin{figure}[!t]
	\begin{minipage}{0.3\textwidth}
	   \includegraphics[width=\textwidth]{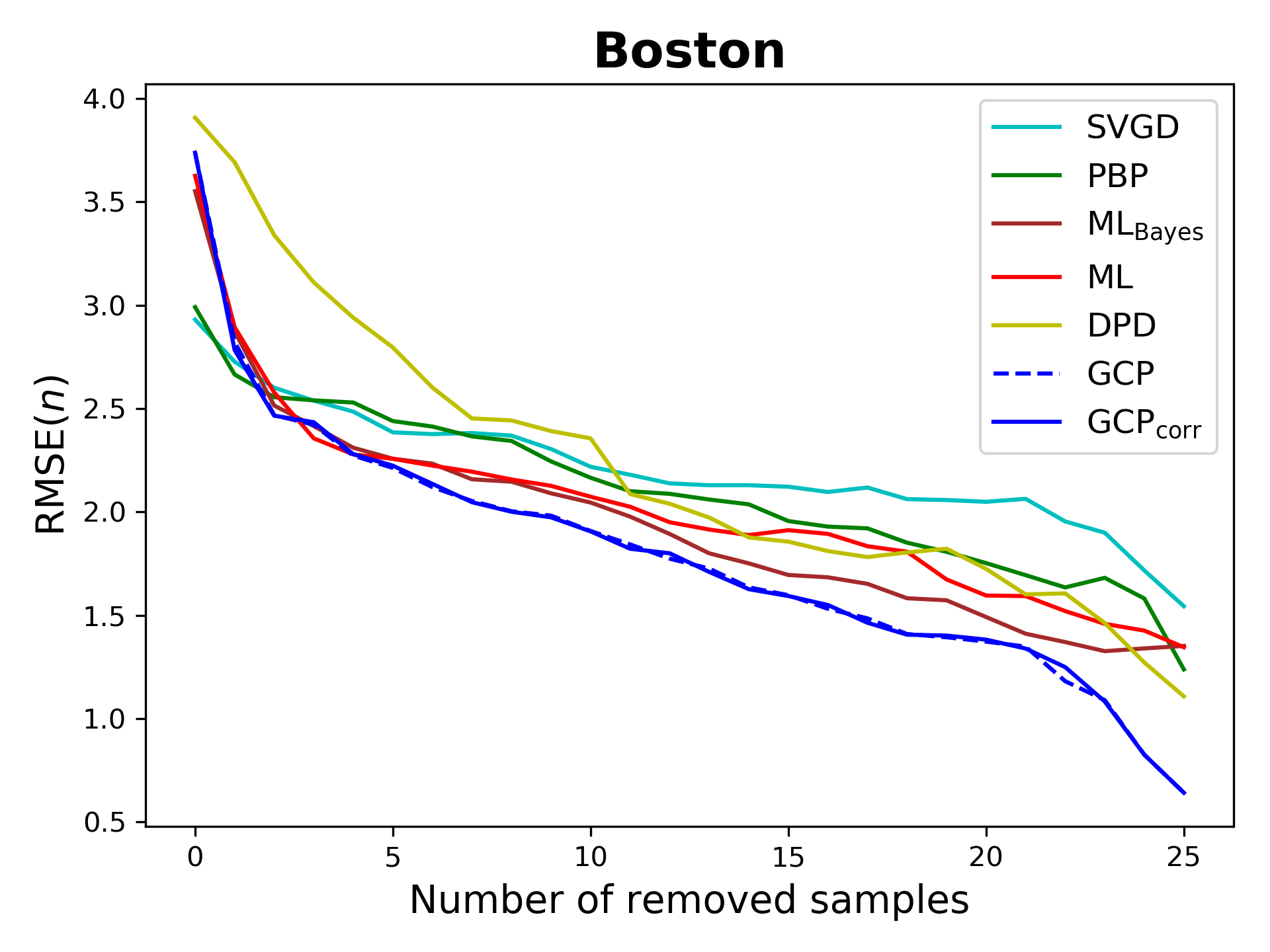}
	\end{minipage}
\hfill
	\begin{minipage}{0.3\textwidth}
       \includegraphics[width=\textwidth]{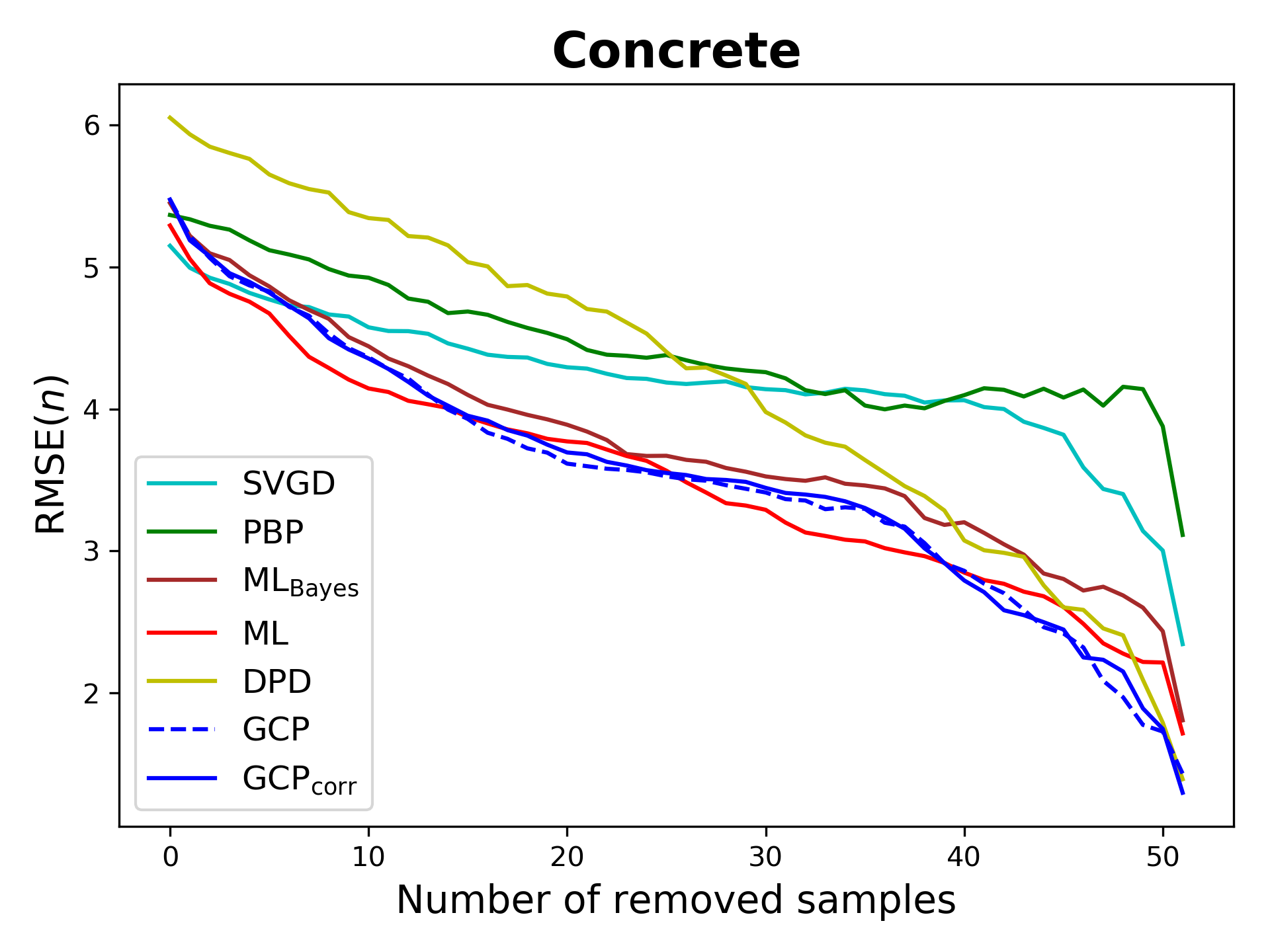}
	\end{minipage}
\hfill
	\begin{minipage}{0.3\textwidth}
       \includegraphics[width=\textwidth]{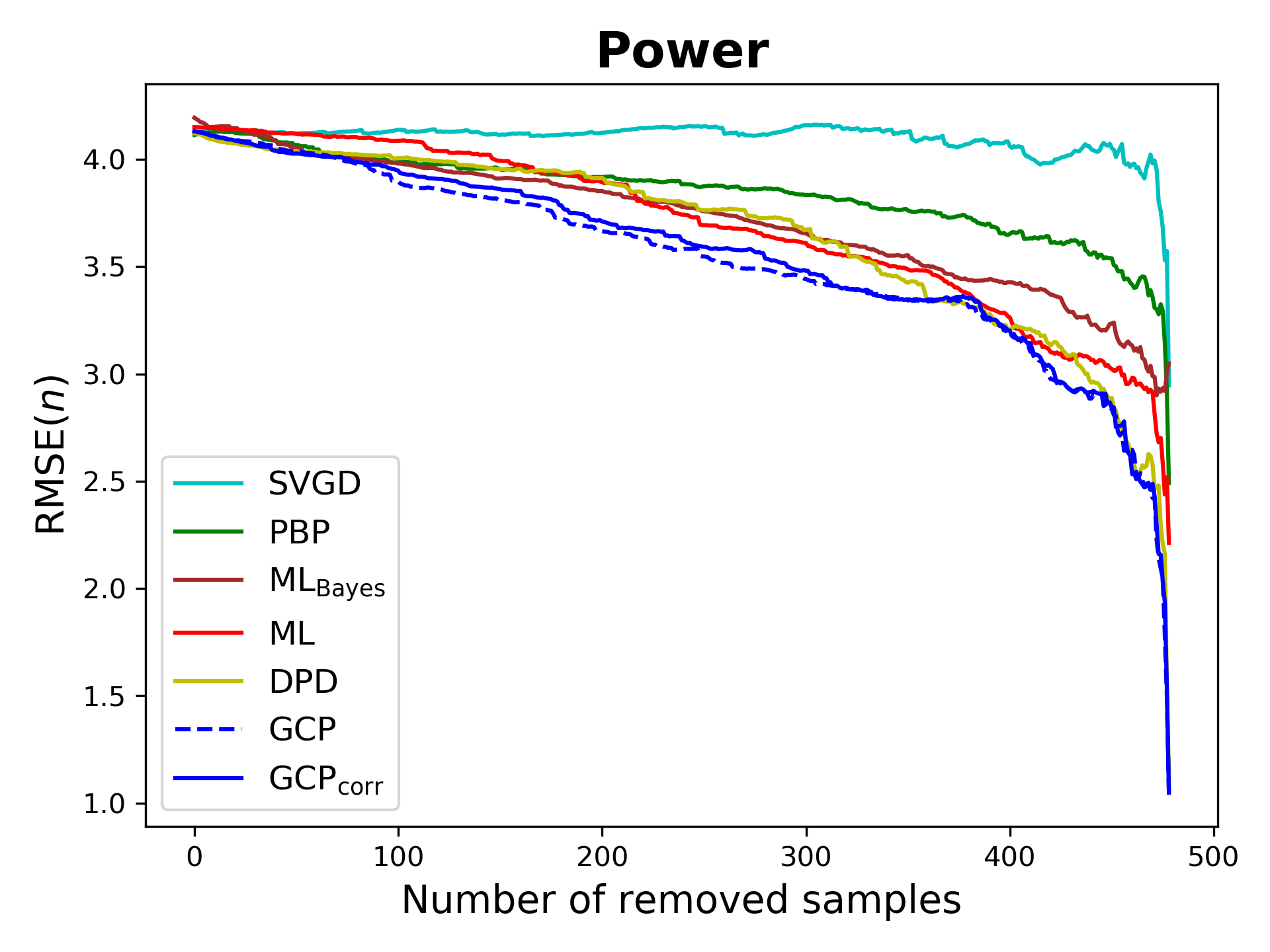}
	\end{minipage}

	\begin{minipage}{0.3\textwidth}
       \includegraphics[width=\textwidth]{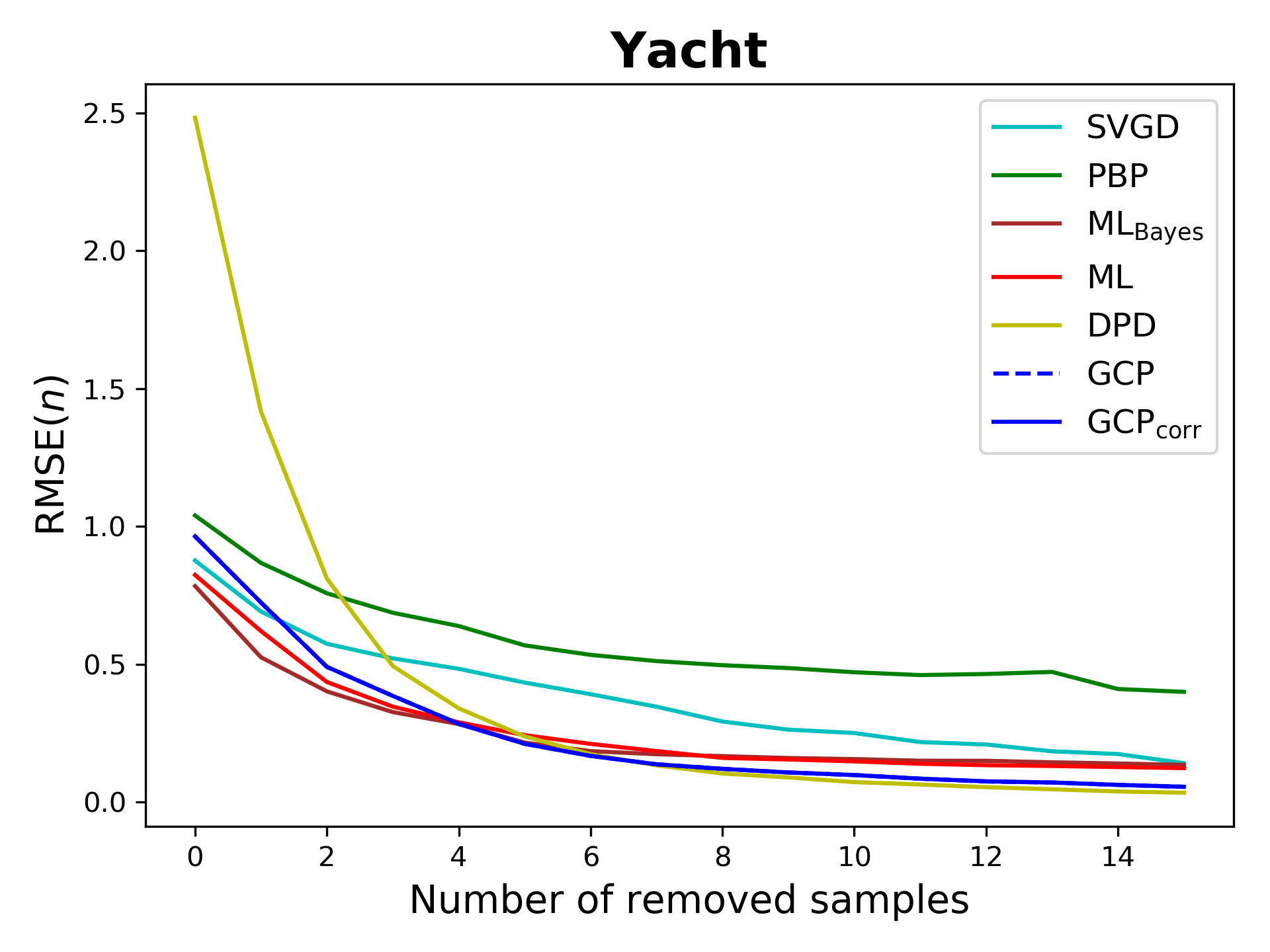}
	\end{minipage}
\hfill
	\begin{minipage}{0.3\textwidth}
       \includegraphics[width=\textwidth]{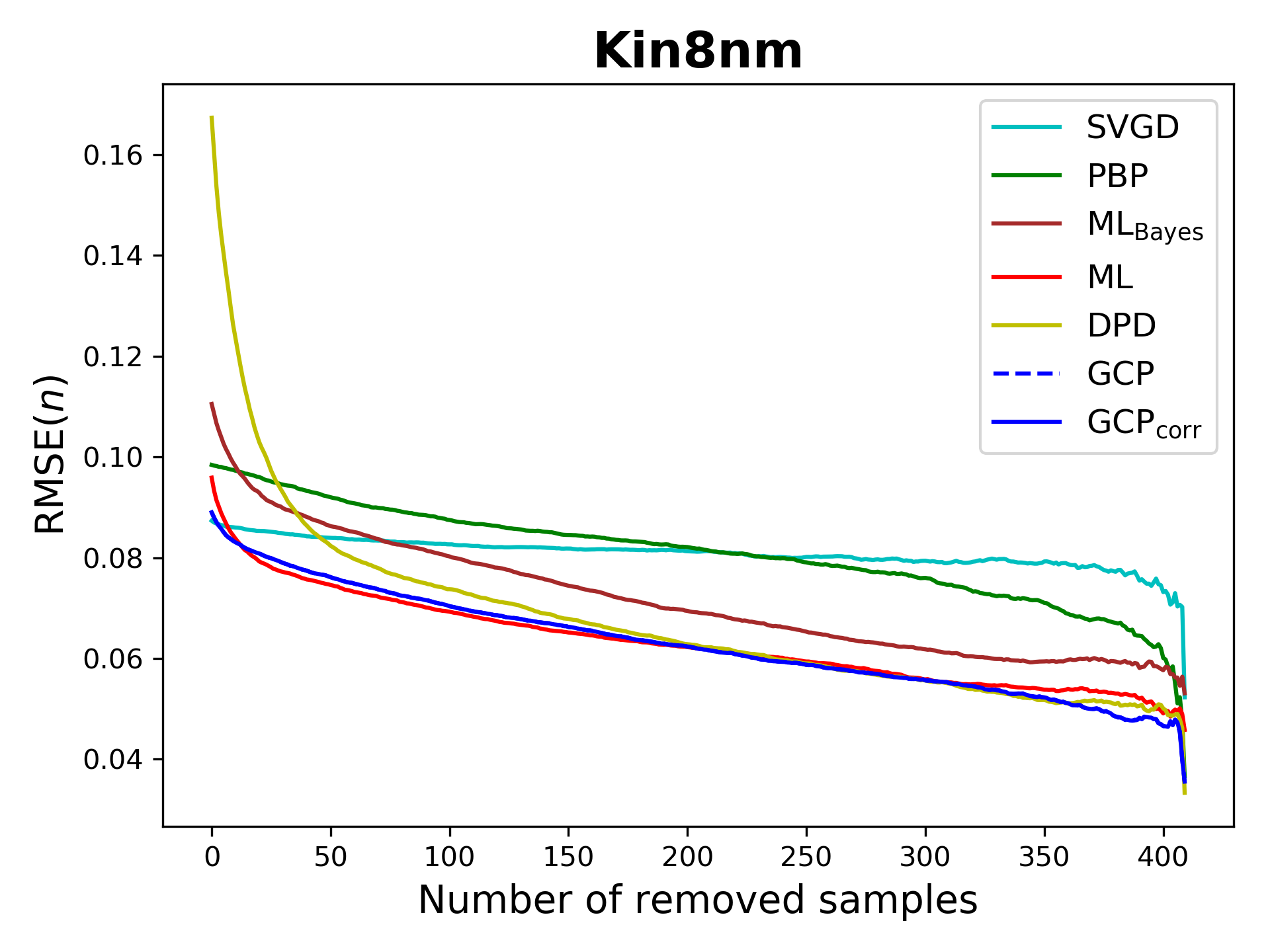}
	\end{minipage}
\hfill
	\begin{minipage}{0.3\textwidth}
       \includegraphics[width=\textwidth]{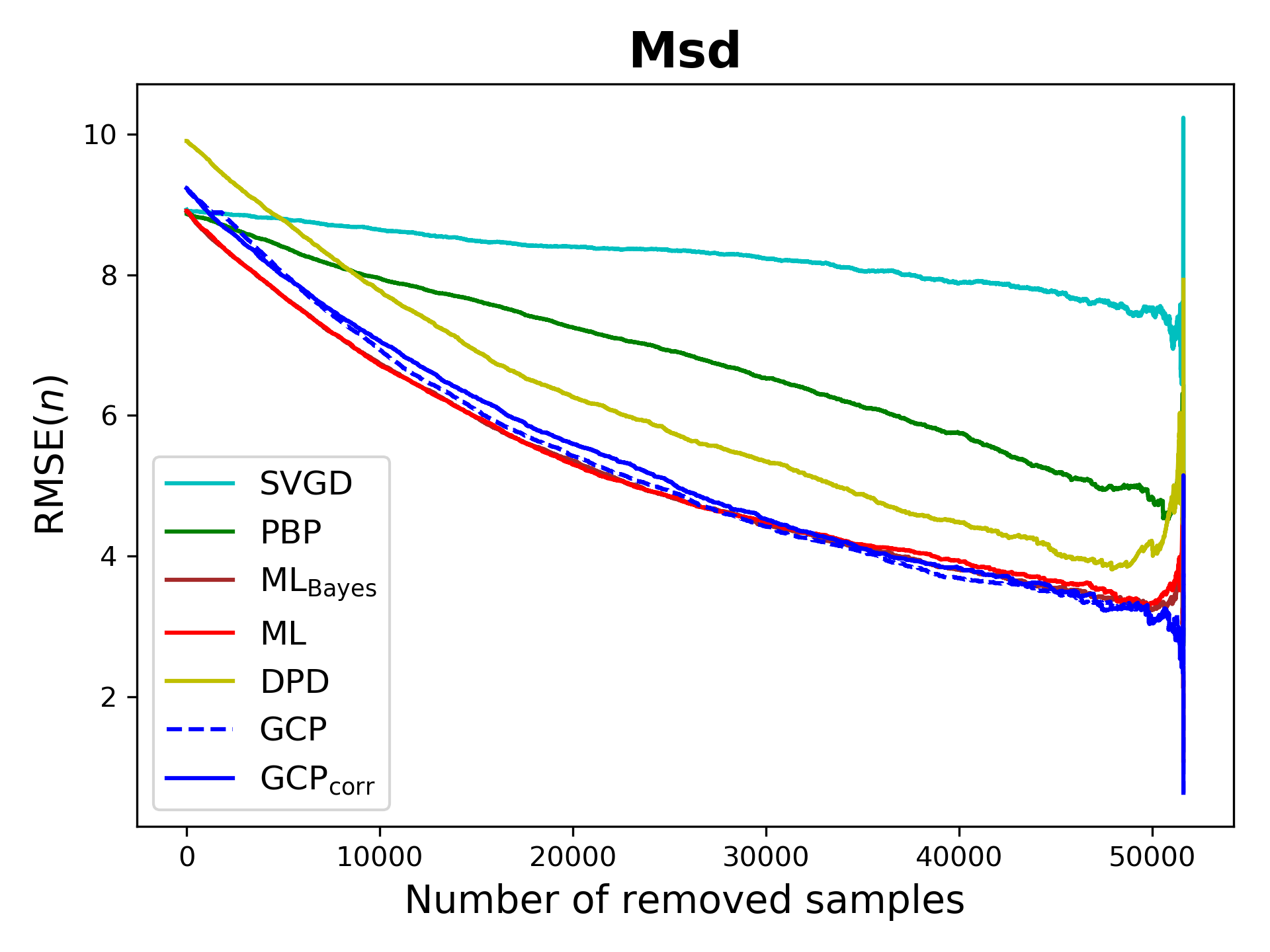}
	\end{minipage}
\caption{The curves ${\rm RMSE}(n)$ for the different methods and data sets from Table~\ref{tableSqErrorBostonConcrete}.}\label{figSqErrorRemovedSampleCurves}
\end{figure}

Figure~\ref{figSqErrorRemovedSampleCurves} shows the curves ${\rm RMSE}(n)$ for the different methods and data sets from Table~\ref{tableSqErrorBostonConcrete}. We see that the curve RMSE$(n)$ typically decays faster for the GCP compared with the other methods. In the absence of outliers, the curves of the GCP and GCP$_{\rm corr}$ practically coincide, yielding very close AUC values.

\subsection{Real world data sets: outliers in the training data sets}\label{subsecRealWorldOutliers}
We analyze the same methods and data sets as in Sec.~\ref{subsecRealWorld}, but now contaminated by outliers. For each training set, we randomly choose 5\% of samples and replace them by outliers. The outliers are sampled from the Gaussian distribution with the mean equal to the mean over all the targets in the original training set and standard deviation equal to ten times the standard deviation over the targets in the original training set. The results are presented in Table~\ref{tableSqErrorBostonConcreteOutliers}. Using the correction formula for $\tVest$ in~\eqref{eqEstMeanVarNetwork} allows the GCP$_{\rm corr}$ to obtain the best AUC values on all data sets. We also note that the Bayesian methods SVGD and PBP are especially sensitive to outliers, which is reflected in their high RMSE and AUC values.

\begin{table}
\resizebox{\textwidth}{!}{%
\begin{tabular}{lcc}
             &         {\bf Boston } &                    \\
\midrule
{} &                  RMSE &                  AUC \\
\midrule
SVGD               &        11.39$\pm$4.46 &        7.39$\pm$1.99 \\
PBP                &         8.64$\pm$2.70 &        5.70$\pm$1.40 \\
ML$_{\rm Bayes}$   &         4.17$\pm$1.28 &        3.56$\pm$1.14 \\
ML                 &  {\bf  3.66$\pm$1.01} &        3.40$\pm$1.50 \\
DPD                &  {\bf  3.74$\pm$1.70} &        2.21$\pm$0.57 \\
GCP                &  {\bf  3.66$\pm$1.54} &        2.83$\pm$1.36 \\
GCP$_{\rm corr}$   &  {\bf  3.66$\pm$1.54} &  {\bf 1.97$\pm$0.50} \\
\bottomrule
\end{tabular}
\
\begin{tabular}{lcc}
      {\bf Concrete } &                    \\
\midrule
                 RMSE &                  AUC \\
\midrule
       16.92$\pm$1.10 &       16.05$\pm$2.78 \\
        9.88$\pm$1.99 &        7.62$\pm$1.23 \\
        6.94$\pm$0.84 &        6.17$\pm$0.93 \\
        6.06$\pm$0.70 &        5.55$\pm$0.99 \\
 {\bf  5.17$\pm$0.81} &  {\bf 3.80$\pm$0.64} \\
        5.54$\pm$0.68 &        5.14$\pm$0.94 \\
        5.54$\pm$0.68 &  {\bf 3.65$\pm$0.69} \\
\bottomrule
\end{tabular}
\
\begin{tabular}{lcc}
      {\bf Power } &                    \\
\midrule
                 RMSE &                  AUC \\
\midrule
        5.53$\pm$0.42 &        5.18$\pm$0.45 \\
        4.60$\pm$0.32 &        4.27$\pm$0.29 \\
        4.67$\pm$0.32 &        4.56$\pm$0.43 \\
        4.73$\pm$0.28 &        4.69$\pm$0.44 \\
 {\bf  4.13$\pm$0.33} &        3.82$\pm$0.50 \\
        4.16$\pm$0.31 &  {\bf 3.65$\pm$0.37} \\
        4.16$\pm$0.31 &  {\bf 3.64$\pm$0.38} \\
\bottomrule
\end{tabular}
\
\begin{tabular}{lcc}
      {\bf Yacht } &                    \\
\midrule
                 RMSE &                  AUC \\
\midrule
       13.37$\pm$5.48 &        6.64$\pm$2.25 \\
      17.55$\pm$29.10 &        8.40$\pm$4.66 \\
        2.62$\pm$1.63 &        1.35$\pm$0.53 \\
        2.52$\pm$1.61 &        1.50$\pm$0.98 \\
        1.37$\pm$0.68 &  {\bf 0.27$\pm$0.13} \\
 {\bf  1.09$\pm$0.56} &        0.47$\pm$0.33 \\
 {\bf  1.09$\pm$0.56} &  {\bf 0.29$\pm$0.14} \\
\bottomrule
\end{tabular}
\
\begin{tabular}{lcc}
      {\bf Kin8nm } &                    \\
\midrule
                 RMSE &                  AUC \\
\midrule
       0.22$\pm$0.03 &         0.17$\pm$0.02 \\
       0.15$\pm$0.01 &         0.13$\pm$0.01 \\
       0.13$\pm$0.01 &         0.12$\pm$0.01 \\
       0.14$\pm$0.01 &         0.13$\pm$0.01 \\
       0.15$\pm$0.02 &  {\bf  0.07$\pm$0.01} \\
 {\bf 0.10$\pm$0.01} &         0.08$\pm$0.01 \\
 {\bf 0.10$\pm$0.01} &  {\bf  0.07$\pm$0.00} \\
\bottomrule
\end{tabular}
\
\begin{tabular}{lcc}
      {\bf MSD } &                    \\
\midrule
                 RMSE &                  AUC \\
\midrule
 {   9.56$\pm$NA} &  {  9.00$\pm$NA} \\
 {  \bf 9.05$\pm$NA} &  {  8.46$\pm$NA} \\
 {   9.17$\pm$NA} &  {  8.32$\pm$NA} \\
 {   9.10$\pm$NA} &  {  8.37$\pm$NA} \\
 {   9.80$\pm$NA} &  {  6.10$\pm$NA} \\
 {   9.31$\pm$NA} &  {  8.16$\pm$NA} \\
 {   9.31$\pm$NA} &  { \bf 5.57$\pm$NA} \\
\bottomrule
\end{tabular}

}
\caption{RMSE and AUC for the different data sets with 5\% of outliers.}\label{tableSqErrorBostonConcreteOutliers}
\end{table}

\begin{figure}[!t]
	\begin{minipage}{0.3\textwidth}
	   \includegraphics[width=\textwidth]{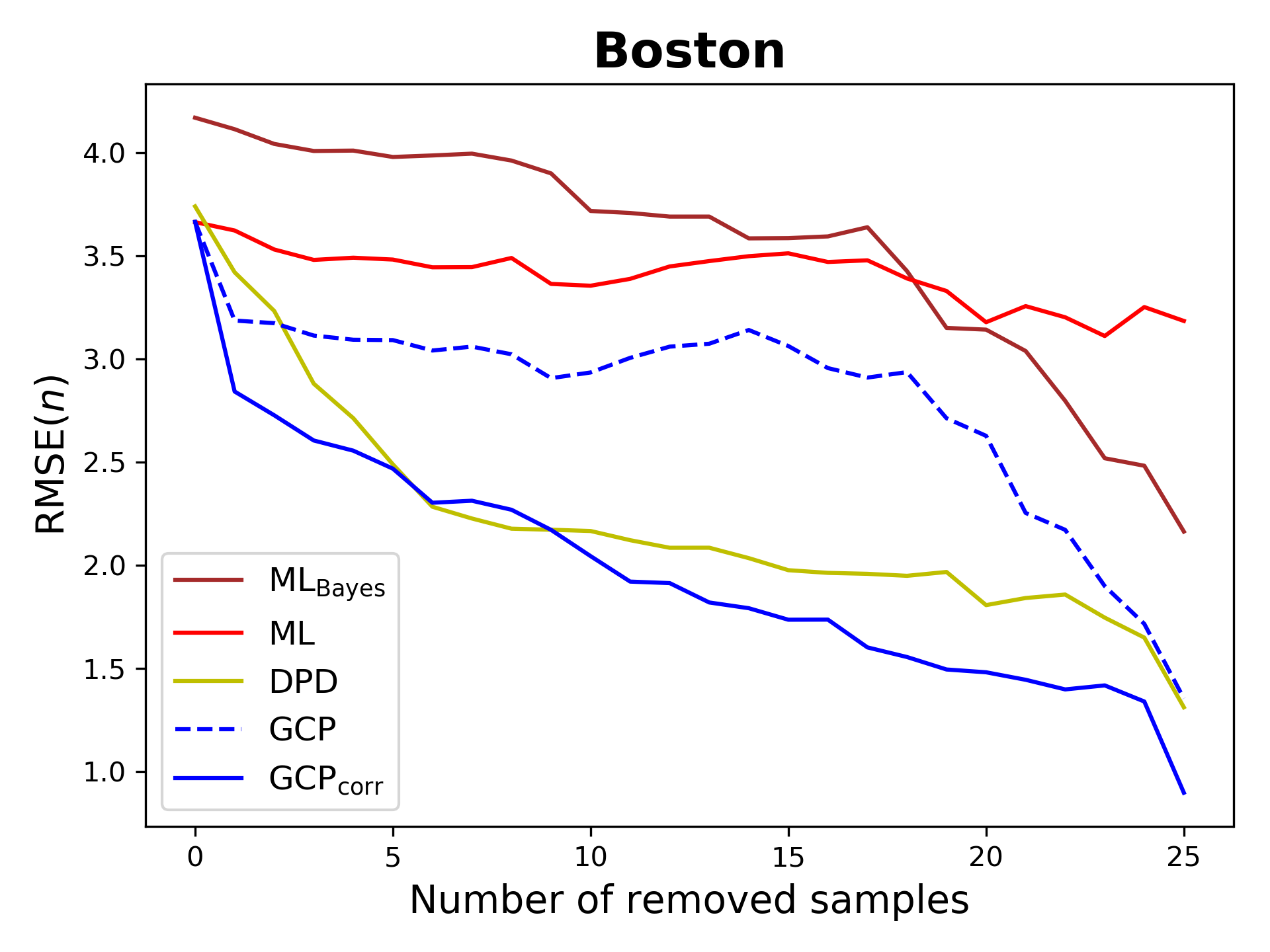}
	\end{minipage}
\hfill
	\begin{minipage}{0.3\textwidth}
       \includegraphics[width=\textwidth]{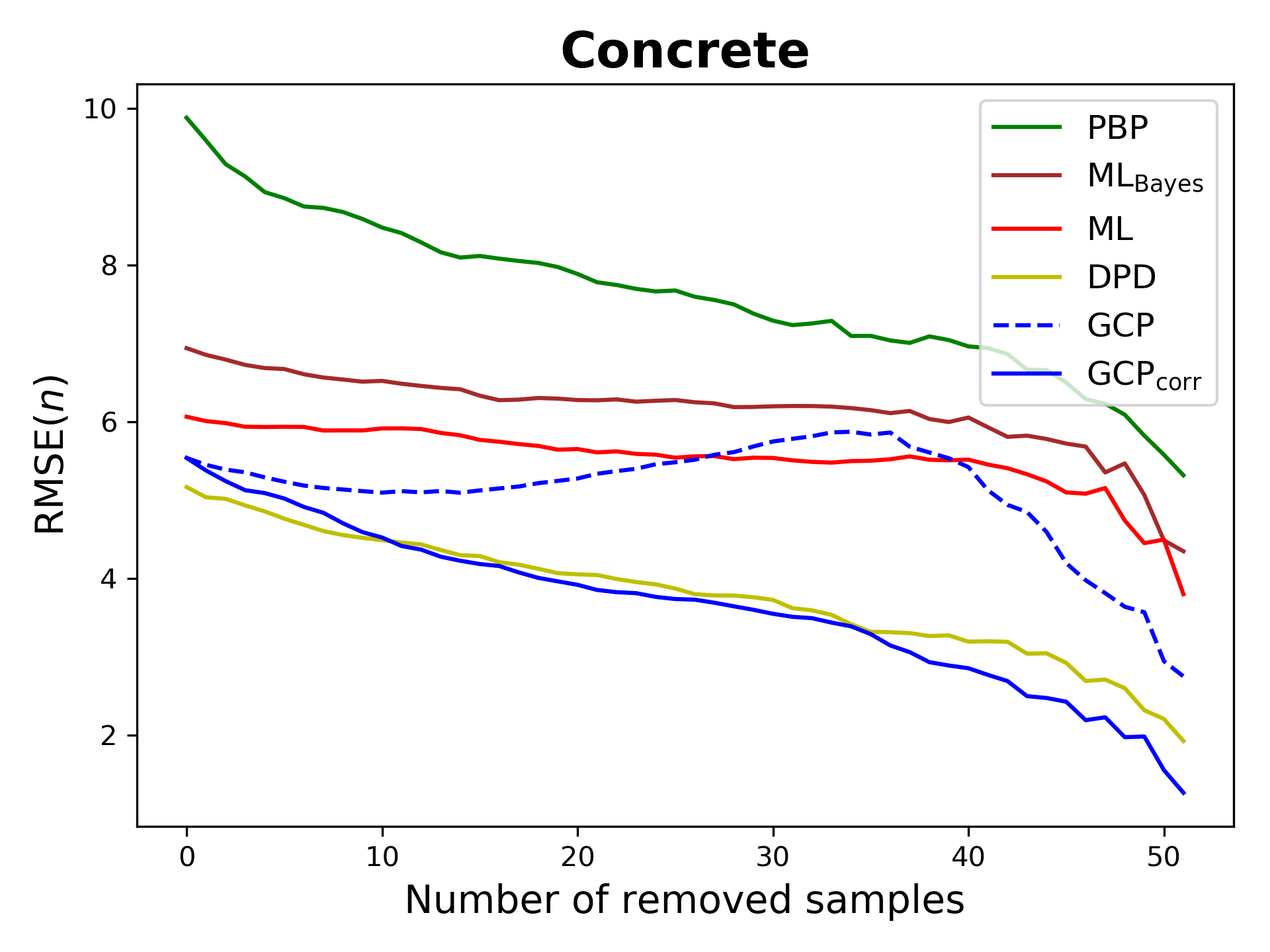}
	\end{minipage}
\hfill
	\begin{minipage}{0.3\textwidth}
       \includegraphics[width=\textwidth]{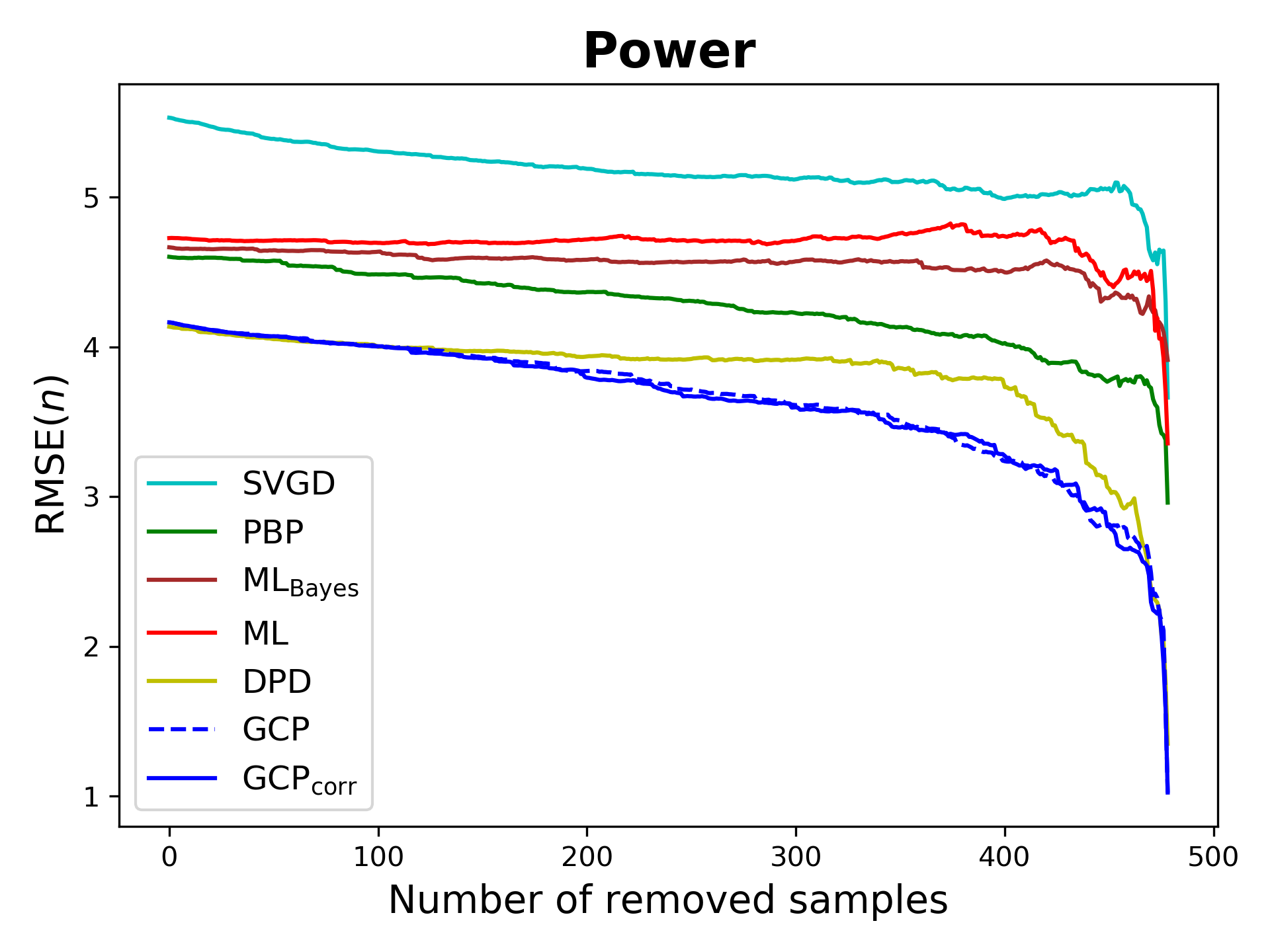}
	\end{minipage}

	\begin{minipage}{0.3\textwidth}
       \includegraphics[width=\textwidth]{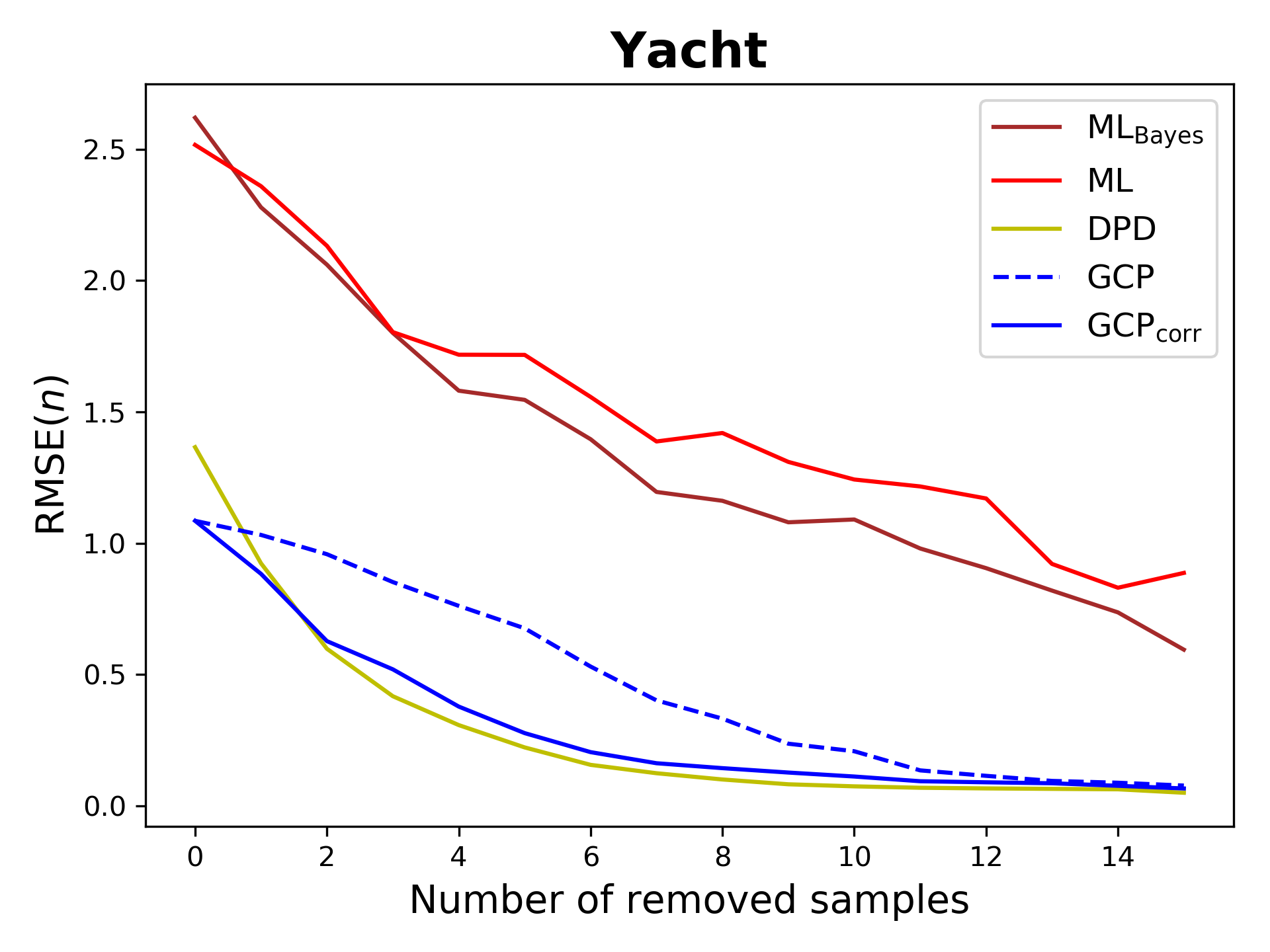}
	\end{minipage}
\hfill
	\begin{minipage}{0.3\textwidth}
       \includegraphics[width=\textwidth]{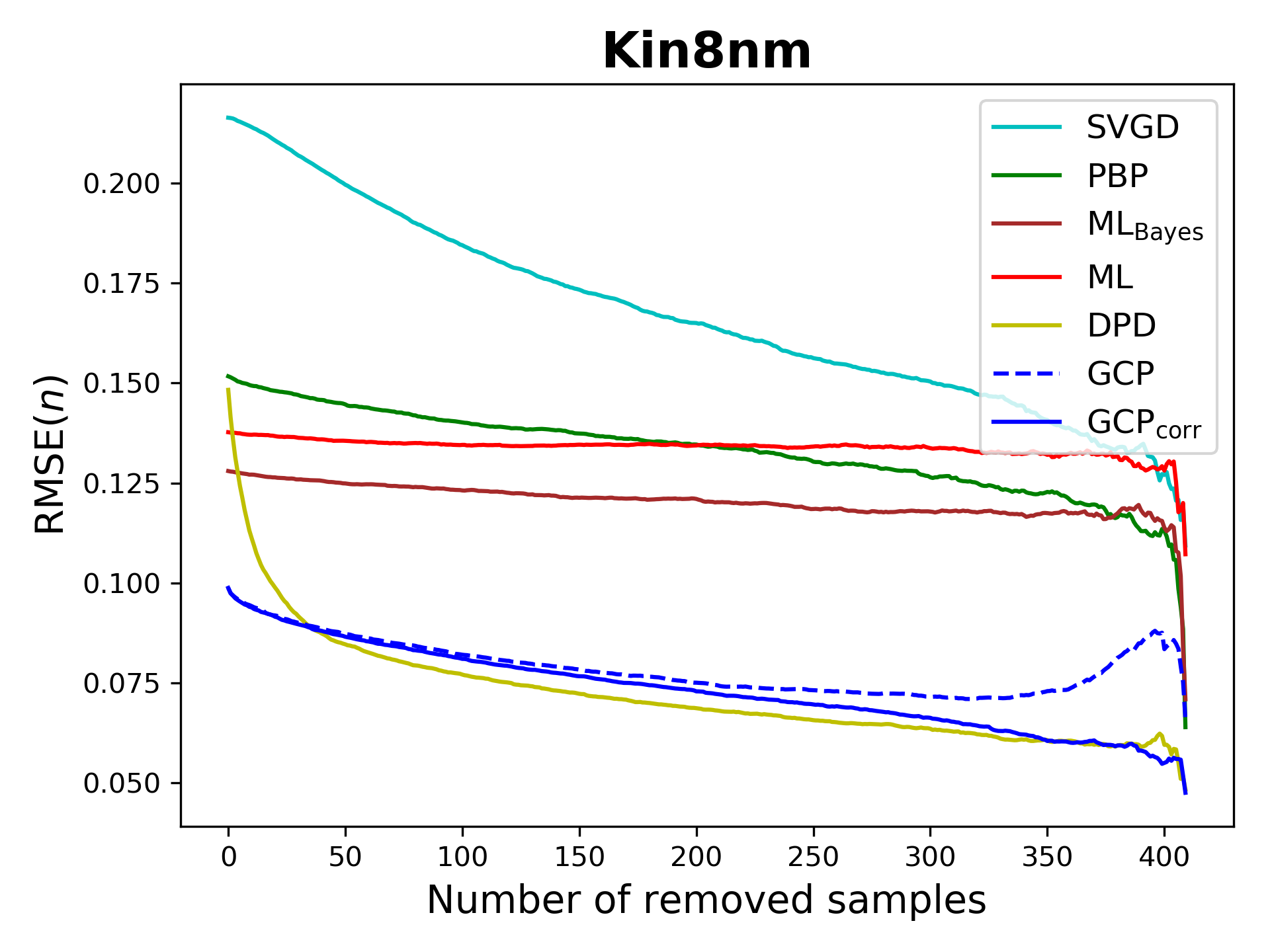}
	\end{minipage}
\hfill
	\begin{minipage}{0.3\textwidth}
       \includegraphics[width=\textwidth]{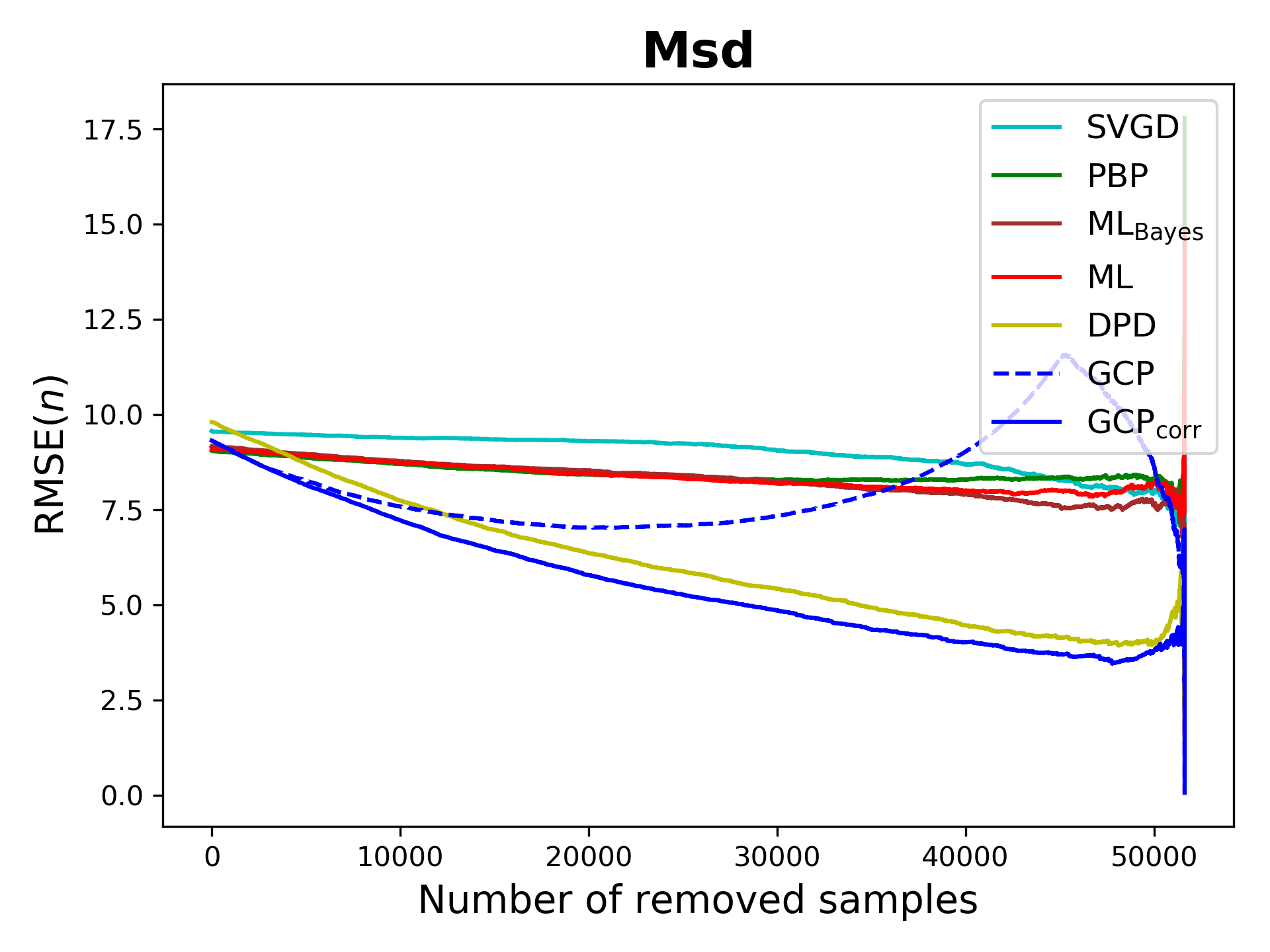}
	\end{minipage}
\caption{The curves ${\rm RMSE}(n)$ for the different methods and data sets with 5\% of outliers from Table~\ref{tableSqErrorBostonConcreteOutliers}. The curves for the PBP and SVGD are not shown on some of the plots because of their high RMSE and AUC values.}\label{figSqErrorRemovedSampleCurvesOutliers}
\end{figure}

Figure~\ref{figSqErrorRemovedSampleCurvesOutliers} shows the curves ${\rm RMSE}(n)$ for the different methods and data sets from Table~\ref{tableSqErrorBostonConcreteOutliers}. The curves of the GCP$_{\rm corr}$ are typically significantly below the corresponding curves of other methods, including the GCP without correction.

\section{Conclusion}\label{secConclusion}

Our goal was to approximate ground truth probability distributions by parametrically defined distributions. For their unknown parameters,  we introduced a prior distribution, whose parameters are learned by neural networks with deterministic weights. In such a setting, one cannot directly update the prior's parameters by the Bayesian rule, but one should rather update the network's weights. Hence, we proposed to replace a full Bayesian update of prior's parameters by a gradient descent step in the direction of minimizing the KL divergence from the posterior to the prior  distribution, which we called the GCP update. We showed that the GCP update is equivalent to the gradient ascent step that maximizes the likelihood of the predictive distribution. Interestingly, this result holds in general, independently of whether the posterior and prior distributions belong to the same family or not.

Next, we concentrated on the case where the ground truth distribution is normal with unknown mean and variance. A natural choice for the prior is the normal-gamma distribution. We obtained a dynamical system for its parameters that approximates the corresponding GCP update and analyzed it in detail. It revealed the convergence of the prior's parameters which is quite different from that for the standard Bayesian update, although in both cases the predictive Student's t-distribution converges to the ground truth normal distribution.

Furthermore, we analyzed how the GCP interacts with outliers in the training set. In the presence of outliers, the prior's parameter $\alpha$ (half the number of degrees of freedom of the predictive Student's t-distribution) does not tend to infinity any more. On one hand, this allows for a much better estimate of the mean of the ground truth normal distribution, compared with the ML method. On the other hand, this leads to overestimation of the variance of the ground truth distribution. We obtained, for the first time, an explicit formula that allows one to correct the estimate of the variance and recover the ground truth variance of the normal distribution.

Finally, we validated the GCP neural network on synthetic and six real-world data sets and compared it with the ML, DPD, SVGD, and PBP neural networks. We analyzed both clean data sets and data sets contaminated by 5\% of outliers. We measured the  trade-off between properly learning the mean and the variance (reflected in the AUC values) and the overall error (RMSE). The GCP was the only method that demonstrated the best AUC values {\em simultaneously} for clean and contaminated data sets.

To conclude, we indicate several directions of future research:
\begin{enumerate}

\item In case where the ground truth distribution of $\by$ is multivariate Gaussian with unknown mean and precision matrix, the conjugate prior is given by a family of normal-Wishart distributions. Due to Lemma~\ref{lNablawKy} and Remark~\ref{rGCPupdatePredictiveDistr} (item~\ref{rGCPupdatePredictiveDistr2}), the gradient conjugate prior update is equivalent to maximizing the likelihood of the predictive distribution, namely, multivariate Student distribution. Hence, one can obtain an analog of system~\eqref{eqODE4}, which will be a gradient system, too, and whose dynamics will be robust against outliers. Its analysis should allow one to derive a correction formula generalizing~\eqref{eqStandardCPEstimate} and  to reconstruct the ground truth  multivariate Gaussian distribution.

  \item
 A rigorous mathematical analysis of the influence of outliers on the dynamics of the prior's parameters seems to be feasible. One can relate the percentage of the outliers and a type of distribution they come from with the dynamical system~\eqref{eqODE4}, in which the expectations will be taken with respect to the new distribution (mixture of normal and the one from which the outliers are sampled). Further comparison with Bayesian methods based on robust divergences~\cite{Futami2017} and generalizing the DPD method is also needed.

  \item Section~\ref{secDynamicsFixedAlpha} shows that one can fix $\alpha$ and still recover the ground truth normal distribution, while Sec.~\ref{secRoleFixedAlpha} indicates how different values of $\alpha$ may influence the learning speed in clean and noisy regions. The influence of $\alpha$ on the fit of the GCP neural networks for real-world data sets would be an interesting practical question. Our preliminary analysis showed that fixing large $\alpha$ was beneficial for the largest MSD data set. For example, fixing $\alpha=30$ yielded ${\rm RMSE}=8.89$ and ${\rm AUC}=5.13$ (cf. Table~\ref{tableSqErrorBostonConcrete}).

  \item It is worth checking the GCP networks for other choices of  ground truth and prior distributions.

  \item The use of ensembles of MLs (called deep ensembles) was recently proposed in~\cite{Lakshminarayanan16,Lakshminarayanan17}. It is worth studying ensembles of GCPs capturing both aleatoric and epistemic uncertainty and being robust against outliers. Another possibility to capture epistemic uncertainty is to treat weights of the networks as random variables and learn their posterior distribution as is usually done for Bayesian neural networks.
\end{enumerate}

\appendix

\section{Properties of the function $A(\alpha)$: proof of Lemma~\ref{lAProperties}}\label{appendix}

For $\alpha>0$ and $A\in[0,\alpha]$, we study  equation~\eqref{eqDKbetanu02}, which is equivalent to the following:
\begin{equation}\label{eqE0}
E(\alpha,A):=\frac{1}{(2\pi)^{1/2}}\int_{-\infty}^{\infty} \frac{(2\alpha+1)z^2}{2(\alpha-A)+z^2}e^{-\frac{z^2}{2}}dz - 1 = 0.
\end{equation}

\begin{lemma}\label{lEalpha1}
  For each $\alpha>0$, equation~\eqref{eqE0} has a unique root $A=A(\alpha)\in(0,\alpha)$.   Furthermore,
  $$
  A(\alpha)=\alpha- \frac{4}{\pi} \alpha^2 + o(\alpha^2)\ \text{as }\alpha\to 0.
  $$
\end{lemma}
\proof
{\bf 1.} Note that $E(\alpha,A)$ is increasing with respect to $A\in(0,\alpha)$ and $E(\alpha,\alpha)=2\alpha>0$. Hence, it remains to show that $E(\alpha,0)<0$. We have
\begin{equation}\label{eqEalpha0}
E(\alpha,0)=\frac{1}{(2\pi)^{1/2}}\int_{-\infty}^{\infty} \frac{(2\alpha+1)z^2}{2\alpha+z^2}e^{-\frac{z^2}{2}}dz - 1 =  \frac{1}{(2\pi)^{1/2}}\int_{-\infty}^{\infty} \frac{z^2}{2\alpha+z^2}(1-z^2)e^{-\frac{z^2}{2}}dz.
\end{equation}
Now the inequality $E(\alpha,0)<0$ follows from~\eqref{eqEalpha0} and the monotonicity of $\frac{z^2}{2\alpha+z^2}$. Indeed,
$$
\begin{aligned}
\int_{0}^{1} \frac{z^2}{2\alpha+z^2}(1-z^2)e^{-\frac{z^2}{2}}dz & < \int_{0}^{1} \frac{1}{2\alpha+1}(1-z^2)e^{-\frac{z^2}{2}}dz \\
& =
\int_{1}^{\infty} \frac{1}{2\alpha+1}(z^2-1)e^{-\frac{z^2}{2}}dz < \int_{1}^{\infty} \frac{z^2}{2\alpha+z^2}(z^2-1)e^{-\frac{z^2}{2}}dz,
\end{aligned}
$$
where we have used the equality
$$
\int_{0}^{\infty} e^{-\frac{z^2}{2}}dz = \int_{0}^{\infty} z^2 e^{-\frac{z^2}{2}}dz = \left(\frac{\pi}{2}\right)^{1/2}.
$$

{\bf 2.} Now we prove the asymptotics of $A(\alpha)$. Using the function $F(x)$ defined in~\eqref{eqFAalpha}, we rewrite equation~\eqref{eqE0} in the form
\begin{equation}\label{eqDKbetanu02'}
(2\alpha+1)(\alpha-A)F(\alpha-A)-\alpha=0.
\end{equation}
Using the expansion of ${\rm erfc}(x)$ around $0$ (see~\cite[Sec.~7.1.6]{Abramowitz}) and formula~\eqref{eqFAalpha}, we have for all $x>0$
\begin{equation}\label{eqDKbetanu03'}
x F(x) = \frac{\pi^{1/2}+\delta}{2}x^{1/2},\quad \delta=o(1)\ \text{as } x\to 0.
\end{equation}
Now, for each $\alpha>0$, we represent $A = \alpha - k^2\alpha^2$ and prove that $k=\frac{2}{\pi^{1/2}}+o(1)$ as $\alpha\to 0$. Combining the representation of $A$ with~\eqref{eqDKbetanu02'} and~\eqref{eqDKbetanu03'}, we obtain
$$
\left(-1+\frac{k\pi^{1/2}}{2}+\frac{k\delta}{2}\right)+k(\delta+\pi^{1/2})\alpha = 0.
$$
Obviously, if $\alpha=\delta=0$, we have $k=\frac{2}{\pi^{1/2}}$. Hence, by the implicit function theorem, $k=\frac{2}{\pi^{1/2}}+o(1)$ as $\alpha,\delta\to 0$. Recalling that $\delta=o(1)$ as $\alpha\to 0$, we complete the proof.
\endproof

\begin{lemma}\label{lEalpha2}
  For each $\alpha\ge 1$, equation~\eqref{eqE0} has a unique root $A=A(\alpha)\in(0,1)$. Furthermore,
  $$
  A(\alpha) = 1 - \frac{3}{2\alpha} + o\left(\frac{1}{\alpha}\right) \ \text{as }\alpha\to \infty.
  $$
\end{lemma}
\proof
{\bf 1.} In the proof of Lemma~\ref{lEalpha1}, we have shown that $E(\alpha,0)<0$. Due to the monotonicity of $E(\alpha,A)$ with respect to $A\in(0,\alpha)$, it remains to show that $E(\alpha,1)>0$. Using that $2\alpha - 2\ge 0$, we have
$$
\begin{aligned}
E(\alpha,1) & =\frac{1}{(2\pi)^{1/2}}\int_{-\infty}^{\infty} \frac{(2\alpha+1)z^2}{2\alpha-2+z^2}e^{-\frac{z^2}{2}}dz - 1 =  \frac{1}{(2\pi)^{1/2}}\int_{-\infty}^{\infty} \frac{z^2(3-z^2)}{2\alpha-2+z^2} e^{-\frac{z^2}{2}}dz\\
&= \frac{1}{(2\pi)^{1/2}} \int_{-\infty}^\infty \int_0^\infty z^2 (3 - z^2) e^{-z^2 /2 - ( 2 \alpha - 2 + z^2) \xi} d\xi\, dz  \\
&= \frac{1}{(2\pi)^{1/2}} \int_0^\infty e^{- ( 2 \alpha - 2) \xi} d\xi \int_{-\infty}^\infty z^2 (3 - z^2) e^{-(\xi + 1/2) z^2} dz  \\
&= 3   \int_0^\infty \exp(- ( 2 \alpha - 2) \xi) \left( \frac{1}{(2 \xi + 1)^{3/2}} - \frac{1}{(2 \xi + 1)^{5/2}} \right)d\xi >0.
\end{aligned}
$$

{\bf 2.} Now we prove the asymptotics of $A(\alpha)$.  Using the expansion of ${\rm erfc}(x)$ around $\infty$ (see~\cite[Sec.~7.1.23]{Abramowitz}) and formula~\eqref{eqFAalpha}, we have for all $x>0$
\begin{equation}\label{eqDKbetanu04'}
x F(x) = \frac{1}{2}-\frac{1}{4x} + \frac{3}{8 x^2} - \frac{15+\delta}{16 x^3},\quad \delta=o(1)\ \text{as } x\to \infty.
\end{equation}
Now, for each $\alpha>0$, we representing $A = 1 - {k}/{\alpha}$ and prove that $k=3/2 + o(1)$ as $\alpha\to \infty$. Combining the representation of $A$ with~\eqref{eqDKbetanu02'} and~\eqref{eqDKbetanu04'}, we obtain
$$
\frac{-12+38 k z^2-28 k^2 z^3+8 k^3 z^4+8 k-28 k z-33 z-2 \delta+16 k^2 z^2-z \delta}{(1-z+k z^2)^3} = 0,\quad z:=\frac{1}{\alpha}.
$$
Obviously, if $\alpha=\delta=0$, we have $k=3/2$. Hence, by the implicit function theorem, $k=3/2+o(1)$ as $z,\delta\to 0$. Recalling that $\delta=o(1)$ as $z\to 0$, we complete the proof.
\endproof

\begin{lemma}\label{lAODE}
   The function $A(\alpha)$ satisfies the differential equation in~\eqref{eqAODE}.
\end{lemma}
\proof
Denoting by $E_\alpha$ and $E_A$ the partial derivatives of the function $E(\alpha,A)$ with respect to $\alpha$ and $A$, respectively, and using the implicit function theorem, we have
\begin{equation}\label{eqA'}
A' = -E_A^{-1} E_\alpha.
\end{equation}
Since $A(\alpha)$ satisfies the equation in~\eqref{eqE0}, we obtain
\begin{equation}\label{eqEalpha}
\begin{aligned}
E_\alpha &= -2 \frac{1}{(2\pi)^{1/2}} \int_{-\infty}^\infty \frac{(2 \alpha + 1) z^2}{(2 (\alpha - A) + z^2)^2} e^{-z^2/2} dz + \frac{1}{(2\pi)^{1/2}} \int_{-\infty}^\infty \frac{2 z^2}{(2 (\alpha - A) + z^2)} e^{-z^2/2} dz\\
& = - E_A + \frac{2}{2\alpha+1}
\end{aligned}
\end{equation}
Now we calculate $E_A$ and integrate by parts:
$$
\begin{aligned}
E_A = \frac{1}{(2\pi)^{1/2}} \int_{-\infty}^\infty \frac{(2 \alpha + 1) (2 z)}{(2 (\alpha - A) + z^2)^2} z e^{-z^2/2} dz
&= \frac{1}{(2\pi)^{1/2}} \int_{-\infty}^\infty \frac{(2 \alpha + 1)(1-z^2)}{(2 (\alpha - A) + z^2)} e^{-z^2/2} dz.
\end{aligned}
$$
Again using the equality in~\eqref{eqE0}, we obtain
\begin{equation}\label{eqEA}
\begin{aligned}
E_A & =
  \frac{1}{(2\pi)^{1/2}} \int_{-\infty}^\infty \frac{2 \alpha + 1}{(2 (\alpha - A) + z^2)} e^{-z^2/2} dz - 1\\
&=   \frac{1}{2 (\alpha - A)} \frac{1}{(2\pi)^{1/2}} \int_{-\infty}^\infty \frac{(2 \alpha + 1)(2(\alpha - A) + z^2 - z^2)}{(2 (\alpha - A) + z^2)} e^{-z^2/2} dz -1 \\
&=  \frac{2\alpha + 1 - 1}{2 (\alpha - A)} -1 = \frac{\alpha}{\alpha - A} - 1.
\end{aligned}
\end{equation}
Combining~\eqref{eqA'}--\eqref{eqEA}, we obtain~\eqref{eqAODE}.
\endproof

\begin{lemma}\label{lAmonotone}
  For all $\alpha>0$, we have $A'(\alpha)>0$  and
\begin{equation}\label{eqAMonotone}
\frac{2\alpha}{2\alpha+3} < A(\alpha).
\end{equation}
\end{lemma}
\proof
It suffices to show that the right-hand side of~\eqref{eqAODE} is positive for all $\alpha>0$, which is equivalent to~\eqref{eqAMonotone}.
We consider the function
$$
g(\alpha):=A(\alpha)-\frac{2\alpha}{2\alpha+3}
$$
and show that $g(\alpha)>0$ for all $\alpha>0$. Assume this is not true. Since $g(\alpha)>0$ for all sufficiently small $\alpha>0$ (due to the asymptotics in Lemma~\ref{lEalpha1}) and $\lim\limits_{\alpha\to\infty}g(\alpha)=0$ (due to the asymptotics in Lemma~\ref{lEalpha2}), this would imply that
\begin{equation}\label{eqg'0}
 g'(\alpha)=0,\quad g(\alpha)\le 0\quad\text{for some } \alpha>0.
\end{equation}
Using the fact that $\left(\frac{2\alpha}{2\alpha+3}\right)'>0$ and applying Lemma~\ref{lAODE}, we have
$$
g'(\alpha)< A'(\alpha) = \frac{2\alpha(A-1)+3A}{(2\alpha+1)A}.
$$
Since $A(\alpha)\le \frac{2\alpha}{2\alpha+3}$ for $\alpha$ in~\eqref{eqg'0}, we obtain
$$
g'(\alpha) < \frac{1}{(2\alpha+1)A}\left(2\alpha\left(\frac{2\alpha}{2\alpha+3}-1\right) + \frac{6\alpha}{2\alpha+3} \right) = 0,
$$
which contradicts~\eqref{eqg'0}.
\endproof 

\section{Proof of Lemma~\ref{lDKbetanu}}\label{appendixLemmaDKbetanu}

{\bf 1.} Let us show that $\dot\nu=0$ if $\sigma=\sigma_0(\alpha)$. Due to~\eqref{eqDKnuExp},
$$
\bbE\left[\frac{\partial K}{\partial \nu}\right]\bigg|_{\sigma=\sigma_0(\alpha)} =  \frac{1}{2\nu(\nu+1)}\left(\frac{1}{(2\pi)^{1/2}}\int_{-\infty}^{\infty} \frac{(2\alpha+1)z^2}{2(\alpha-A)+z^2}e^{-\frac{z^2}{2}}dz - 1\right).
$$
Now the equation $\bbE\left[\frac{\partial K}{\partial \nu}\right]\big|_{\sigma=\sigma_0(\alpha)}=0$ can be rewritten as follows:
\begin{equation}\label{eqDKbetanu01}
\int_{-\infty}^{\infty} \frac{(2\alpha+1)z^2}{2(\alpha-A)+z^2}e^{-\frac{z^2}{2}}dz - 1=0,
\end{equation}
which is equivalent to~\eqref{eqDKbetanu02}.
By Lemma~\ref{lAProperties}, it has a unique root $A=A(\alpha)\in(0,\min(\alpha,1))$.

Since $\bbE\left[\frac{\partial K}{\partial \nu}\right]$ is decreasing with respect to $\sigma$ due to~\eqref{eqDKnuExp}, the assertions about~$\nu$ in~\eqref{eqDKbetanu1} follow.

{\bf 2.}  Next, we show that $\dot\beta=0$ if $\sigma=\sigma_0(\alpha)$. Due to~\eqref{eqDKbetaExp}, \eqref{eqDKbetanu02}, and~\eqref{eqFAalpha},
$$
\begin{aligned}
\bbE\left[\frac{\partial K}{\partial \beta}\right]\bigg|_{\sigma=\sigma_0(\alpha)} &= \frac{\nu+1}{\nu V}
\left(\frac{1}{(2\pi)^{1/2}} \int_{-\infty}^{\infty} \frac{2\alpha+1}{2(\alpha-A)+z^2}e^{-\frac{z^2}{2}}dz - \frac{\alpha}{\alpha-A}\right)\\
&= \frac{\nu+1}{\nu V}\left( (2\alpha+1) F(\alpha-A) - \frac{\alpha}{\alpha-A}\right) = 0.
\end{aligned}
$$
Since the expression in the brackets in~\eqref{eqDKbetaExp} is   increasing with respect to $\sigma$,  the assertions about $\beta$ in~\eqref{eqDKbetanu1} follow.

\section{Proof of Lemma~\ref{lCurvesNuBeta}}\label{appendixLemmaCurvesNuBeta}

Using the equality $\frac{1}{(2\pi)^{1/2}}\int_{-\infty}^{\infty} e^{-z^2/2}dz=1$, we rewrite~\eqref{eqDKnuExp} as follows:
$$
\bbE\left[\frac{\partial K}{\partial \nu}\right] =
-\frac{1}{\nu(\nu+1)}\left(\frac{\alpha+1/2}{(2\pi)^{1/2}}\int_{-\infty}^{\infty} \frac{1}{1 + \frac{V}{\sigma} \frac{z^2}{2}} e^{-\frac{z^2}{2}}dz  - \alpha\right).
$$
Combining this relation with~\eqref{eqDKbetaExp} shows that the points $(\nu(t),\beta(t))$ belong to the integral curves of the differential equation~\eqref{eqCurvesNuBetaODE}, because
\[
\frac{d \beta}{d \nu} = \frac{\bbE\left[\frac{\partial K}{\partial \beta}\right]}{\bbE\left[\frac{\partial K}{\partial \nu}\right]} = -\frac{\nu ( \nu+1)}{\beta}
\]
Separating variables in this equation, one can see that the integral curves are given by~\eqref{eqCurvesNuBeta}.

\section{Proof of Lemma~\ref{lDKalpha}}\label{appendixLemmaDKalpha}

{\bf 1.} Due to~\eqref{eqDKalphaExp},
\begin{equation}\label{eqDKalpha1}
\begin{aligned}
 \bbE\left[\frac{\partial K}{\partial \alpha}\right]\bigg|_{\sigma=\sigma_\varkappa(\alpha)} &= \frac{1}{(2\pi)^{1/2}}\int_{-\infty}^{\infty} \ln\left(1+\frac{1}{\left(1-\frac{\varkappa}{\alpha^2}\right)(\alpha-A)}\frac{z^2}{2}\right)e^{-\frac{z^2}{2}}dz \\
 & + \Psi(\alpha)-\Psi\left(\alpha+\frac{1}{2}\right).
\end{aligned}
\end{equation}
Using that $\Psi(\alpha)-\Psi(\alpha+1/2)\to 0$ and $A(\alpha)\to 1$ as $\alpha\to\infty$, we see that
\begin{equation}\label{eqDKalphaInfty}
 \lim\limits_{\alpha\to\infty} \bbE\left[\frac{\partial K}{\partial \alpha}\right]\bigg|_{\sigma=\sigma_\varkappa(\alpha)} = 0.
\end{equation}

{\bf 2.}  To complete the proof of~\eqref{eqDKalpha}, it suffices (due to~\eqref{eqDKalphaInfty}) to show that the derivative of the right-hand side in~\eqref{eqDKalpha1} is positive for $\varkappa=0$ and all $\alpha>0$.
We denote the derivative of the right-hand side in~\eqref{eqDKalpha1} by $G_\varkappa(\alpha)$. To calculate it, we set
$$
B(\alpha) := 1-\frac{\varkappa}{\alpha},\quad \gamma(\alpha)=\frac{1}{2B(\alpha-A)}.
$$
Then
\begin{equation}\label{eqDKalpha1'}
\frac{\partial}{\partial\alpha}\ln(1+\gamma z^2) = -\frac{z^2}{1+\gamma z^2}\left(\frac{(1-A')}{2B(\alpha-A)^2} + \frac{B'}{2B^2(\alpha-A)}\right).
\end{equation}

If $\varkappa=0$, then $B=1$, $B'=0$, and (due to~\eqref{eqE0})
\begin{equation}\label{eqDKalpha1''}
\frac{1}{(2\pi)^{1/2}}\int_{-\infty}^{\infty} \frac{z^2}{1+\gamma z^2} = \frac{1}{(2\pi)^{1/2}}\int_{-\infty}^{\infty} \frac{z^2}{1+\frac{1}{2(\alpha-A)} z^2} = \frac{2(\alpha-A)}{2\alpha+1}.
\end{equation}
Using~\eqref{eqDKalpha1'} and~\eqref{eqDKalpha1''}, we obtain
\begin{equation*}
G_0(\alpha) = -\frac{1-A'}{(2\alpha+1)(\alpha-A)} + \Psi'(\alpha)-\Psi'\left(\alpha+\frac{1}{2}\right),
\end{equation*}
or, using Lemma~\ref{lAProperties} (item~\ref{lAProperties4}),   equivalently,
\begin{equation*}
G_0(\alpha) = -\frac{2}{(2\alpha+1)^2 A} + \Psi'(\alpha)-\Psi'\left(\alpha+\frac{1}{2}\right).
\end{equation*}
Due to the inequality $A(\alpha)>\frac{2\alpha}{2\alpha+3}$ (see the first inequality in Lemma~\ref{lAProperties}, item~\ref{lAProperties3}),
$$
G_0(\alpha) > \tilde G_0(\alpha),
$$
where
\begin{equation}\label{eqDeriv2alphaTildeG0'}
\tilde G_0(\alpha) = -\frac{2\alpha+3}{(2\alpha+1)^2 \alpha} + \Psi'(\alpha)-\Psi'\left(\alpha+\frac{1}{2}\right).
\end{equation}
Therefore, for the proof of~\eqref{eqDKalpha} it remains to show that
\begin{equation}\label{eqG0alphaPositive}
\tilde G_0(\alpha)>0\quad\text{for all } \alpha>0.
\end{equation}


{\bf 2.1.} First, we prove~\eqref{eqG0alphaPositive} for large $\alpha$.
Using the asymptotics~\cite[Sec.~6.4.12]{Abramowitz}
\begin{equation}\label{eqAsympPsi'infty}
\Psi'(\alpha) - \Psi'\left(\alpha+\frac{1}{2}\right) = \frac{1}{2\alpha^2} + \frac{1}{4\alpha^3}+ O\left(\frac{1}{\alpha^5}\right)\quad\text{as }\alpha\to \infty
\end{equation}
 we obtain from~\eqref{eqDeriv2alphaTildeG0'}
\begin{equation}\label{eqG0alphaPositiveLargeAlpha}
\begin{aligned}
\tilde G_0(\alpha) & =  \frac{3}{8\alpha^4} + O\left(\frac{1}{\alpha^5}\right) > 0\quad\text{for all sufficiently large}\ \alpha.
\end{aligned}
\end{equation}

{\bf 2.2.} Due to~\eqref{eqG0alphaPositiveLargeAlpha}, to complete the proof of~\eqref{eqG0alphaPositive} it now suffices to show that
$$
\tilde G_0(\alpha) - \tilde G_0(\alpha+1) > 0 \quad\text{for all } \alpha>0.
$$
Applying the recurrence relation $\Psi'(z+1)=\Psi'(z) - 1/z^2$ (see~\cite[Sec.~6.4.6]{Abramowitz}), we obtain from~\eqref{eqDeriv2alphaTildeG0'}
$$
\begin{aligned}
\tilde G_0(\alpha)- \tilde G_0(\alpha+1) & = -\frac{2\alpha+3}{(2\alpha+1)^2 \alpha} + \frac{2(\alpha+1)+3}{(2(\alpha+1)+1)^2 (\alpha+1)}    +\frac{1}{\alpha^2} -  \frac{1}{(\alpha+1/2)^2}\\
& = \frac{3(4\alpha+3)}{(2\alpha+1)(\alpha+1)(2\alpha+3)^2 \alpha^2}>0,
\end{aligned}
$$
which proves~\eqref{eqG0alphaPositive} and thus completes the proof of~\eqref{eqDKalpha}.

{\bf 3.} Now consider the case $\varkappa>0$. Note that the expression in the brackets in~\eqref{eqDKalpha1'} is positive due to Lemma~\ref{lAProperties} and the fact that $B'(\alpha)>0$. Hence,  we obtain from~\eqref{eqDKalpha1'}
$$
\begin{aligned}
\frac{\partial}{\partial\alpha}\ln(1+\gamma z^2) & = -2B(\alpha-A)\frac{z^2}{2B(\alpha-A)+z^2}\left(\frac{1-A'}{2B(\alpha-A)^2} + \frac{B'}{2B^2(\alpha-A)}\right)\\
&< - \frac{z^2}{2(\alpha-A)+z^2}\left(\frac{1-A'}{\alpha-A} + \frac{B'}{B}\right) .
\end{aligned}
$$
Combining the latter inequality with~\eqref{eqDKalpha1''} and  using Lemma~\ref{lAProperties} (item~\ref{lAProperties4}), we have
$$
  G_\varkappa(\alpha) <  -\frac{2}{ (2\alpha+1)^2 A }-\frac{B'}{B (2\alpha+1)}+ \Psi'(\alpha)-\Psi'\left(\alpha+\frac{1}{2}\right).
$$
Additionally using the asymptotics in~\eqref{eqAsympPsi'infty} and the expansion of $A(\alpha)$ in~\eqref{eqPropertiesA} as $\alpha\to\infty$, we obtain
$$
G_\varkappa(\alpha)  < -\frac{\varkappa}{2\alpha^3} + o\left(\frac{1}{\alpha^3}\right)<0\quad\text{for all sufficiently large } \alpha.
$$

\section{Hyperparameters}\label{appendixHyperparameters}

When we fit different methods on the real world data sets, we normalize them so that the input features and the targets
have zero mean and unit variance in the training set. We used minibatch 5 on Boston, Concrete, and Yacht, minibatch 10 on Power and Kin8nm, and minibatch 5000 on MSD. We used Adam (with $\beta_1=0.9$, $\beta_2=0.999$), RmsProp (with $\rho=0.5$), or Nesterov momentum (with momentum $0.9$) optimizers for fitting the ML, DPD, and GCP methods. In case of each optimizer, we performed a grid search for the learning rate in the range $\{0.00002, 0.00005, 0.0001, 0.0002, 0.0007, 0.001, 0.005\}$ and for the dropout rate in the range $\{0,0.1,0.2,0.3,0.4\}$. The optimizers and the parameters yielding the best AUC are presented in Tables~\ref{tableHyperparametersBostonConcrete} and~\ref{tableHyperparametersYacht}. We trained the SVGD, using the authors code\footnote{https://github.com/DartML/Stein-Variational-Gradient-Descent}, for 3000 epochs with the learning rate chosen by the grid search in the range \{0.00005,0.0001,0.0005,0.001,0.005,0.001\}. The optimal learning rate was 0.001 for Boston,	0.005 for Concrete, 0.005 for	Power, 0.001 for Yacht, 0.005 for Kin8nm, and 0.005 for MSD. We trained the PBP, using the authors' code\footnote{https://github.com/HIPS/Probabilistic-Backpropagation} for 40 epochs as recommended in~\cite{HernandezLobato15} (the learning rate need not be fine tuned because it is absent in the PBP as such). For the ML$_{\rm drop}$, we used the same hyperparameters as for the ML, except for the dropout rate, which was set to $0.4$ and used for both training and prediction. For prediction, we used 50 samples as suggested in~\cite{KendallGal17}.


\begin{table}
\resizebox{\textwidth}{!}{%
\begin{tabular}{lcccc}
{} & {\bf Boston} & {} & {} & {}\\
\toprule
{} & {\bf Optimizer} & \begin{tabular}{@{}c} {\bf Learning} \\ {\bf rate} \end{tabular} & {\bf Dropout} & \begin{tabular}{@{}c} {\bf Number}\\ {\bf of epochs} \end{tabular}  \\
\midrule
{\bf ML}                       &    Adam &         0.0001  & 0.4 & 700   \\
{\bf ML$_{\rm Bayes}$}         &    Adam &         0.0001  & 0.4 & 700   \\
{\bf DPD}                      &    Nesterov &     0.00002 & 0.4 & 5000 \\
{\bf GCP}                      &    Adam &         0.0001  & 0.3 & 700  \\
\bottomrule
\end{tabular}
\
\begin{tabular}{cccc}
 {\bf Concrete} & {} & {} & {}\\
\toprule
 {\bf Optimizer} & \begin{tabular}{@{}c} {\bf Learning} \\ {\bf rate} \end{tabular} & {\bf Dropout} & \begin{tabular}{@{}c} {\bf Number}\\ {\bf of epochs} \end{tabular}  \\
\midrule
    Adam &     0.0001  & 0.1 & 800  \\
    Adam &     0.0001  & 0.1 & 800  \\
    Nesterov & 0.00001 & 0.1 & 5000 \\
    Adam & 0.0001 & 0.1 & 1000 \\
\bottomrule
\end{tabular}
\
\begin{tabular}{cccc}
 {\bf Power} & {} & {} & {}\\
\toprule
 {\bf Optimizer} & \begin{tabular}{@{}c} {\bf Learning} \\ {\bf rate} \end{tabular} & {\bf Dropout} & \begin{tabular}{@{}c} {\bf Number}\\ {\bf of epochs} \end{tabular}  \\
\midrule
     Adam & 0.00005 & 0 & 150 \\
     Adam & 0.00005 & 0.2 & 150 \\     
     Adam & 0.0001  & 0 & 400  \\
     Adam & 0.00005 & 0 & 150 \\
\bottomrule
\end{tabular}
}
\caption{Optimizers and hyperparameters for the ML, DPD, and GCP methods for the Boston, Concrete, and Power data sets.}\label{tableHyperparametersBostonConcrete}
\end{table}
\begin{table}
\resizebox{\textwidth}{!}{%
\begin{tabular}{lcccc}
{} & {\bf Yacht} & {} & {} & {}\\
\toprule
{} & {\bf Optimizer} & \begin{tabular}{@{}c} {\bf Learning} \\ {\bf rate} \end{tabular} & {\bf Dropout} & \begin{tabular}{@{}c} {\bf Number}\\ {\bf of epochs} \end{tabular}  \\
\midrule
{\bf ML}                       &    Adam &         0.0001 & 0.1 & 2000   \\
{\bf ML$_{\rm Bayes}$}         &    Adam &         0.0001 & 0.1 & 4000   \\
{\bf DPD}                      &    Adam     &         0.0002 & 0.1 & 2500 \\
{\bf GCP}                      &    RmsProp &         0.001 & 0.1 & 1000  \\
\bottomrule
\end{tabular}
\
\begin{tabular}{cccc}
 {\bf Kin8nm} & {} & {} & {}\\
\toprule
 {\bf Optimizer} & \begin{tabular}{@{}c} {\bf Learning} \\ {\bf rate} \end{tabular} & {\bf Dropout} & \begin{tabular}{@{}c} {\bf Number}\\ {\bf of epochs} \end{tabular}  \\
\midrule
    Adam & 0.0002 & 0 & 200  \\
    Adam & 0.0002 & 0.1 & 200  \\    
    Adam & 0.0001 & 0 & 400 \\
    Nesterov & 0.0007 & 0 & 250 \\
\bottomrule
\end{tabular}
\
\begin{tabular}{cccc}
 {\bf MSD} & {} & {} & {}\\
\toprule
 {\bf Optimizer} & \begin{tabular}{@{}c} {\bf Learning} \\ {\bf rate} \end{tabular} & {\bf Dropout} & \begin{tabular}{@{}c} {\bf Number}\\ {\bf of epochs} \end{tabular}  \\
\midrule
     Adam & 0.005 & 0.1 & 150  \\
     Adam & 0.005 & 0.1 & 150  \\     
     Adam & 0.005 & 0.1 & 100 \\
     Adam & 0.001 & 0.1 & 200 \\
\bottomrule
\end{tabular}
}
\caption{Optimizers and hyperparameters for the ML, DPD, and GCP methods for the Yacht, Kin8nm, and MSD data sets.}\label{tableHyperparametersYacht}
\end{table}

\bigskip

{\bf Acknowledgements}. Both authors would like to thank the DFG project SFB 910. The
research of the first author was also supported by the DFG Heisenberg Programme and by the ``RUDN University Program 5-100''. The authors are grateful to anonymous referees for their comments and suggestions, which significantly clarified and improved our presentation.

\end{document}